\documentclass[12pt,a4paper,reqno]{amsart}
\usepackage{amsfonts,amsthm,amsmath,amssymb,
amsbsy,euscript,wasysym}

\newtheorem{theor}{Theorem}
\newtheorem{criterion}[theor]{Criterion}
\theoremstyle{definition}

\newtheorem{state}[theor]{Proposition}
\newtheorem{proposition}[theor]{Proposition}
\newtheorem{lemma}[theor]{Lemma}
\newtheorem{cor}[theor]{Corollary}
\newtheorem{define}{Definition}

\newtheorem{OpenProblem}{Open problem}
\newtheorem{example}{Example}[section]
\newtheorem{counterexample}[example]{Counterexample}
\newtheorem*{convent}{Convention}
\newtheorem*{notation}{Notation}
\newtheorem{model}{Model}
\theoremstyle{remark}
\newtheorem{rem}{Remark}[section]

\newcommand{\pinner}{\mathbin{\mathchoice
   {\hbox{\vrule width0.6em depth0pt height0.4pt
   \vrule width0.4pt depth0pt height0.8ex}}
   {\hbox{\vrule width0.6em depth0pt height0.4pt
   \vrule width0.4pt depth0pt height0.8ex}}
   {\hbox{\kern0.14em
   \vrule width0.48em depth0pt height0.4pt
   \vrule width0.4pt depth0pt height0.6ex\kern0.14em}}
   {\hbox{\kern0.1em
   \vrule width0.39em depth0pt height0.4pt
   \vrule width0.4pt depth0pt height0.5ex\kern0.1em}}}}

\newcommand{\cEv}{\partial}

\newcommand{\BBR}{\mathbb{R}}
\newcommand{\BBF}{\mathbb{F}}
\newcommand{\BBS}{\mathbb{S}}\newcommand{\BBT}{\mathbb{T}}
\newcommand{\BBZ}{\mathbb{Z}}

\newcommand{\cA}{{{\EuScript A}}}

\newcommand{\bcP}{{\boldsymbol{\mathcal{P}}}}
\newcommand{\bcQ}{{\boldsymbol{\mathcal{Q}}}}
\newcommand{\bcR}{{\boldsymbol{\mathcal{R}}}}

\newcommand{\cC}{\mathcal{C}}

\newcommand{\cF}{\mathcal{F}}
\newcommand{\cH}{\mathcal{H}}
\newcommand{\cI}{\mathcal{I}}

\newcommand{\cP}{\mathcal{P}}\newcommand{\cQ}{\mathcal{Q}}

\newcommand{\cX}{{\EuScript X}}    
\newcommand{\boldb}{{\boldsymbol{b}}}
\newcommand{\bolds}{{\boldsymbol{s}}}

\newcommand{\ba}{{\boldsymbol{a}}}
\newcommand{\bb}{{\boldsymbol{b}}}

\newcommand{\bp}{{\boldsymbol{p}}}
\newcommand{\bq}{{\boldsymbol{q}}}

\newcommand{\bs}{{\boldsymbol{s}}}

\newcommand{\bx}{{\boldsymbol{x}}}
\newcommand{\bby}{{\boldsymbol{y}}}
\newcommand{\bz}{{\boldsymbol{z}}}

\newcommand{\bQ}{{\boldsymbol{Q}}}

\newcommand{\bpi}{{\boldsymbol{\pi}}}

\newcommand{\binfty}{\pmb{\infty}}
\newcommand{\BOne}{{\boldsymbol{1}}}
\newcommand{\BTwo}{{\boldsymbol{2}}}

\newcommand{\bun}{\mathbf{1}}

\newcommand{\got}{\mathfrak{t}}
\newcommand{\gotht}{\mathfrak{t}}

\newcommand{\veps}{\varepsilon}
\newcommand{\vph}{\varphi}
\newcommand{\dd}{\partial}
\newcommand{\Id}{{\mathrm d}}

\newcommand{\vx}{{\vec{\mathrm{x}}}}
\newcommand{\bvx}{{\vec{\mathbf{x}}}}

\newcommand{\piNC}{\pi_{\text{\textsc{nc}}}}
\newcommand{\MnnC}{M^n_{\text{\textsc{nc}}}}
\newcommand{\cAZO}{\cA^{(0|1)}}
\newcommand{\bpiNC}{\bpi_{\text{\textsc{nc}}}}
\newcommand{\bpiNCZO}{\bpi^{(0|1)}_{\text{\textsc{nc}}}}

\renewcommand{\le}{\leqslant}
\renewcommand{\ge}{\geqslant}

\DeclareMathOperator{\Span}{span}
\DeclareMathOperator{\spanOp}{span}

\DeclareMathOperator{\img}{im}

\DeclareMathOperator{\Mat}{Mat}

\DeclareMathOperator{\tr}{tr}
\DeclareMathOperator{\Map}{Map}
\DeclareMathOperator{\grad}{grad}

\newcommand{\GH}[1]{|{#1}|}
\newcommand{\schouten}[1]{\lshad {#1} \rshad}

\newcommand{\Free}{\text{\textsf{Free}}\,}

\DeclareMathOperator{\dvol}{d
vol}
\DeclareMathOperator{\supp}{supp}
\DeclareMathOperator{\jet}{jet}

\newcommand{\lshad}{[\![}
\newcommand{\rshad}{]\!]}
\newcommand{\ov}{\overline}
\newcommand{\nC}{{\text{\textsc{nc}}}}
\newcommand{\KdV}{{\text{KdV}}}

\newcommand{\eqdef}{\mathrel{\stackrel{\text{\textup{def}}}{=}}}

\newcommand{\ib}[3]{ \{\!\{ {#1},{#2} \}\!\}_{{#3}} }

\newcommand{\Mars}{\mars}
\newcommand{\Venus}{\venus}

\newcommand{\by}[1]{\textit{{#1}}}
\newcommand{\jour}[1]{\textit{{#1}}}
\newcommand{\vol}[1]{\textbf{{#1}}}
\newcommand{\book}[1]{\textrm{{#1}}}

\title[
Variational formal noncommutative symplectic (super)\/geometry%
]{The calculus of multivectors\\[3pt] on noncommutative jet spaces}

\date{18 December 2017}

\author[Arthemy Kiselev]{Arthemy V. Kiselev}
\thanks{\textit{Address}: Johann Ber\-nou\-lli Institute for Mathematics and Computer Science, University of Groningen,
P.O.~Box 407, 9700~AK Groningen, The Netherlands.
\quad\textit{E-mail}: \texttt{A.V.Kiselev\symbol{"40}rug.nl}
}



\begin{document}
\begin{abstract}
The Leibniz rule for derivations is 
invariant under cyclic permutations of 
co\/-\/multiples within the arguments of derivations. 
We explore the implications of this 
prin\-ci\-ple: in effect, 
we construct a class of noncommutative 
bundles in which the sheaves of algebras of walks along a tesselated affine manifold form the base, whereas the fibres are free associative algebras or, at a later stage, such algebras quotients over the linear relation of equivalence under 
cyclic shifts. 
The calculus of 
variations is 
developed on the infinite jet spaces over such noncommutative bundles.


In the frames of such 
field\/-\/theoretic extension of the Kon\-tse\-vich for\-mal non\-com\-mu\-ta\-ti\-ve sym\-p\-lec\-tic (super)\/geometry, we prove the main properties of the Batalin\/--\/Vilkovisky Laplacian and 
Schouten bracket.
We show as by\/-\/product that the structures which arise in the classical variational Poisson geometry of infinite\/-\/di\-men\-sio\-nal integrable systems 
do actually not refer to 
the graded com\-mu\-ta\-ti\-vi\-ty assumption.
\end{abstract}
\maketitle

\enlargethispage{\baselineskip}

\tableofcontents

\subsection*{Introduction}
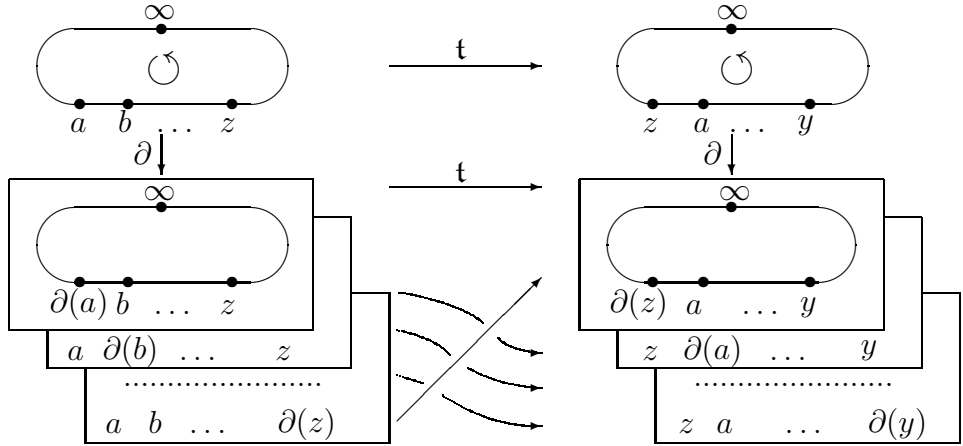
\begin{figure}[htb]
\begin{center}
\unitlength=1mm
\linethickness{0.4pt}
\begin{picture}(125.00,56.33)
\put(0.00,15.00){\line(1,0){40.00}}
\put(40.00,15.00){\line(0,1){20.00}}
\put(40.00,35.00){\line(-1,0){40.00}}
\put(0.00,35.00){\line(0,-1){20.00}}
\put(5.00,15.00){\line(0,-1){5.00}}
\put(5.00,10.00){\line(1,0){40.00}}
\put(45.00,10.00){\line(0,1){20.00}}
\put(45.00,30.00){\line(-1,0){5.00}}
\put(10.00,10.00){\line(0,-1){10.00}}
\put(10.00,0.00){\line(1,0){40.00}}
\put(50.00,0.00){\line(0,1){20.00}}
\put(50.00,20.00){\line(-1,0){5.00}}
\put(20.17,26.33){\oval(33.00,10.00)[]}
\put(20.17,50.00){\oval(33.00,10.00)[]}
\put(20.00,55.00){\circle*{1.33}}
\put(9.33,45.00){\circle*{1.33}}
\put(15.67,45.00){\circle*{1.33}}
\put(29.33,45.00){\circle*{1.33}}
\put(9.33,21.33){\circle*{1.33}}
\put(15.67,21.33){\circle*{1.33}}
\put(29.33,21.33){\circle*{1.33}}
\put(20.00,31.33){\circle*{1.33}}
\put(95.00,31.33){\circle*{1.33}}
\put(20.00,41.00){\vector(0,-1){5.33}}
\put(50.00,50.00){\vector(1,0){20.00}}
\put(50.00,34.00){\vector(1,0){20.00}}
\put(51.00,3.00){\vector(1,1){19.00}}
\bezier{48}(51.33,20.00)(59.00,18.67)(62.33,16.00)
\bezier{32}(51.00,15.00)(56.00,14.00)(58.00,11.67)
\bezier{12}(51.00,9.00)(53.33,8.33)(54.33,8.00)
\bezier{56}(57.00,6.33)(63.67,2.00)(70.00,2.33)
\bezier{40}(60.33,10.00)(64.33,7.67)(70.00,7.33)
\bezier{24}(64.33,14.00)(65.67,12.33)(70.00,12.00)
\put(69.00,12.00){\vector(1,0){1.00}}
\put(68.67,7.33){\vector(1,0){1.33}}
\put(68.67,2.33){\vector(1,0){1.33}}
\put(96.50,50.00){\oval(33.00,10.00)[]}
\put(75.00,35.00){\line(1,0){40.00}}
\put(115.00,35.00){\line(0,-1){20.00}}
\put(115.00,15.00){\line(-1,0){40.00}}
\put(75.00,15.00){\line(0,1){20.00}}
\put(80.00,15.00){\line(0,-1){5.00}}
\put(80.00,10.00){\line(1,0){40.00}}
\put(120.00,10.00){\line(0,1){20.00}}
\put(120.00,30.00){\line(-1,0){5.00}}
\put(85.00,10.00){\line(0,-1){10.00}}
\put(85.00,0.00){\line(1,0){40.00}}
\put(125.00,0.00){\line(0,1){20.00}}
\put(125.00,20.00){\line(-1,0){5.00}}
\put(95.00,26.33){\oval(32.67,10.00)[]}
\put(84.67,45.00){\circle*{1.33}}
\put(91.33,45.00){\circle*{1.33}}
\put(105.33,45.00){\circle*{1.33}}
\put(84.67,21.33){\circle*{1.33}}
\put(91.33,21.33){\circle*{1.33}}
\put(105.33,21.33){\circle*{1.33}}
\put(95.00,41.00){\vector(0,-1){5.33}}
\put(95.00,55.00){\circle*{1.33}}
\put(17.67,56.33){\makebox(0,0)[lb]{$\infty$}}
\put(93.00,56.33){\makebox(0,0)[lb]{$\infty$}}
\put(58.67,51.00){\makebox(0,0)[lb]{$\mathfrak{t}$}}
\put(58.67,34.67){\makebox(0,0)[lb]{$\mathfrak{t}$}}
\put(18.00,47.50){\makebox(0,0)[lb]{{\Large$\circlearrowleft$}}}
\put(93.33,47.50){\makebox(0,0)[lb]{{\Large$\circlearrowleft$}}}
\put(8.00,41.33){\makebox(0,0)[lb]{$a$}}
\put(14.33,41.33){\makebox(0,0)[lb]{$b$}}
\put(19.33,41.33){\makebox(0,0)[lb]{$\dots$}}
\put(27.67,41.33){\makebox(0,0)[lb]{$z$}}
\put(83.33,41.33){\makebox(0,0)[lb]{$z$}}
\put(90.00,41.33){\makebox(0,0)[lb]{$a$}}
\put(94.67,41.33){\makebox(0,0)[lb]{$\dots$}}
\put(103.67,40.83){\makebox(0,0)[lb]{$y$}}
\put(16.33,37.00){\makebox(0,0)[lb]{$\dd$}}
\put(91.33,37.00){\makebox(0,0)[lb]{$\dd$}}
\put(17.67,32.33){\makebox(0,0)[lb]{$\infty$}}
\put(93.33,32.33){\makebox(0,0)[lb]{$\infty$}}
\put(5.33,16.3){\makebox(0,0)[lb]{$\dd(a)$}}
\put(14.00,17.17){\makebox(0,0)[lb]{$b$}}
\put(19.00,17.17){\makebox(0,0)[lb]{$\dots$}}
\put(27.67,17.17){\makebox(0,0)[lb]{$z$}}
\put(79.00,16.3){\makebox(0,0)[lb]{$\dd(z)$}}
\put(89.00,17.17){\makebox(0,0)[lb]{$a$}}
\put(96.33,17.17){\makebox(0,0)[lb]{$\dots$}}
\put(104.00,16.67){\makebox(0,0)[lb]{$y$}}
\put(7.67,11.00){\makebox(0,0)[lb]{$a$}}
\put(12.33,10.13){\makebox(0,0)[lb]{$\dd(b)$}}
\put(22.33,11.00){\makebox(0,0)[lb]{$\dots$}}
\put(35.00,11.00){\makebox(0,0)[lb]{$z$}}
\put(83.33,11.00){\makebox(0,0)[lb]{$z$}}
\put(88.67,10.13){\makebox(0,0)[lb]{$\dd(a)$}}
\put(100.00,11.00){\makebox(0,0)[lb]{$\dots$}}
\put(112.00,11.00){\makebox(0,0)[lb]{$y$}}
\put(90.00,7.67){\makebox(0,0)[lb]{$.......................$}}
\put(15.00,7.67){\makebox(0,0)[lb]{$.......................$}}
\put(12.67,1.33){\makebox(0,0)[lb]{$a$}}
\put(18.33,1.33){\makebox(0,0)[lb]{$b$}}
\put(23.67,1.33){\makebox(0,0)[lb]{$\dots$}}
\put(35.33,0.46){\makebox(0,0)[lb]{$\dd(z)$}}
\put(88.00,1.33){\makebox(0,0)[lb]{$z$}}
\put(93.00,1.33){\makebox(0,0)[lb]{$a$}}
\put(102.00,1.33){\makebox(0,0)[lb]{$\dots$}}
\put(113.00,0.46){\makebox(0,0)[lb]{$\dd(y)$}}
\end{picture}
\caption{The cyclic\/-\/shift invariance of derivations.}\label{FigIntro}
\end{center}
\end{figure}

Let $\BBF$ be a free 
algebra over $\Bbbk\mathrel{{:}{=}}\BBR$
and suppose 
$a_1$,\ $\ldots$,\ $a_k\in
\BBF$. Denote by~$\circ$ the associative multiplication in~$\BBF$
and by~$\got$ the counterclockwise cyclic shift of co\/-\/multiples in the product~$a_1\circ\ldots\circ a_k$,
\[
\got\,(a_1\circ\ldots\circ a_{k-1}\circ a_k) \eqdef
a_k\circ a_1\circ\ldots a_{k-1}.
\]
For the sake of definition, now assume that a given derivation $\dd\colon\BBF\to\BBF$ is such that its values at $a_1$,\ $\ldots$,\ $a_k$ do not leave that set.
By the Leibniz rule, the derivation is cyclic\/-\/shift invariant:
\begin{equation}\label{EqDerCycle}
\dd\bigl(\got\,(a_1\circ\ldots\circ a_k)\bigr) =
\got\,\bigl(\dd(a_1\circ\ldots\circ a_k)\bigr).
\end{equation}
Indeed, both sides of the above equality are given by the sum
\[
\dd(a_k)\circ a_1\circ\ldots\circ a_{k-1} + a_k\circ \dd(a_1)\circ\ldots\circ a_k +
a_k\circ a_1\circ\ldots\circ\dd(a_{k-1}),
\]
up to a 
sequential order in which these $k$~summands follow each other (see Fig.~\ref{FigIntro}).
This observation is generalised in an obvious way to the case where the elements of algebra~$\BBF$ are graded by some Abelian group, each element~$a_1$,\ $\ldots$,\ $a_k$~is homogeneous with respect to the grading, and $\dd\colon\BBF\to\BBF$ is a graded derivation (i.e.\ not necessarily preserving the set~$\{a_1$,\ $\ldots$,\ $a_k\}$ at hand).

How much (graded-) commutativity is really needed to make the calculus of variations in the  Lagrangian and Hamiltonian formalisms work, thus allowing for the Batalin\/--\/Vilkovisky technique for quantisation of gauge systems~--- and creating a 
cohomological approach to the complete integrability of infinite\/-\/dimensional KdV-\/type systems\,?\footnote{We refer to~\cite{BV1981,BV1983,CattaneoFelderCMP2000,GitmanTyutin,HenneauxTeitelboim,
SchwarzCMP1993,WittenAntibracket,ZinnJustin1975} or~\cite{gvbv} and 
to~\cite{DeSoleKacCMP2012,DeSoleKacJap2013,DubrZhang2001,LiuZhang2009AdvM,
LiuZhangCMP2013,Magri1979,Olver1993,OlverSokolovCMP1998}
respectively (see also~\cite{Prague2011,Galli10} in both contexts).}

We 
claim that it is not the restrictive assumption of commutativity that shows through \emph{arbitrary} permutations~--- but it is the linear equivalence $a\sim\got(a)$ of words~$a$, written in a given 
alphabet, with respect to the \emph{cyclic} permutations~$\got$ that is sufficient for the structures of the calculus of iterated variations to be well defined. Introduced in this cyclic\/-\/invariant setup, the Batalin\/--\/Vilkovisky Laplacian~$\Delta$ and variational Schouten bracket~$\lshad\,,\,\rshad$ are 
proven to satisfy the main identities such as the cocycle condition~$\Delta^2=0$, see 
(\ref{EqDeviationDerivationIntro}--\ref{EqDeltaSquareIntro}) below.
Both the definitions and assertions are then literally valid in the sub\/-\/class of graded\/-\/commutative geometries; the reason why 
is that the latter can be obtained from the former by using the postulated commutativity reduction at the end of the day when the proof is over.

The idea to establish 
the formal noncommutative symplectic geometry on the cyclic invariance, generalising the geometry of commutative symplectic manifolds, was introduced by Kontsevich in~\cite{KontsevichCyclic}, cf.~\cite{Looijenga} and references therein.
The quotient spaces of cyclic words were employed as target sets for maps from usual manifolds in~\cite{OlverSokolovCMP1998} by Olver and Sokolov (cf.\ Model~\ref{ModelMatr} on p.~\pageref{ModelMatr} below); 
several integrable equations of KdV-\/type were 
recovered in such noncommutative set\/-\/up.\footnote{Noncommutative extensions of classical infinite\/-\/dimensional systems can acquire new components that are invisible in the commutative world: e.g., there appear --\,often, through 
nonlocalities\,-- the terms that contain the commutants $a_i\circ a_j-a_j\circ a_i$.}
Variations arising in the 
Poisson or Schouten brackets for integral functionals, their calculus was then pursued along the lines of~\cite{Olver1993}. 
The paper~\cite{OlverSokolovCMP1998} initiated a 
classification and study of evolutionary ODE and PDE systems on associative algebras, 
which required the calculation of standard geometric structures for such models in jet spaces (e.g., see~\cite{Treves2007
} in this context).

In this paper we futher that approach to noncommutative jet spaces.\footnote{We note that the positive 
differential order calculus on infinite jet spaces lies far beyond the bare tensor calculus on usual commutative manifolds; for instance, compare~\cite{Vaintrob} with~\cite{Galli10} 
or contrast~\cite{AKZS} vs~\cite{Norway} 
and~\cite{SchwarzCMP1993
} vs~\cite{gvbv}.}
Continuing the line of reasoning from~\cite{gvbv,sqs13,prg15} where the intrinsic regularisation of Batalin\/--\/Vilkovisky formalism is revealed, we 
verify the main identities for~$\Delta$ and~$\lshad\,,\,\rshad$ in the variational noncommutative set\/-\/up of (homogeneous) local functionals~$F$,\ $G$,\ $H$:
\begin{subequations}\label{EqAllIntro}
\begin{align}
\Delta(F\times G) &= \Delta F\times G + (-)^{|F|}\lshad F,G\rshad
+(-)^{|F|}F\times\Delta G,\label{EqDeviationDerivationIntro}
\\
\lshad F,G\times H\rshad &= \lshad F,G\rshad\times H 
+ (-)^{(|F|-1)\cdot|G|} G\times\lshad F,H\rshad,\label{EqSchoutenOnProductIntro}
\\
\Delta\bigl(\lshad F,G\rshad\bigr) &= \lshad\Delta F,G\rshad 
+(-)^{|F|-1}\lshad F,\Delta G\rshad,
\label{EqZimesIntro}
\\
\text{Jacobi}\bigl(\lshad\,,\,\rshad\bigr)&=0\qquad \Longleftrightarrow\qquad 
\Delta^2 = 0.\label{EqDeltaSquareIntro}
\end{align}
\end{subequations}
It is quite paradoxical that for a long time, these identities were proclaimed to be valid just formally~\cite{GomisParisSamuel,HenneauxTeitelboim}; 
for it was believed that the Batalin\/--\/Vilkovisky technique would necessarily contain some divergencies or ``infinite constants'', 
whereas their manual regularisation appealed to surreal principles like ``$\boldsymbol{\delta}(0)\mathrel{{:}{=}}0$'' for the Dirac $\boldsymbol{\delta}$-\/function (see~\cite{gvbv} and references therein for discussion on the history of the problem).

The notion of associative algebra structures itself has deserved much attention in the 
mathematical physics literature, e.g., in relation to the 
Yang\/--\/Baxter equation. 
Such structures arise naturally in the topological context;
the calculus of cyclic words 
serves the alphabet of homotopy group generators.
Likewise, the 
multiplication in homology gives rise to the Gromov\/--\/Witten potential solving the WDVV~equations, see~\cite{DubrZhang2001} and~\cite{DVVWDVV,WittenWDVV}, cf.~\cite{ManinFmanifolds}.
Another construction, which will be discussed in Remark~\ref{RemV1V2} on p.~\pageref{RemV1V2} below, stems from the calculation of matrix integrals in the Batalin\/--\/Vilkovisky framework~\cite{Bar2007,BarCRM2010}. 
Furthermore, associative but not necessarily commutative $\star$-\/products are obtained --\,on finite\/-\/dimensional affine manifolds\,-- by using the deformation quantisation procedure~\cite{KontsevichFormality}, cf.~\cite{cpp,dq15} and Model~\ref{ModelStar} on p.~\pageref{ModelStar}.
Now we study the extent to which the differential calculus can be developed on the basis of 
associative algebra structures as input data.%
\footnote{An alternative approach to 
noncommutativity suggests that 
manifolds --\,and derivative objects such as the fibre bundles\,-- are determined as the spectra of associative noncommutative algebras. 
Provided that the algebras are `smooth', they 
are viewed as the algebras of smooth functions on the objects which they determine. Nowadays, noncommutative geometry \`a\ la Connes~\cite{ConnesBook} is a well\/-\/established domain. 
However, we keep the framework closer to the needs which one encounters in
a 
class of path-{} and loop\/-\/based QFT 
models
~\cite{LoopReview,LoopQG,Wilson1975}). Let us therefore 
study the language of closed strings of symbols~-- written around the circles and encoding paths in the granulated space~$M^n$
(see Model~\ref{ModelBV} 
on p.~\pageref{ModelBV} 
below).}%
\\
\centerline{\rule{1in}{0.7pt}}

\noindent%
This paper consists of three parts. 
In Ch.~\ref{SecStatic} we introduce the static set\/-\/up of noncommutative infinite jet (super-)\/spaces. Based on the algorithmic construction of parity\/-\/odd 
Laplacian~$\Delta$ and variational Schouten bracket~$\lshad\,,\,\rshad$, the calculus of iterated variations of local functionals --\,i.e., kinematics\,-- is developed in Ch.~\ref{SecKinematics}. 
Such BV\/-\/geometry of local functionals is then 
contrasted in Ch.~\ref{SecDynamics} with the noncommutative 
Poisson formalism, that is, the dynamics determined by variational multi\/-\/vectors.

The text is structured as follows.
The commutative but not associative algebra~$\cA$ of cyclic words written in the alphabet~$\langle a^i\rangle$ of a free associative algebra is introduced in~\S\ref{SecAlgebra}. 
The generators~$a^i$ themselves are viewed in~\S\ref{SecAs} as 
words written in the alphabet $\langle\vx^{\pm1}_i\rangle$ of edges in the adjacency graph for a cell\/-\/complex tiling of the substrate manifold~$M^n$, which is introduced in~\S\ref{SecWhy}. 
The alphabets $\langle\vx^{\pm1}_i\rangle$ and~$\langle a^i\rangle$ provide
the respective noncommutative analogues of base and fibre in the bundle~$\piNC$:
the base is the sheaf of [unital extensions of] free associative algebras generated by $\langle\vx^{\pm1}_i\rangle$ for a crystal tiling of~$M^n$, whereas the fibres of~$\piNC$ are
[the unital extension of] the algebra~$\cA$ of cyclic words written in the alphabet~$\langle a^i\rangle$ (see the figure on p.~\pageref{FigBundlePi}).
The jet space $J^\infty(\piNC 
)$ of sections is built in~\S\ref{SecJets};
various elements of the jet\/-\/space language are then recovered.
In particular, as soon as the notion of variational (co)vec\-tors is available, 
we show 
why the Substitution Principle works for (non)\/commutative identities in total derivatives. 

The second part begins with the definition of noncommutative analogue for the variational cotangent bundle over the infinite jet space $J^\infty(\piNC 
)$, see~\S\ref{SecCotangent}. The sections target algebra alphabet~$\langle a^i\rangle$ is doubled by using the canonical pairs~$\langle a^i,a^\dagger_i\rangle$; 
sign convention~\eqref{EqNormalizeShifts} for the two ordered couplings of the virtual 
variations~$\delta\ba$ and~$\delta\ba^\dagger$ ensures 
the matching of signs in all the structures that are defined in what follows.
In the meantime (see~\S\ref{SecNeighbour}), 
the $\BBZ_2$-\/parity reversion 
$\Pi\colon a^\dagger_i\rightleftarrows b_i$ acts on the dual symbols~$\ba^\dagger$, producing the parity\/-\/odd slots~$\boldb$.
Now, the geometric approach of~\cite{gvbv} to iterated variations works in the noncommutative set\/-\/up of evaluation \emph{maps} $\ba=\bs(\bx,\bvx^{\pm1})$ and \emph{antimaps} $\ba^\dagger=\bs^\dagger(\bx,\bvx^{\pm1})$ using the sheaf over~$M^n$ (see Fig.~\ref{FigSetup} on p.~\pageref{FigSetup}).
Therefore, while giving the operational definition of BV-\/Laplacian~$\Delta$ in~\S\ref{SecLaplacian}, we focus on the unlock\/-\/and\/-\/join reconfigurations of cyclic words. The variational Schouten bracket~$\lshad\,,\,\rshad$ is a derivative structure, that is, it is determined by the parity\/-\/odd operator~$\Delta$ via its action on products, 
as in~\eqref{EqDeviationDerivationIntro} above.\footnote{In geometric terms, the bracket~$\lshad\,,\,\rshad$ of cyclic word\/-\/valued functionals is encoded by the standard topological pair of pants~$\BBS^1\times\BBS^1\to\BBS^1$ that links the cycles. 
In fact, this topological procedure 
also underlies 
each of the following 
structures and operations 
in the differential calculus under study: 
\begin{itemize}
\item multiplication~$\times$ of cyclic words and word\/-\/valued function(al)s, 
\item termwise action of derivations (e.g., in~\eqref{EqTermwiseActDeriv}), including 
\item the commutation of vector fields,~--- and also
\item evaluation of multi\/-\/vectors at the tuples of covectors (see~\eqref{EvalkVector}): in particular,
\item the Poisson bracket of Hamiltonian functionals.
\end{itemize}
Indeed, all of the above 
amounts to the 
detach\/-\/and\/-\/join picture 
$\BBS^1\times\BBS^1\to\BBS^1$.}
Then we confirm that the variational Schouten bracket~$\lshad\,,\,\rshad$ is shifted\/-\/graded skew\/-\/symmetric and satisfies the Jacobi identity.
The two structures~$\Delta$ and $\lshad\,,\,\rshad$ endow the ring of local functionals with the structure of differential 
graded Lie algebra.

The third part of this text narrates on the noncommutative variational Poisson formalism. The notion of noncommutative variational multi\/-\/vectors is introduced in~\S\ref{SecMultVect}. We recall that not every grading\/-\/homogeneous integral functional over the infinite jet superspace 
$J^\infty(\bpiNCZO 
)$, canonically extended in Ch.~\ref{SecKinematics}, would be a well defined variational multi\/-\/vector containing the respective number of parity\/-\/odd slots~$\boldb$. Remark~\ref{RemABAB} on p.~\pageref{RemABAB} is a key to that concept. Specifically, by viewing the variational multi\/-\/vectors as maps that take the respective tuples of --\,possibly, exact\,-- variational covectors to the top\/-\/degree 
horizonal cohomology space of cyclic word\/-\/valued integral functionals, we analyse in~\S\ref{SecDB} the 
geometry of iterated variations that arise in the derived brackets encoding such maps. 
We 
discover that the calculus of noncommutative variational multivectors is the paradigm of steps and stops. Finally, we arrive at the definition of Poisson brackets. In~\S\ref{SecPBr} we study the geometry of differential forms that stands behind the criterion under which the variational noncommutative bi\/-\/vectors are Poisson, i.e.\ endow the space of noncommutative Hamiltonians with the variational Poisson brackets. (In particular, the Helmholtz lemma is proved in the setting of cyclic words.)


\newpage
\section{
The nature of associative symbols}\label{SecStatic}
\subsection{
The algebra $\cA$ of cyclic words}\label{SecAlgebra}
In this section we introduce the main object to consider in the 
future reasoning. Namely, by starting with a non\/-\/commutative free associative algebra, we define the commutative but not associative unital algebra $\cA$ 
of cyclic words written in the free algebra's alphabet. Note that for the sake of 
clarity,
neither of these two
algebras is graded; however, in what follows we shall extend
the alphabet by using 
symbols of $\BBZ_2$-\/valued parity.

Throughout this 
text, the ground field $\Bbbk$ is the field $\BBR$ of real numbers.\\
\centerline{\rule{3cm}{0.4pt}.}

\noindent%
Consider the free associative algebra $\Free(a^1,\,\dots\,,a^m)$ with $m$ generators $a^1,\,\dots\,,a^m$; let $m<\infty$ for definition. 
(One may presently think that the free algebra at hand is not necessarily unital.)
Denote by~$\circ$ the multiplication in that algebra. By definition, put
\begin{equation}\label{EqDefShift}
\gotht(a^i)=a^i,\qquad
\gotht\,(a^{i_1}\circ\,\dots\,\circ a^{i_{\lambda}}) \mathrel{{:}{=}}
a^{i_{\lambda}}\circ a^{i_1}\circ\,\dots\,\circ a^{i_{\lambda-1}},\quad \lambda>1;
\end{equation}
otherwise speaking, the operator~$\gotht$ is the counterclockwise cyclic permutation of symbols in a homogeneous word of length~$\lambda>0$.

Introduce the linear equivalence relation $\sim$ on $\Free(a^1,\,\dots\,,a^m)$ by setting%
\footnote{It is readily seen that
$a^{i_1}\circ\ldots\circ a^{i_{\lambda}}=\gotht^{\lambda-1}\,(\gotht\,(a^{i_1}\circ\ldots\circ a^{i_{\lambda}}))$ so that $a\sim a$
and $\gotht(a)\sim a$, whence the transitive relation $\sim$ is reflexive and symmetric indeed.}
$$a\sim\gotht(a),$$
where $a$ is a homogeneous word as in~\eqref{EqDefShift}, and then extending $\sim$ onto the algebra by linearity:
$a\sim a'$ and $b\sim b'$ implies $a+b\sim a'+b'$. For instance, one has that%
\footnote{We emphasize that the cyclic invariance itself does \emph{not} imply the commutativity: even though
$a^i\sim a^i$ and $a^i\circ a^j\sim a^j\circ a^i$ one has that $a^i\circ a^j\circ a^k\nsim a^i\circ a^k\circ a^j$
unless some of the indexes coinside.}
\[
a^1+a^2\circ a^3+a^1\circ a^2\circ a^3\sim a^1+a^3\circ a^2+ a^3\circ a^1\circ a^2.
\]
Notice also that 
$$a\sim \gotht(a)\sim\,\dots\,\sim\gotht^{\lambda(a)-1}(a)\sim
\frac1{\lambda(a)}\sum_{i=1}^{\lambda(a)}\gotht^{i-1}(a)$$
for any word $a$ of length $\lambda(a)>0$; by convention, a word of zero length is an element of the ground field $\Bbbk$,
see~\eqref{EqDefLength0} below.


We denote by $\cA$ the quotient $\Free(a^1,\,\ldots\,,a^m)/\sim$, that is, $\cA$ is the vector space of (formal sums of) cyclic
words such that each homogeneous component $a^{i_1}\circ\ldots\circ a^{i_{\lambda}}$ can be read starting from any
letter $a^{i_{\alpha}}$ for $1\le\alpha\le\lambda$. 
Let us denote by~$(a)\in\cA$
the equivalence class of an element $a\in\Free(a^1,\,\dots\,,a^m)$ under
cyclic permutations of symbols in all its homogeneous components (i.\,e.\ in all its ``words'' in proper sense).

Now we endow the vector space~$\cA$ of cyclic words with the algebra 
structure~$\times$. Consider the equivalence classes~$(a_1)$ and~$(a_2)$ of two homogeneous elements $a_1,a_2\in\Free(a^1,\,\dots\,,a^m)$ of postive lengths $\lambda(a_1)$ and~$\lambda(a_2)$, respectively. Let their product~be
\begin{equation}\label{EqDefMult}
(a_1)\times(a_2)\mathrel{\stackrel{\text{def}}{=}}
\frac1{ 
\lambda(a_1) 
\cdot 
\lambda(a_2) 
}
\Bigl(\sum_{i=1}^{\lambda(a_1)}\sum_{j=1}^{\lambda(a_2)}\gotht^{i-1}\,(a_1)\circ\gotht^{j-1}\,(a_2)\Bigr),
\end{equation}

\hangindent=-5cm\hangafter=4
\noindent%
where the equivalence class in the right\/-\/hand side is normalized 
in such a way that the definition correlates 
with the commutative set\/-\/up (should it be recovered postfactum);
now extend the product onto~$\cA$ by (bi-)\/linearity.
The definition of operation~$\times$ says that, each homogeneous string of symbols in the first co\/-\/multiple read, time
after time starting from every next letter, it is then pasted --~time after time in its turn~-- in between every two
consecutive letters occurring in each homogeneous string 
contained in the second co\/-\/multiple. Sure, this is the classical topological pair of pants $\BBS^1\times\BBS^1\to\BBS^1$
in which every symbol in either of the factors has the right to be read first, see the figure.
{\unitlength=1mm
\begin{picture}(0,0)(-26,-2.5)
\put(0,0){\oval(4,6)[t]}
\put(8.5,0){\oval(13,3)}
\put(-8.5,0){\oval(13,3)}
\put(0,20){\oval(15,3)}
\qbezier(15,0)(15,10)(10.75,10)
\qbezier(7.5,20)(7.5,10)(10.75,10)
\qbezier(-15,0)(-15,10)(-10.75,10)
\qbezier(-7.5,20)(-7.5,10)(-10.75,10)
\end{picture}%
}

\begin{state}
Multiplication~\eqref{EqDefMult} on $\cA$ is commutative. 
\end{state}

\begin{proof}
Notice that not only the necklace $(a_1)$ is unlocked at all possible multiplication signs~$\circ$ and joined to $(a_2)$
in between each pair of adjacent symbols in that word but, as one 
shifts the symbols in $(a_2)$ around the circle, 
exactly the same is done with respect to the insertion of $\gotht^{j-1}\,(a_2)$ into~$(a_1)$.
\end{proof}

However, it is readily seen that the symbols in homogeneous strings in $(a_1)$ and~$(a_2)$ always stay next to each other in the nested 
product $\bigl((a_1)\times(a_2)\bigr)\times(a_3)$, whereas they are separated by the symbols from $(a_3)$ in at least
one homogeneous term in $(a_1)\times\bigl((a_2)\times(a_3)\bigr)$, provided that the alphabet contains at least two different letters
and the 
length 
of the word~$a_3$ is greater than~one.
\footnote{Obviously, the associativity equation for $\times$ can be satisfied incidentally, for a 
special choice of the three co\/-\/multiples.}

\begin{proposition}
If $m\geqslant2$ so that the letters~$a^1$ and~$a^2$ are distinct in
the alphabet, multiplication~\eqref{EqDefMult} on~$\cA$ is not associative: 
\begin{equation}\label{EqNassoc}
\bigl((a_1)\times(a_2)\bigr)\times(a_3)\nsim(a_1)\times\bigl((a_2)\times(a_3)\bigr),
\end{equation}
see the figure below.\footnote{%
Let us recall that in Nature, not all processes are 
associative. For example, take a proton~$\mathsf{p}^+$, another proton, and a neutron~$\mathsf{n}^0$. Letting their strong interaction events be arranged using 
\[
\bigl((\cdot \times \cdot)\times \cdot\bigr)\colon 
\mathsf{p}^+ \sqcup \mathsf{p}^+ \sqcup \mathsf{n}^0 \longmapsto
(\mathsf{p}^+ \times \mathsf{p}^+)\times\mathsf{n}^0 = \mathsf{p}^+\sqcup\mathsf{p}^+\sqcup\mathsf{n}^0 \longmapsto \mathsf{p}^+ \sqcup \mathsf{D}_2^1,
\] 
one obtains the input objects intact after the first interaction event. But the arrangement \[\bigl(\cdot\times(\cdot\times\cdot)\bigr)\colon
\mathsf{p}^+ \sqcup \mathsf{p}^+ \sqcup \mathsf{n}^0 \longmapsto
\mathsf{p}^+\times(\mathsf{p}^+\times\mathsf{n}^0) = \mathsf{p}^+\times\mathsf{D}_2^1 = \text{\textsf{He}}{}_3^2\]
produces helium\/-\/3 via deuterium. This fusion process is not associative.
}
\begin{figure}[h]
{\unitlength=0.7mm
\centerline{%
\begin{picture}(33,25)(0,5.5)
\put(0,0){\oval(4,6)[t]}
\put(8.5,0){\oval(13,3)}
\put(-8.5,0){\oval(13,3)}
\qbezier(-15,0)(-15,10)(-10.75,10)
\qbezier(-10.75,10)(0,20)(0,30)
\put(10,30){\oval(20,4)}
\put(26.5,0){\oval(13,3)}
\qbezier(15,0)(15,10)(10.75,10)
\qbezier(10.75,10)(15,14)(15,15)
\qbezier(15,15)(20,10)(20,0)
\qbezier(33,0)(33,10)(28.75,15)
\qbezier(20,30)(26.5,15)(28.75,15)
\end{picture}\qquad $\neq$\qquad{\ }
\begin{picture}(33,0)(-30,5.5)
\put(0,0){\oval(4,6)[t]}
\put(8.5,0){\oval(13,3)}
\put(-8.5,0){\oval(13,3)}
\qbezier(15,0)(15,10)(10.75,10)
\qbezier(10.75,10)(0,20)(0,30)
\put(-10,30){\oval(20,4)}
\put(-26.5,0){\oval(13,3)}
\qbezier(-15,0)(-15,10)(-10.75,10)
\qbezier(-10.75,10)(-15,14)(-15,15)
\qbezier(-15,15)(-20,10)(-20,0)
\qbezier(-33,0)(-33,10)(-28.75,15)
\qbezier(-20,30)(-26.5,15)(-28.75,15)
\end{picture}
}}%
\end{figure}
\end{proposition}

\begin{counterexample}[``abba'']
Let $a_1\mathrel{{:}{=}}a^{\BOne}$, $a_2\mathrel{{:}{=}}a^{\BOne}$, 
and $a_3\mathrel{{:}{=}}a^{\BTwo}a^{\BTwo}$. 
Then $(a_1)\times(a_2)=(a^{\BOne}\circ a^{\BOne})$ so that these two copies of the letter $a^{\BOne}$ always stay next to each other
in any product of $(a_1)\times(a_2)$ with any other word. 
On the other hand (see Fig.~\ref{FigNassocCounterEx}),%
\begin{figure}[h]
\begin{center}
\unitlength=1mm
\linethickness{0.4pt}
\begin{picture}(110.00,17.5)
\put(15.00,10.00){\circle{10.00}}
\put(30.00,10.00){\circle{10.00}}
\put(45.00,10.00){\circle{10.00}}
\put(75.00,10.00){\circle{10.00}}
\put(90.00,10.00){\circle{10.00}}
\put(105.00,10.00){\circle{10.00}}
\bezier{40}(17.67,14.33)(22.33,12.33)(27.33,14.33)
\bezier{40}(32.67,14.33)(37.67,12.00)(42.67,14.33)
\bezier{44}(77.33,14.33)(82.67,12.33)(87.67,14.33)
\bezier{44}(92.67,14.33)(97.67,12.00)(102.33,14.33)
\bezier{40}(17.67,5.67)(22.67,7.67)(27.33,5.67)
\bezier{44}(32.67,5.67)(37.33,8.00)(42.67,5.67)
\bezier{44}(77.33,5.67)(82.67,7.67)(87.67,5.67)
\bezier{44}(92.67,5.67)(97.33,8.00)(102.33,5.67)
\put(10.00,10.00){\circle*{1.33}}
\put(30.00,5.00){\circle*{1.33}}
\put(45.00,5.00){\circle*{1.33}}
\put(45.00,15.00){\circle*{1.33}}
\put(70.00,10.00){\circle*{1.33}}
\put(90.00,5.00){\circle*{1.33}}
\put(90.00,15.00){\circle*{1.33}}
\put(105.00,5.00){\circle*{1.33}}
\put(5.33,7.67){\makebox(0,0)[lb]{$a^{\BOne}$}}
\put(27.67,0.00){\makebox(0,0)[lb]{$a^{\BOne}$}}
\put(43.33,0.00){\makebox(0,0)[lb]{$a^{\BTwo}$}}
\put(43.33,17.00){\makebox(0,0)[lb]{$a^{\BTwo}$}}
\put(15,10){\makebox(0,0)[cc]{\large$\circlearrowleft$}}
\put(45,10){\makebox(0,0)[cc]{\large$\circlearrowleft$}}
\put(56.33,8.33){\makebox(0,0)[lb]{$\ne$}}
\put(65.33,7.67){\makebox(0,0)[lb]{$a^{\BOne}$}}
\put(75,10){\makebox(0,0)[cc]{\large$\circlearrowleft$}}
\put(88.67,17){\makebox(0,0)[lb]{$a^{\BTwo}$}}
\put(88.67,0.00){\makebox(0,0)[lb]{$a^{\BTwo}$}}
\put(103.33,0.00){\makebox(0,0)[lb]{$a^{\BOne}$}}
\put(105,10){\makebox(0,0)[cc]{\large$\circlearrowleft$}}
\end{picture}
\caption{The letters $a^{\BOne}$ are (not) separated by the letters $a^{\BTwo}$.
}\label{FigNassocCounterEx}
\end{center}
\end{figure}
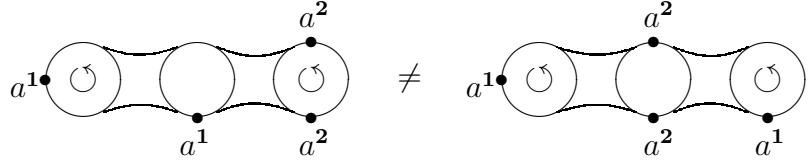
the word $(a_2)\times(a_3)$ is equal to $(a^{\BTwo}a^{\BOne}a^{\BTwo})$, 
whence the nested product $(a_1)\times\bigl((a_2)\times(a_3)\bigr)$ contains the term
$\tfrac{1}{3}a^{\BOne}a^{\BTwo}a^{\BOne}a^{\BTwo}$, which is absent in the left\/-\/hand side of~\eqref{EqNassoc} for these $a_1$,\ $a_2$,~$a_3$.
\end{counterexample}

\begin{convent}
By interpreting the ground field~$\Bbbk$ as the linear span of the zero\/-\/length word~$\bun$
and its equivalence class~$(\bun)$, we extend the commutative algebra of cyclic words to $\cA\oplus\Bbbk\cdot(\bun)$, now
endowed with the multiplication~$\times$ such that, in 
agreement with the vector space structure of~$\cA$, 
formula~\eqref{EqDefMult} is extended by 
\begin{equation}\label{EqDefLength0}
(k)\times(a)\mathrel{\stackrel{\text{def}}{=}} k\cdot(a)
\end{equation}
for any~$k\in\Bbbk$ and all cyclic words~$(a)$. 
Allowing for the slightest abuse of notation, we continue denoting by~$\cA$
the unital algebra of cyclic words that contains 
such zero\/-\/length but non\/-\/empty strings of symbols.
\end{convent}




\begin{OpenProblem}[prime decomposition]
Is there a way to detect that a given sum~$(a)\in\cA$ of several cyclic words is the product~$(b)\times(c)$ of two shorter cyclic words~$(b)$,\ $(c)\in\cA$ of positive length\,?
\end{OpenProblem}

Let us give several examples of natural constructions of 
the algebra~$\cA$ that contains nonnegative\/-\/length cyclic words written in an 
alphabet~$a^1$,\ $\ldots$,\ $a^m$. 
By realising every such algebra as fibre in a 
bundle~$\piNC$ over a given manifold~$M^n$ (e.g., in the trivial bundle over a finite\/-\/dimensional affine real manifold, cf.~\S\ref{SecWhy}), we shall proceed in~\S\ref{SecJets}
with the construction of the space $J^\infty(\piNC)$ of infinite jets of sections for such bundles~$\piNC$.

\begin{model}\label{ModelMatr}
Consider the algebra $\Mat(\mathsf{n},\BBR)$ of square matrices of size~$\mathsf{n}\times\mathsf{n}$ with real entries. Roughly speaking, as $\mathsf{n}\to+\infty$, the matrix multiplication~$\circ$ will never become commutative (yet it always stays associative). For definition, let $m\mathrel{{:}{=}}\mathsf{n}^2$ be the dimension of entire matrix algebra and choose a basis~$a^1$,\ $\ldots$,\ $a^m$ in~it.
Although this $\BBR$-\/algebra is not free, we still introduce the linear equivalence relation~$\sim$ on the vector space of words written in the alphabet~$\ba=\langle a^1$,\ $\ldots$,\ $a^m\rangle$, which yields the cyclic word algebra~$\cA$.

Because the matrix multiplication is not commutative, the content of every 
cyclic word $(a)=(a^{i_1}\circ\dots\circ a^{i_\lambda})$ of length~$\lambda>0$, viewed as the actual product of~$\lambda$ matrices going in a specified sequential order, can take up to~$\lambda$ different values, namely,
\begin{equation}\label{EqManyValues}
a^{i_1}\circ\dots\circ a^{i_\lambda},\qquad
\gotht\bigl(a^{i_1}\circ\dots\circ a^{i_\lambda}\bigr),\quad \ldots, \quad
\gotht^{\lambda-1}\bigl(a^{i_1}\circ\dots\circ a^{i_\lambda}\bigr).
\end{equation}
The value depends on the place where the multiplication is started along the orientation of the cycle (see the figure).
\begin{figure}[htb]
{\unitlength=1mm
\centerline{
\begin{picture}(20,20)(-10,-10)
\put(0,0){\circle*{1.5}}
\put(-8,0){\line(1,0){16}}
\put(0,-10){\line(0,1){20}}
\put(2,9){$\lambda=1$}
\end{picture}
\qquad
\begin{picture}(20,20)(-10,-10)
\put(0,4){\circle*{1.5}}
\put(0,-4){\circle*{1.5}}
\put(-8,0){\line(1,0){16}}
\put(0,-10){\line(0,1){20}}
\put(2,9){$\lambda=2$}
\end{picture}
\qquad
\begin{picture}(20,20)(-10,-10)
\put(0,5){\circle*{1.5}}
\put(0,-5){\circle*{1.5}}
\put(0,0){\circle*{1.5}}
\put(-8,0){\line(1,0){16}}
\put(0,-10){\line(0,1){20}}
\put(2,9){$\lambda=3$}
\end{picture}
\qquad
\begin{picture}(20,20)(-10,-10)
\put(0,7.5){\circle*{1.5}}
\put(0,-7.5){\circle*{1.5}}
\put(0,2.5){\circle*{1.5}}
\put(0,-2.5){\circle*{1.5}}
\put(-8,0){\line(1,0){16}}
\put(0,-10){\line(0,1){20}}
\put(2,9){$\lambda=4$}
\end{picture}
\quad\raisebox{10mm}{\text{etc.}}
}
}
\end{figure}
This effect --\,the value of a word~$(a)$ of length $\lambda>0$ can co\/-\/exist in $\mathsf{s}\leqslant\lambda$ realisations\,-- will be natural to the other two models which we consider below. Reproduced verbatim by the star\/-\/product~$\star$ in Model~\ref{ModelStar}, such value multiplicity can be suppressed ($1\leqslant \mathsf{s}\leqslant\lambda$ so that the first equality is attained and the last inequality is strict if~$\lambda>1$) in the model of walks, e.g., along closed contours~$a^i$ from a point to itself within a given manifold (see~\S\ref{SecWhy}).

Now let $M^n$~be a real manifold and $\piNC\colon M^n\times\cA\to M^n$ be the trivial bundle.
By construction, sections of $\piNC$~viewed as noncommutative bundle are obtained as follows.
First, let $\bs=\bigl(s^1(\bx)$,\ $\ldots$,\ $s^m(\bx)\bigr)$ be a tuple of 
functions from~$C^\infty(M^n\to\BBR)$ (e.g., compact\/-\/supported over~$M^n$). A tuple~$\bs$ chosen, over every~$\bx\in M^n$ the $i$th generator~$a^i$ of the matrix algebra~$\Mat(\mathsf{n},\BBR)$ is taken with the real coefficient~$s^i(\bx)$. Likewise, every product $a^{i_1}\circ\ldots\circ a^{i_\lambda}$ acquires the coefficient $s^{i_1}\cdot\ldots\cdot s^{i_\lambda}$. Finally, such coefficient is passed through~$\sim$ to the quotient~$\cA$ modulo the linear equivalence, pointwise over~$\bx\in M^n$.
So, all cyclic words in~$\cA$ are weighted by smooth real coefficients, depending on points of~$M^n$, in such a way that the multiplication of cyclic words is respected by those weights.
\end{model}

\begin{model}\label{ModelStar}
Likewise, let $M^n$~be a finite\/-\/dimensional affine real Poisson manifold 
and $\star=\cdot+\hbar\,\{\,,\,\}_{\cP}+\bar{o}(\hbar)$ be the arising 
associative non\/-\/commutative star\/-\/product in the unital algebra $C^\infty(M^n\to\BBR)[[\hbar]]$ of formal power series (see~\cite{KontsevichFormality}; an expansion~$\star$ mod~$\bar{o}(\hbar^4)$ is given in~\cite{cpp}). 
Keeping in mind the linearity of~$\star$ over
~$\hbar$, suppose $a^1$,\ $\ldots$,\ $a^m\in C^\infty(M^n\to\BBR)[[\hbar]]$. Using the addition and $\star$-\/product, generate from this (in)finite alphabet and~$\hbar$ a unital subalgebra of nonnegative\/-\/length words~$1$,\ $\hbar$,\ $\ldots$,\ $a^i$,\ $\ldots$,\ $a^{i_1}\star\dots\star a^{i_\lambda}$,\ $\ldots$, and pass to the quotient algebra~$\cA$ of cyclic words. (Our earlier remark that every such homogeneous word $\bigl(a^{i_1}\star\dots\star a^{i_\lambda}\bigr)$ can co\/-\/exist in up to~$\lambda$ different values is still in order.) Now, the construction of the noncommutative bundle~$\piNC$ of cyclic\/-\/word algebras~$\cA$ over the affine manifold~$M^n$ at hand is immediate; 
its section is a choice which function from~$C^\infty(M^n\to\BBR)[[\hbar]]$ each element~$a^i$ of the alphabet is equal to. Whenever all the elements of the alphabet are compact\/-\/supported over the base manifold~$M^n$, so are all the cyclic words.
\end{model}

An outline of the third model is stretched over several sections; it will be concluded on p.~\pageref{ModelBV} by comparing the result with the standard graded\/-\/commutative geometry of the Batalin\/--\/Vilkovisky (BV) superbundle~$\boldsymbol{\zeta}^{(0|1)}$. Let us specify at once that the sheaf~$\MnnC$ of algebras of walks (introduced in~\S\ref{SecWhy}) and realisation of sections in~$\piNC$ as the (cyclic) word algebra mappings 
in~\S\ref{SecAs} are pertinent to this model. At the same time, the construction of the symplectic\/-\/dual variables~$a^\dagger_i$ in~\S\ref{SecCotangent} and of their parity\/-\/odd neighbours~$b_i$ (see~\S\ref{SecNeighbour}) is common to all the 
models.\footnote{\label{FootWhyNoEdgesInDensity}%
This is why from \S\ref{SecLocalFunctionals} onwards, we shall assume that densities of integral functionals over the jet superspace $J^\infty(\bpiNCZO)$ do not depend explicitly on the edge alphabet~$\bvx^{\,\pm1}$ of a tiling of the base manifold~$M^n$ underlying the noncommutative superbundle~$\bpiNCZO$. Indeed, the availability of such edge alphabet is a feature of the third model, which we presently discuss.}

\subsection{The sheaves of algebras of walks}\label{SecWhy}
In this section we motivate the construction of the algebra~$\cA$ that contains nonnegative\/-\/length cyclic words written in the alphabet~$a^1$,\ $\ldots$,\ $a^m$. By introducing several new elements into the picture now, in~\S\ref{SecJets} 
we shall recover the notion of space of infinite jets $J^\infty(\piNC)$ of sections of the noncommutative bundle~$\piNC$ in which the algebra~$\cA$ provides the fibres.

\begin{figure}[hb]
\begin{center}
{\unitlength=1mm
\linethickness{0.4pt}
\begin{picture}(146.67,40.33)(0,-3)
\put(5.00,5.00){\line(1,0){45.00}}
\put(25.00,15.00){\line(1,0){30.00}}
\put(5.00,25.00){\line(1,0){55.00}}
\put(15.00,35.00){\line(1,0){40.00}}
\put(25.00,15.00){\line(1,1){25.00}}
\put(25.00,15.00){\line(-1,-1){13.67}}
\put(11.33,28.67){\line(1,-1){27.33}}
\put(35.00,5.00){\line(1,1){23.33}}
\put(35.00,5.00){\line(-1,-1){4}}
\put(15.00,5.00){\line(-1,1){13.67}}
\put(15.00,5.00){\line(1,-1){4.00}}
\put(15.00,25.00){\line(1,1){15.33}}
\put(15.00,25.00){\line(-1,-1){13.33}}
\put(45.00,15.00){\line(-1,1){25.00}}
\put(45.00,15.00){\line(1,-1){5.00}}
\put(55.00,25.00){\line(-1,1){14.33}}
\put(55.00,25.00){\line(1,-1){4.00}}
\put(13.33,13.67){\makebox(0,0)[lb]{$\Delta_1$}}
\put(23.00,27.67){\makebox(0,0)[lb]{$\Delta_6$}}
\put(23.00,19.33){\makebox(0,0)[lb]{$\Delta_2$}}
\put(23.00,7.33){\makebox(0,0)[lb]{$\Delta _5$}}
\put(32.33,17.33){\makebox(0,0)[lb]{$\Delta_3$}}
\put(32.33,10.00){\makebox(0,0)[lb]{$\Delta_4$}}
\put(22.00,-5.00){\makebox(0,0)[lb]{$\{a\}$}}
\put(0,0){
\begin{picture}(0,0)(-10,0)
\put(58.33,19.33){\vector(1,0){16.67}}
\put(62.33,20.00){\makebox(0,0)[lb]{dual}}
\end{picture}}
\put(0,0){
\begin{picture}(0,0)(-20,0)
\put(80.00,15.00){\line(2,1){10.00}}
\put(90.00,20.00){\line(2,-1){10.00}}
\put(100.00,15.00){\line(0,-1){5.00}}
\put(100.00,10.00){\line(-2,-1){10.00}}
\put(90.00,5.00){\line(0,-1){5.00}}
\put(90.00,5.00){\line(-1,1){9.67}}
\put(80.33,14.67){\line(-4,-3){5.33}}
\put(80.00,15.00){\line(-2,1){5.00}}
\put(90.00,20.00){\line(0,1){5.00}}
\put(90.00,25.00){\line(-1,1){5.00}}
\put(90.00,25.00){\line(2,1){10.00}}
\put(100.00,30.00){\line(0,1){7.67}}
\put(100.00,30.00){\line(2,-1){10.00}}
\put(110.00,25.00){\line(0,-1){5.00}}
\put(110.00,20.00){\line(-2,-1){10.00}}
\put(110.00,25.00){\line(1,1){5.00}}
\put(110.00,20.00){\line(5,-3){5.00}}
\put(100.00,10.00){\line(2,-1){5.00}}
\put(80,15){\circle*{1}}
\put(90,20){\circle*{1}}
\put(90,25){\circle*{1}}
\put(90,5){\circle*{1}}
\put(100,10){\circle*{1}}
\put(100,15){\circle*{1}}
\put(100,30){\circle*{1}}
\put(110,20){\circle*{1}}
\put(110,25){\circle*{1}}
\put(84.00,24){\makebox(0,0)[lb]{$\Delta_6$}}
\put(73.00,13){\makebox(0,0)[lb]{$\Delta _1$}}
\put(84.00,3.00){\makebox(0,0)[lb]{$\Delta_5$}}
\put(88.33,15.00){\makebox(0,0)[lb]{$\Delta_2$}}
\put(97.5,16.67){\makebox(0,0)[lb]{$\Delta_3$}}
\put(98,5){\makebox(0,0)[lb]{$\Delta_4$}}
\put(90.67,-5.00){\makebox(0,0)[lb]{$\{b\}$}}
\end{picture}}
\end{picture}
\caption{A fragment of cell\/-\/complex tiling~(\textit{a}) and
its adjacency graph~(\textit{b}).
}\label{FigVoronoi}%
}\end{center}
\end{figure}
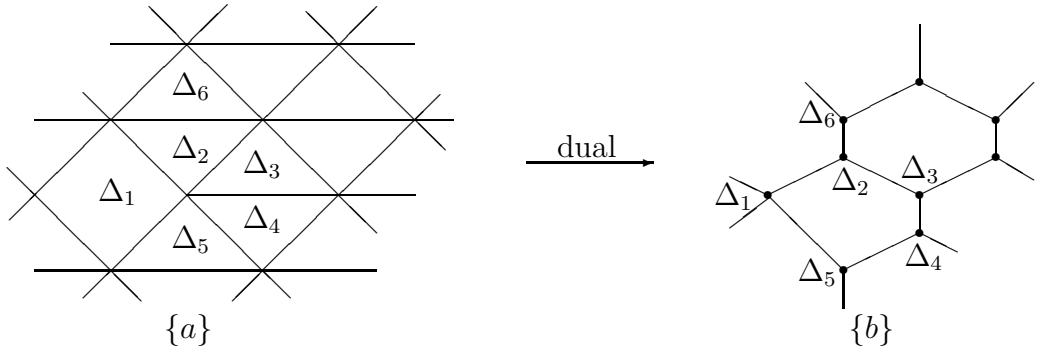

Let $M^n$ be an oriented affine 
real manifold of positive dimension~$n$. Suppose that a 
tiling of the manifold~$M^n$ is given, that is, $M^n$~is realised
by $M^n=\cup_{\alpha\in\mathcal{I}}\overline{\Delta}_{\alpha}$ via\footnote{The closure
$\overline{\Delta}_{\alpha}$ of each cell $\Delta_\alpha$ is taken with respect to the 
Euclidean topology on the manifold~$M^n$ under study.} 
the complex of cells $\Delta_{\alpha}$ of dimension $n$, see Fig.~\ref{FigVoronoi}a. 
(Of course, the manifold~$M^n$ can be topologically nontrivial: e.g., roll the plane to a cylinder, respecting --\,in one of the many possible ways\,-- a given regular crystal structure on~$\mathbb{E}^2$.)
We remark also that the choice of a tiling can be not unique for a given manifold~$M^n$. Next, construct the tiling adjacency graph: each cell $\Delta_{\alpha}$
represented by the vertex in the dual picture (see Fig.~\ref{FigVoronoi}b), two vertices are connected by the edge iff the respective cells in the tiling are adjacent through a common face of lower dimension
\footnote{The \emph{discrete} adjacency table, 
\emph{finite} for every vertex $\Delta_{\alpha}$ in the dual complex, is the
main profit that one gains by taking the 
tiling of space, 
however 
tiny be the diameter of each cell with respect to a given 
distance function on~$M^n$.
The property of base manifold~$M^m$ to be affine, that is, to admit a flat structure (consisting of an atlas of charts with affine transition maps) is natural in 
this context. 
Namely. affine reparametrisations within a tiling domain amount to a change of frame's reference to a point which marks that domain.
} (that of $n-1$).

\begin{define}
Two binary operations are defined for paths along the edges between adjacent cells in a tiling: namely, the formal addition~$+$ and multiplication~$\circ$.
Whenever a (connected component of a) path~$b$ continues a (connected part of a) path~$a$, we write $a\xrightarrow{} b$. Suppose onward that $a$,\ $b$, and~$c$ are connected paths.
If $a\xrightarrow{} b$, then $a\circ b$ is the connected path obtained by using the concatenation; otherwise, we set $a\circ b=a+b$. 

The respective neutral elements for~$+$ and~$\circ$ are~$0$ and the null path $\bullet=\bun$.
\end{define}

The addition~$+$ is commutative and associative;
clearly, the multiplication of paths is 
not always commutative.

\begin{lemma}
At the same time, the multiplication of (sums of) paths is not associative. 
\end{lemma}

(Of course, if $a\xrightarrow {}b$ and $b\xrightarrow {}c$, then $(a\circ b)\circ c= a\circ(b\circ c)$.) 

\begin{proof}
Namely, if $a\xrightarrow {}b$ and $a\xrightarrow {}c$ but $b\not\xrightarrow {}c$, then $(a\circ b)\circ c = a\circ b+c$, yet the associator right\/-\/hand side is different: $a\circ(b\circ c)=a\circ(b+c)=a\circ b+a\circ c$.
\end{proof}

Finally, let us inspect the distributivity law.
\begin{itemize}
\item If $a\xrightarrow{} c$ and $b\xrightarrow{} c$, then $(a+b)\circ c= a\circ c + b\circ c$.
\item If $a\xrightarrow{} c$ but $b\not\xrightarrow{} c$, then $(a+b)\circ c = a\circ c+b+c=a\circ c+b\circ c$ as well.
\item If $a\not\xrightarrow{} c$ nor $b\not\xrightarrow{} c$, then $(a+b)\circ c=a+c+b+c=a\circ c + b\circ c$.
\end{itemize}
And now, the other way round:
\begin{itemize}
\item If $a\xrightarrow{} b$ and $a\xrightarrow{} c$, then $a\circ(b+c)= a\circ c + a\circ c$.
\item If $a\xrightarrow{} b$ but $a\not\xrightarrow{} c$, then $a\circ(b+c) = a\circ b+a+c=a\circ b+a\circ c$ as well.
\item If $a\not\xrightarrow{} b$ nor $a\not\xrightarrow{} c$, then $a\circ(b+c)=a+b+a+c=a\circ b + a\circ c$.
\end{itemize}

The \emph{unital algebra of walks} is the vector space of formal sums (with respect to the addition~$+$) of paths that are multiplied by using the concatenation~$\circ$; both the operations are proclaimed $\Bbbk$-\/linear.

Without any extra assumptions made about the tiling, the cells adjacency table and the portrait of edges in the dual graph are \emph{local}. Indeed, a
\textit{quasi}crystal structure of the cell complex realisation of~$M^n$ could contain defects. Consequently, the larger
an open domain $U\subseteq M^n$~is, the larger can be the alphabet of edges which are used to encode paths within~$U$ as words. 

For the sake of definition, we 
assume 
that the substrate manifold's 
tiling is globally regular,
so that the crystal structure $\{\Delta_{\alpha}\}$ is formed by (in)\/finite replication of a finite union of cells. 

\begin{define}
The edge alphabet is any minimal (i.e.\ without repetitions) subset of the set of edges such that every walk between cells in a given tiling of~$M^n$ can be expressed using that subset. Up to a permutation of edges and up to a choice --\,which direction of an $i$th edge is denoted 
by the symbol~$\vx_i^{\,+1}$ and the other by~$\vx_i^{\,-1}$,\,-- every alphabet~$\bvx^{\,\pm1}$ consists of
\begin{enumerate}
\item all the edges connecting the cells in their finite union which is replicated so that the tiling is made, and
\item the edges which interconnect that generating union of cells with all those replicas which are adjacent to that union (cf.~\cite{ConwaySloane,LoopReview}).
\end{enumerate}
\end{define}

\noindent\begin{minipage}[b]{117mm}\noindent%
\begin{example}
Consider the honeycomb 
tiling 
of the plane, see the figure. 
The regularity assumption makes the alphabet~$\bvx^{\pm1}$
finite even if the tiling of the (non)\/compact manifold~$M^n$ is infinite.
We denote by~$N$ the cardinality of the set of generators, so that the chosen alphabet is
$\bvx^{\,\pm1}=\langle\vx_1^{\,\pm1}$, $\ldots$, $\vx_N^{\,\pm1}\rangle$.
Let us remember that the number~$N$ of elementary displacements~$\vx_i$ depends on the choice of a tiling for the affine manifold~$M^n$ of dimension~$n$.
Now, the price that one 
pays is that the coding of edges can no longer be referred to any specific cell, hence a presence of irregular, non\/-\/periodic defects is no longer possible.
\end{example}
\end{minipage}
{\unitlength=1mm
\linethickness{0.4pt}
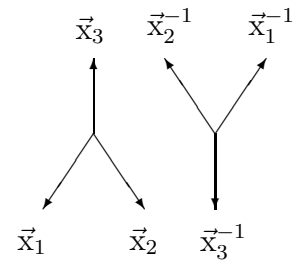
\begin{picture}(0,0)(114,-5)
\put(124.00,20.00){\vector(0,1){10.00}}
\put(124.00,20.00){\vector(-2,-3){6.67}}
\put(124.00,20.00){\vector(2,-3){6.67}}
\put(140.00,20.00){\vector(0,-1){10.00}}
\put(140.00,20.00){\vector(-2,3){6.67}}
\put(140.00,20.00){\vector(2,3){6.67}}
\put(121.67,32.00){\makebox(0,0)[lb]{$\vx_3$}}
\put(131.00,32.00){\makebox(0,0)[lb]{$\vx_2^{-1}$}}
\put(144.33,32.00){\makebox(0,0)[lb]{$\vx_1^{-1}$}}
\put(114.00,4.00){\makebox(0,0)[lb]{$\vx_1$}}
\put(128.67,4.00){\makebox(0,0)[lb]{$\vx_2$}}
\put(138.00,3.33){\makebox(0,0)[lb]{$\vx_3^{-1}$}}
\end{picture}}%


From now on, let an alphabet $\bvx^{\,\pm1}=\langle\vx_i^{\,\pm1}\rangle$ be fixed for a given crystal tiling of the affine manifold~$M^n$ under study. For every value of the index~$i$, 
the symbols~$\vx_i^{\,+1}$ and~$\vx_i^{\,-1}$ denote the edges passed in the adjacency graph in either of the two directions.\footnote{A possibility to walk every edge, hence every path \emph{backwards} --\,along the respective reverses~$\vx_i^{\,\mp1}$, reading the words right to left,\,--
is a forerunner of the introduction of canonical conjugate symbols~$a^\dagger_j$, which are responsible for the dual, parity\/-\/odd part of the 
picture. This will be discussed in~\S\ref{SecCotangent} and~\S\ref{SecNeighbour}, see Fig.~\ref{FigSetup} on p.~\pageref{FigSetup} in particular.}

Over the substrate manifold~$M^n$ let us construct the almost constant sheaf (see~\cite{ManinSchemesQG,MumfordRedBook}) of unital extensions 
$\Bbbk\cdot\BOne\oplus\Free_{\Bbbk}\bigl(\vx^{\pm1}_1,\ldots,\vx^{\pm1}_N\bigr)$ of free algebras generated by~$\vx_i^{\,+1}$ and~$\vx_i^{\,-1}$. 
The sheaf is glued --\,from such unital algebras of walks\,-- over open subsets $U\subseteq M^n$.
For all pairs $U_j\subseteq U_i$ of non\/-\/empty open subsets of the set~$M^n$ with a chosen topology (e.g., the Euclidean one), the restriction homomorphisms are the identity mapping unless $U_j\subseteq \Delta_\alpha$ for some cell, marked by~$\alpha\in\cI$ in the tiling; in that case, the sheaf structure over~$U_j$ is the null path component~$\Bbbk\cdot\bun$ and the restriction mapping is the canonical projection.
Over the empty subset of~$M^n$, the sheaf structure is empty by definition.

\begin{notation}
This sheaf over~$M^n$ will be denoted by~$\MnnC$; it remembers the topology on the substrate manifold and it carries the finite alphabet $\bvx^{\pm1}$ of the $N$~edges that interconnect 
cells in (the replicas of) a fundamental domain in the tiling.
\end{notation}

\subsection{The formal fibre~$\cA$ of~$\piNC$}\label{SecAs}
We start building a noncommutative analogue~$\bpiNC$ of the variational cotangent bundle over~$M^n$ and then, using that noncommutative object, 
\begin{minipage}[b]{73mm}
the analogue~$\bpiNCZO$ of the Batalin\/--\/Vil\-ko\-vi\-sky superbundle over
the space\/-\/time. Recalling from~\S\ref{SecAlgebra} 
the construction of the algebra~$\cA$ of cyclic words, we notice that 
whenever such algebra is realised as a fibre, it suffices to evaluate the generators~$a^i$ of the free algebra. Such ``sections'' (that is, the generator evaluation mappings)
are then extended onto (the quotient of) the target space $\Free(a^1,\,\dots\,,a^m)$ by using both the multiplication~$\circ$ and addition~$+$. 
Indeed, consider the evaluation map 
$\bs{\bigr|}_U \colon \Free(a^1,\,\dots\,,a^m) \to \MnnC{\bigr|}_U$ which, at every 
point~$\bx$ in a chart $U\subseteq M^n$ 
within the substrate manifold~$M^n$, takes
\end{minipage}\label{FigBundlePi}
{\unitlength=1mm
\begin{picture}(85,0)(-5,-2)
\put(35,16){\circle*{1}}
\put(35.5,14.5){$\bx$}
\multiput(6,19)(2.5,2.5){3}{\line(1,0){54}}
\multiput(6,19)(6,0){10}{\line(1,1){5}}
\qbezier
(7,5)(30,25)(50,5)
\put(7,5){\line(3,1){12}}
\put(50,5){\line(3,1){12}}
\qbezier
(19,9)(35,26)(62,9)
\put(58,4){$M^n$}
\put(65,15){$\left.\rule{0pt}{10mm}\right\}$}
\put(70,15){$\MnnC$}
\put(35,32){\vector(0,-1){5}}
\put(37,28){$\piNC$}
\put(35,34){\circle*{1}}
\put(38,33){$\oplus\:\langle\bun\rangle_{\Bbbk}$}
\put(35,34){\line(-2,5){8}}
\put(35,34){\line(1,3){12}}
\put(27,54){\line(5,4){20}}
\put(35,47){\vector(0,1){10}}
\put(31.5,48){\vector(1,1){8}}
\put(36.5,56){$\ba$}
\put(27,54){\line(-1,2){8}}
\put(19,70){\line(1,0){28}}
\put(21,71.5){$\text{length}\,(a)\geqslant 1$}
\qbezier(30,30)(25,36)(30,42)
\put(30,42){\vector(1,1){0}}
\put(30,30){\vector(1,-1){0}}
\put(26,35){\llap{$\bolds$}}
\qbezier(53,60)(47,58)(41,60)
\put(41,60){\vector(-2,1){0}}
\put(54,59){$\text{span}_\Bbbk\langle\ba\rangle$}
\put(19,55){\llap{$\left\{\rule{0pt}{15mm}\right.$}}
\put(15.5,58.5){\llap{$\Free_\Bbbk(\ba)$,}}
\put(10.5,53){\llap{then}}
\put(14,47.5){\llap{$/{\sim}\mathrel{{=}{:}}\cA$}}
\qbezier[160]%
(19,70)(27,52)(35,34)
\put(32,17){\line(-3,-1){9}}
\end{picture}%
}
each generator~$a^i$ to a word of positive proper length\footnote{
Obviously, the case where $a^i=s^i(\bx)$ for some~$i$ would be somewhat special: 
the algebra~$\cA$ of nonnegative\/-\/length 
cyclic words was unital by construction, but the above assignment would convert the generator~$a^i$ to the multiple of the 
neutral element at every~$\bx$ in a chart. To exclude this 
situation from the 
study, let us technically assume that the lexicographic length of all 
the word(s) in each 
component~$s^i$ is strictly positive.
Moreover, one should even require that the \emph{walk}~$s^i$ along the edges~$\vx^{\pm1}_i$
of the graph be more than a null path~$\BOne$, for it could be that the walk is contractable: e.\,g., $s^i=\vx_j\circ\vx^{-1}_j=\BOne$.
} 
--\,or to a formal sum of such words\,-- written in the alphabet $\bvx^{\pm1}=\{\vx^{\pm1}_j$,\ $1\leqslant j\leqslant N\}$:
\begin{equation}\label{EqNCsection}
a^i=s^i(\bx,\bvx^{\pm1}),\qquad 1\leqslant i\leqslant m,
\end{equation}
each word taken with a smooth coefficient from~$C^{\infty}(M^n)$. 
Actually, 
formula~\eqref{EqNCsection} is a compact notation:\label{FootCoeffsXs} 
its right\/-\/hand side evaluates at $\bx\in M^n$ the infinitely
many coefficients of $\vx_i^{\,\pm1}$,\ $\vx^{\,\pm1}_i\circ\vx^{\,\pm1}_j$,\
$\vx^{\,\pm1}_i\circ\vx^{\,\pm1}_j\circ\vx^{\,\pm1}_k$, etc.

By construction, the value $a{\bigr|}_\bs$ of a homogeneous word~$a$ written in the alphabet $\ba=\{a^i$,\ $1\leqslant i\leqslant m\}$ is the 
product of the map values at the consecutive letters of that word. 
For instance, we postulate that
\[
\left.\mathstrut(a^i\circ a^j)\right|_{\bs}(\bx,\bvx^{\pm1})=
s^i(\bx,\bvx^{\pm1})\circ s^j(\bx,\bvx^{\pm1});
\]
the multiplication~$\circ$ in the right
\/-\/hand side is 
the 
multiplication in the sheaf of free associative algebras, so that one proceeds recursively.




\begin{convent}[irreducibility]
Let the mapping in~\eqref{EqNCsection} be such that no positive\/-\/length word $a\in\Free(a^1$,\ $\ldots$,\ $a^m)$ is evaluated to a zero proper length word $a{\bigr|}_\bs$ in the algebra of walks (i.e., the null path~$\bun$ with a nonzero coefficient from~$C^\infty(M^n)$).
(For example, we exclude the case where $a^1\mathrel{{:}{=}}\vx_1\circ\vx_2^{\,-1}$ and 
$a^2\mathrel{{:}{=}}\vx_2\circ\vx_1^{\,-1}$, so that $a^1\circ a^2{\bigr|}_\bs=\bun$ at all~$x\in M^n$).
\end{convent}

The construction of sections~\eqref{EqNCsection} is furthered 
to the quotient~$\cA=\Free(a^1$,\ $\ldots$,\ $a^m)/\sim$, which yields the evaluation mapping~$\bs$ from the 
sheaf of algebras~$\cA$ over the commutative manifold~$M^n$ to the sheaf of unital algebras~$\cX(\bvx^{\,\pm1})$ of cyclic words written in the edge alphabet for a given tiling of~$M^n$. (The restriction maps~$r^U_V$ in that sheaf, for $V\subseteq U$ open in~$M^n$, are the identity mapping of~$\cX(\bvx^{\,\pm1})$, except for the constant mapping to the null word~$(\bun)$ over~$V\subseteq\Delta_\alpha$ at some~$\alpha\in\cI$; over~$\varnothing\subset M^n$, the sheaf structure is empty.)

\begin{rem}\label{RemNotAlgHom}
Let us remember that this evaluation mapping~$\bs$ is not a homomorphism of the cyclic word algebras~$\cA$ and~$\cX(\bvx^{\,\pm1})$, respectively. 
The inequality,
\[
\bigl((a_1)\stackrel{\cA}{\times}(a_2)\bigr){\bigr|}_{\bs} (\bx) \neq 
(a_1){\bigl|}_\bs \stackrel{\cX(\bvx^{\,\pm1})}{\times} (a_2){\bigr|}_\bs (\bx),
\]
can occur for some words $(a_1)$,\ $(a_2)\in\cA$ and at some point~$\bx\in M^n$.
Indeed, the multiplication~$\times$ in~$\cA$ unlocks the cyclic words in between the letters~$\ba$ 
that will later be evaluated 
using~\eqref{EqNCsection}, whereas the multiplication~$\times$ in~$\cX(\bvx^{\pm1})$ unlocks the cyclic words
between \emph{every} two consecutive symbols from the edge alphabet~$\bvx^{\,\pm1}$
(see also Remark~\ref{RemTwoMult} on p.~\pageref{RemTwoMult} below).
\end{rem}


\begin{rem}
Evaluation~\eqref{EqNCsection} of a word~$a$ from~$\Free(\ba)$ 
paves 
the way 
(weighted by elements of~$C^{\infty}(M^n)$) 
along the edges~$\vx^{\,\pm1}_i$ of the graph which we started with. 
If the path $a{\bigr|}_{\bs}$ is closed, then it does not matter where one starts reading that cyclic word (now written in the alphabet~$\bvx^{\,\pm1}$); hence the value $(a){\bigr|}_{\bs}(\bvx^{\,\pm1})$ is uniquely defined.
However, the cyclic invariance of the word~$(a)$ does \emph{not} 
imply that the path~$a{\bigr|}_{\bs}$ is closed.%
\footnote{Alternatively, it could require 
some effort to make 
a given value a cyclic cyclic word indeed 
by contracting the graph between the path loose ends.}
Strictly speaking, 
not every word written in the alphabet $\bvx^{\pm1}_i$ encodes some path connecting
cells in the 
tiling. (Still the converse is true: every path is encoded by the respective word and every closed path --\,written by using the alphabet~$\ba$ and map~\eqref{EqNCsection}\,-- is described by the equivalence class of cyclic words.)

It is readily seen that for a word~$a$ of length~$\lambda>0$, the evaluation of~$(a)$ by using~\eqref{EqNCsection} can produce up to $\lambda$~different 
elements in the space of cyclic words~$\cX(\bvx^{\,\pm1})$.
Such co\/-\/existence of the value $(a){\bigr|}_{\bs}$ of a given cyclic word~$(a)$ in several states occurs due to the noncommutativity of the 
concatenation~$\circ$ of 
words in the edge alphabet~$\bvx^{\,\pm1}$.

Let us remember that the same multiple\/-\/value effect can also be produced (moreover, regardless of the availability of an edge alphabet) whenever the multiplication~$\cdot$ of coefficients in~\eqref{EqNCsection} --\,or in~\eqref{EqSetACrossAtX} in what follows\,-- is replaced by using a noncommutative associative star\/-\/product~$\star$ on the affine manifold~$M^n$ (see pp.~\pageref{ModelStar} and~\pageref{RemMergeModels}).
\end{rem}

\begin{rem}
The cyclic shift operation~\eqref{EqDefShift} on $\Free(a^1,\,\dots\,,a^m)$ descends to the identity mapping $(a)\mapsto(a)$ on the algebra~$\cA$ of cyclic words~$(a)$. In what follows --\,in particular, starting from the moment when the alphabet is $\BBZ_2$-\/graded, $\gotht(\gamma_1\circ\ldots\circ\gamma_\lambda)=(-)^{|\gamma_1\circ\ldots\circ\gamma_{\lambda-1}|\cdot|\gamma_\lambda|}\gamma_\lambda\circ\gamma_1\circ\ldots\circ\gamma_{\lambda-1}$ 
for $\lambda\geqslant2$, so that the restriction of~$\gotht$ on~$\cAZO$ is not just the identity\,--
we shall \emph{not} attempt viewing the mapping~$\gotht$ as conjugation $\gamma_\lambda\circ(\gamma_1\circ\ldots\circ\gamma_\lambda)\circ\gamma_\lambda^{-1}$ whenever the rightmost comultiple is well defined.
In terms of the algebra of walks on a lattice such conjugations would mean that the entire contour $\gamma_1\circ\ldots\circ\gamma_\lambda$ is first displaced by~$\gamma_\lambda$, then read in full, and followed by a step back. This is what one should avoid, especially on \emph{irregular} lattices.
Conversely, we shall always view the shift~$\gotht$ 
as a replacement of the marker~$\pmb{\infty}$ at which one begins reading a given cyclic word.
\end{rem}

\begin{rem}[$\BOne({\bx})\in C^{\infty}(M^n)$]
As soon as the unital algebra~$\cA$ of cyclic words is placed over the ``points'' of $\MnnC$ --\,
in earnest, over usual points $\bx\in M^n$ of the substrate manifold\,-- 
the zero\/-\/length words in~$\cA$ are weighted 
pointwise over $M^n$ by elements of the ring $C^{\infty}(M^n)$ that now plays the r\^ole of the ground field $\Bbbk$.
This blow-up $\Bbbk\hookrightarrow C^{\infty}(M^n)$ is standard in the differential calculus on (jet) bundles in the commutative case 
(cf.~\cite{TwelveLectures,Norway,Galli10,Olver1993}
).
\end{rem}


\subsection{The geometry of jet space $J^{\infty}(\piNC%
)$}\label{SecJets}
Now we recall 
the standard construction of infinite jet space 
$J^{\infty}(\piNC
)$ towered over 
the substrate manifold~$M^n$ and 
the sheaf~$\MnnC$.
We emphasize that this construction (local with respect to~$\bx\in U\subseteq M^n$) refers only to the 
affine 
structure on the domain set~$
M^n$ and to the vector space organisation of 
objects over~it.


Expansion~\eqref{EqNCsection} 
yields the infinite jet alphabet which consists of 
$a^i\equiv a^i_{\varnothing}$ and $a^i_{x^j}$,\ $a^i_{x^jx^k}$,\ $\dots$,\ $a_{\sigma}$ for 
$|\sigma|\geqslant0$ over a chart $U\subseteq M^n$ with local coordinates $\bx=(x^1,\dots,x^n)$; here $\sigma$~is a multiindex.
%
%
The evaluation mappings $\ba=\bolds(\bx$,\ $\bvx^{\,\pm1})$ are extended for all~$|\sigma|\geqslant0$ by
$\ba_\sigma=\bigl(\frac{\dd^{|\sigma|}}{\dd\bx^\sigma}\bolds\bigr)(\bx$,\ $\bvx^{\,\pm1})  $
using the jets~$\text{jet}^\infty(\bolds)$. Under the assumption that the base manifold~$M^n$ be affine, the jet letters~$\ba_\sigma$ are well behaved under a change $\bx=\bx(\widetilde{\bx})$ of local coordinates.
Let us denote by~$[\ba]$ the differential dependence on letters $a^i$,\ $a^i_{x^j}$,\ $\dots$,\ $a_{\sigma}$ up to some arbitrarily high 
but always finite order~$|\sigma|<\infty$. The construction of the 
algebra $\cF
(\piNC
)$ 
of cyclic\/-\/word valued functions on $J^{\infty}(\piNC
)$
is standard: namely, it is the inductive limit of filtered algebras (\cite{TwelveLectures,Olver1993}).
Likewise, the \emph{total derivatives} $\frac{\Id}{\Id x^i}$, which we denote synonymically by $D_{x^i}$ for $1\leqslant i\leqslant n$ making no further
distinction between $\left(\frac{\Id}{\Id\bx}\right)^{\sigma}$ and $D_{\bx}^{\sigma}$, are introduced by using the restrictions
of elements $f\in\cF
(\piNC
)$ to `graphs' of~\eqref{EqNCsection}, i.e.
\begin{equation}\label{EqDefTD}
\left.\frac{
{\Id}}{\Id x^i}(f)\right|_{\text{jet}^{\infty}(\ba=\bs(\,\cdot\,,\bvx^{\pm1}))}(\bx_0)
\mathrel{\stackrel{\text{def}}{=}}
\left.\frac{
{\dd}}{\dd x^i}\right|_{\bx_0}
\left(\left.\mathstrut f\right|_{\text{jet}^{\infty}(\ba=\bs(\,\cdot\,,
\bvx^{\pm1}))}\right).
\end{equation}
This determines the usual coordinate expressions for $1\leqslant i\leqslant n$,
\begin{align*}
\frac{\overrightarrow{\Id}}{\Id x^i} &=
\frac{\dd}{\dd x^i}+\sum_{j=1}^m\sum_{|\sigma|\ge0}a^j_{\sigma\cup\{i\}}\,
\frac{\overrightarrow{\dd}}{\dd a^j_{\sigma}}\\
\intertext{which starts at~$\binfty$ and acts 
along the orientation of every cyclic word, and}
\frac{\overleftarrow{\Id}}{\Id x^i}&=
\frac{\dd}{\dd x^i}+\sum_{j=1}^m\sum_{|\sigma|\ge0}\frac{\overleftarrow{\dd}}{\dd a^j_{\sigma}}\,
a^j_{\sigma\cup\{i\}},
\end{align*}
which acts from $\binfty$ clockwise.
Both the operators $\overleftarrow{D}\!{}_{x^i}$ and $\overrightarrow{D}\!{}_{x^i}$
show up, 
first, through the substrate part $\BOne\cdot\dd/\dd x^i$ plus the $m$~sums --\,formally, infinite\,-- of cyclic words such that 
the derivations $\dd/\dd a^j_{\sigma}$ sit in their locks.

\begin{rem}
We have that
\[
\frac{\Id}{\Id x^i} = 1\ %
{{\unitlength=1mm
\begin{picture}(8,8)(-4,-1)
\put(0,0){\circle{8}}
\put(-4.25,0){\circle*{1}}
\put(4.25,0){\circle*{1}}
\put(-1.5,-1){$\circlearrowleft$}
\end{picture}}}\ %
\frac{\dd}{\dd x^i}
+ \sum_{j,\sigma} a^j_{\sigma\cup\{i\}}\ %
\text{{\unitlength=1mm
\begin{picture}(8,8)(-4,-1)
\put(0,0){\circle{8}}
\put(-4.25,0){\circle*{1}}
\put(4.25,0){\circle*{1}}
\put(-1.5,-1){$\circlearrowleft$}
\end{picture}}}\ %
\frac{\dd}{\dd a^j_\sigma},
\]
whence all the terms which are produced from 
the counterclockwise action of $\overrightarrow{\Id}/\Id x^i$ via the Leibniz rule on a given cyclic word $f\in\cF(\piNC)$ have the shape%
\footnote{If the base coordinates~$x^k$ are not considered as symbols of any alphabet at hand, then the entire coefficient $\in C^\infty(M^n)$ of the cyclic word $f\in\cF(\piNC)$ can be placed at the lock~$\binfty$.}
\begin{equation}\label{EqTermwiseActDeriv}
\text{\raisebox{1mm}[10mm][4mm]{
$\cdots+\ \ $
{\unitlength=1mm
\begin{picture}(21,16)(-4,-1)
\put(0,2){\oval(8,8)[t]}
\put(0,-2){\oval(8,8)[b]}
\put(-4,2){\line(0,-1){4}}
\put(4,2){\line(1,0){4}}
\put(4,-2){\line(1,0){4}}
\put(12,2){\oval(8,8)[t]}
\put(12,-2){\oval(8,8)[b]}
\put(16,2){\line(0,-1){4}}
\put(12,6){\circle*{1}}
\put(10.5,7){\tiny$\binfty$}
\put(8,0){\circle*{1}}
\put(9,-1){{\small$x^i$}}
\put(4.75,-1.25){\llap{\tiny $\vec{\dd}/\dd x^i$}}
\put(4.75,0){\vector(1,0){2.5}}
\put(15,-1.5){\vector(0,1){3}}
\put(-4,0){\circle*{1}}
\put(-5,-1.35){\llap{$1$}}
\end{picture}}
$+\cdots+\qquad\ \:$
{\unitlength=1mm
\begin{picture}(21,16)(-4,-1)
\put(0,2){\oval(8,8)[t]}
\put(0,-2){\oval(8,8)[b]}
\put(-4,2){\line(0,-1){4}}
\put(4,2){\line(1,0){4}}
\put(4,-2){\line(1,0){4}}
\put(12,2){\oval(8,8)[t]}
\put(12,-2){\oval(8,8)[b]}
\put(16,2){\line(0,-1){4}}
\put(12,6){\circle*{1}}
\put(10.5,7){\tiny$\binfty$}
\put(8,0){\circle*{1}}
\put(9,-1){{\small$a^j_\sigma$}}
\put(4.95,-1.1){\llap{\tiny $\vec{\dd}/\dd a^j_\sigma$}}
\put(4.75,0){\vector(1,0){2.5}}
\put(15,-1.5){\vector(0,1){3}}
\put(-4,0){\circle*{1}}
\put(-5,-1.35){\llap{$a^j_{\sigma\cup\{i\}}$}}
\end{picture}}
$+\cdots$.}}
\end{equation}
(The derivations proceed along the orientation of the argument~$f$, acting on the symbols in front of which the cyclic word~$f$ is disrupted.)
This shows 
that the operation 
$\Id/\Id x^i\otimes f\mapsto \Id/\Id x^i(f)$ is 
again a topological pair of pants~$\BBS^1\times\BBS^1\to\BBS^1$.
\end{rem}

\newpage
\section{
Differential graded Lie algebra of noncommutative local functionals}\label{SecKinematics}

\subsection{The variational symplectic dual}\label{SecCotangent}
We shall presently extend the alphabet $a^1$,\ $\dots$,\ $a^m$ of the associative algebra $\Free_{\Bbbk}(a^1$,\ $\dots$,\ $a^m)$ which we
started with. Namely, we introduce the new symbols $a^{\dagger}_1$,\ $\dots$,\ $a^{\dagger}_m$ that \emph{ought to be} the
canonical conjugates of the respective variables $a^1$,\ $\dots$,\ $a^m$; 
let us explain what this means. 

First, let us consider the free associative algebra standing alone, that is, \emph{before} the evaluation of generators 
by~\eqref{EqNCsection} under a given map $\bs
$. In this set\/-\/up, there still
remain two ways to understand the nature of new generators $a^{\dagger}_i$, namely, the coarse and fine. 
The former is to proclaim that the vector space 
$V^{\dagger}\mathrel{{:}{=}}\Span_{\Bbbk}(a^{\dagger}_1,\dots,a^{\dagger}_m)$
is dual to the linear span $V\mathrel{{:}{=}}\Span_{\Bbbk}(a^1,\dots,a^m)$ under the $\Bbbk$-\/valued coupling; by construction, 
the elements $a^{\dagger}_i$ specify the basis dual to that of~$a^i$ in~$V$. The new letters are then incorporated 
into the set of generators of (the unital extention of) the associative algebra 
$\Bbbk\cdot1
\oplus\Free_{\Bbbk}(a^1,\dots,a^m;a^{\dagger}_1,\dots,a^{\dagger}_m)$.
This definition is sufficient (which is explained in Chapter~\ref{SecDynamics}) to make the noncommutative variational Poisson formalism work.

The fine approach is as follows; although less is required, it is still enough to construct the (non)\/commutative 
Batalin\/--\/Vilkovisky geometry. Suppose that the generators~$a^i$ of the free associative algebra undergo a 
shift by
$\delta\ba=\delta a^i\cdot\vec e_i$, where the~$m$ vectors~$\vec e_i$  constitute the adapted%
\footnote{\label{FootDiagDeform}%
In other words, only the \emph{diagonal} deformations of the associative algebra generators are now allowed.
This should be expected; for in the commutative BV-\/geometry, 
the variables~$a^i$ and~$b_i=\Pi(a^{\dagger}_i)$, see below,
describe the conjugate field\/-\/antifield or ghost\/-\/antighost pairs that stem from the different generations of Noether's
identities between the Euler\/--\/Lagrange equations of motion. Hence by construction, the variables~$a^i$ or~$b_i$ at
different values of the index~$i$ are fibre coordinates in different vector bundles, merged later to their Whitney sum
(see~\cite[\S2,\:6,\:11]{TwelveLectures} or~\cite{Prague2011} and references therein).}
basis in~$T_{\ba}V$, each of them pointing along the respective generator in the vector space 
$V=\Span_{\Bbbk}(a^1,\,\dots\,,a^m)$.
Likewise, consider the adapted basis $\vec{e}^{{}\,\dagger,i}$ in the tangent space $T_{\ba^{\dagger}}V^{\dagger}$
at the point $\ba^{\dagger}$ of the vector space $V^{\dagger}=\Span_{\Bbbk}(a^{\dagger}_1,\,\dots\,,a^{\dagger}_m)$.
We require that the frame $\vec{e}^{{}\,\dagger,i}$ be $\Bbbk$-dual to the frame $\vec{e}_i$, $1\le i\le m$, so that
the 
variation $\delta\ba^{\dagger}=\delta a^{\dagger}_i\cdot\vec{e}^{{}\,\dagger,i}$ is the canonical conjugate of the diagonal deformation $\delta\ba=\delta a^i\cdot\vec{e}_i$, see~\eqref{EqTwoCouplings} 
and~\eqref{EqNormalizeShifts} below.

\begin{rem}
In the second approach, we do not proclaim that the new symbols $a^{\dagger}_i$ are the duals of the old generators $a^i$
(or their inverses, or \emph{re}ver\-ses, cf.~\eqref{EqReverseLoop}
). In other words, we do not use the isomorphism
between the vector space $V^{\dagger}=\Span_{\Bbbk}(a^{\dagger}_1,\dots,a^{\dagger}_m)$ and the vector space
$T_{\ba^{\dagger}}V^{\dagger}$ tangent to it at a point. Note that the 
left\/-\/hand side of the isomorphism
$V^{\dagger}\simeq T_{\ba^{\dagger}}V^{\dagger}$ exploits the \emph{global} vector\/-\/space organisation of $V^{\dagger}$ whereas the right\/-\/hand side refers to its \emph{local} portrait near the point~$\ba^{\dagger}$.
This is what the Batalin\/--\/Vilkovisky and Poisson formalisms really need.
\end{rem}

So, we extend the set $a^1$,\ $\ldots$,\ $a^m$ of generators by the symbols $a^\dagger_1$, $\ldots$,\ $a^\dagger_m$: at every~$i$, the new symbol~$a^\dagger_i$ matches the respective generator~$a^i$ in the above sense. As before, we take the quotient of the free algebra
$\Free_\Bbbk\bigl(a^1$,\ $\dots$,\ $a^m$; $a^{\dagger}_1$,\ $\dots$,\ $a^{\dagger}_m\bigr)$
over the linear relation~$\sim$ of equivalence under cyclic shifts. Thus we obtain the unital commutative non\/-\/associative algebra $\bigl(\Bbbk\cdot\bun\oplus
\Free_\Bbbk(a^1$,\ $\dots$,\ $a^m$; $a^{\dagger}_1$,\ $\dots$,\ $a^{\dagger}_m)\bigr)/\sim$ of cyclic words (written now in the double alphabet).\footnote{The space of free algebra generators is, strictly speaking, not the direct sum $\text{span}_\Bbbk\langle\ba\rangle\oplus\text{span}_\Bbbk\langle\ba^\dagger\rangle$ because under a rescaling of the generators~$\ba$, the \emph{dual} letters~$\ba^\dagger$ can be rescaled inverse proportionally.%
}
We postulate that the resulting algebra of cyclic words becomes the fibre in the noncommutative bundle~$\bpiNC$ over the sheaf~$\MnnC$ of algebras of walks along a given tiling of the substrate affine manifold~$M^n$.

\begin{rem}\label{RemBVOddSection}
Let us examine how the noncommutative sections~\eqref{EqNCsection}, which evaluate~$a^i$ to~$s^i(\bx,\bvx^{\,\pm1})$ over~$\bx\in M^n$, can be extended to the double alphabet evaluation using sections~$(\bs,\bs^\dagger)$ of~$\bpiNC$.
The guiding principle that one must keep in mind is that in the Batalin\/--\/Vilkovisky (BV) formalism, the quantum action functional is constrained by the narual postulate $\langle 1 \rangle = 1$ for the averaging over sections of the BV~superbundle.
This condition implies that the objects in that formalism are effectively independent of a choice of sections by using which the new, dual variables could be evaluated at~$\bx\in M^n$. Hence the generators~$a^\dagger_i$ could acquire \emph{whatever} values; indeed, no physics depends on them at the end of the day. (If so, leaving the respective components of the sections 
unspecified would be another 
option.) 

However, we are also free to assign the values $\ba=\bs(\bx,\bvx^{\pm1})$ and
$\ba^{\dagger}=\bs^{\dagger}(\bx,\bvx^{\pm1})$ in a way we choose.\footnote{The \emph{fourth} scenario is specific to the (non)\/commutative variational Poisson formalism, in the frames of
which the symbols~$\ba^{\dagger}$ play the r\^oles of placeholders for the variational covectors that are not exact;
but still, the isomorphism $V^{\dagger}\simeq T_{\ba^{\dagger}}V^{\dagger}$ is explicitly used in the assignment $\ba^{\dagger}\mathrel{{:}{=}}\bp$ (we shall discuss this in Chapter~\ref{SecDynamics}).}
\end{rem}

\begin{convent}
For a given section~$\bolds$ of~$\piNC$,
\begin{subequations}\label{EqSuperSection}
\begin{align}
a^i &= s^i(\bx,\bvx^{\pm1})=
\sum_J f^{i,J}(\bx)\,\vx^{\,\alpha(1)}_{j_1}\circ\ldots\circ\vx^{\,\alpha(\lambda)}_{j_{\lambda}},\qquad f^{i,J}\not\equiv0,
\\
\intertext{we set the respective components of $\bs^{\dagger}$ equal 
to the sum of formal 
\emph{re}ver\-ses for each nonzero, homogeneous word in $\bs$,}
a_i^{\dagger} &\mathrel{{:}{=}}
 s_i^{\dagger}(\bx,\bvx^{\pm1})=
\sum_J \frac{1}{f^{i,J}(\bx)}\,\vx^{\,-\alpha(\lambda)}_{j_{\lambda}}
\circ\ldots\circ\vx^{\,-\alpha(1)}_{j_1},\label{EqSetACrossAtX}
\end{align}
\end{subequations}
where, at every 
point~$\bx\in U\subseteq M^n$, the sum is taken over the indexes~$J$ such that
the coefficients~$f^{i,J}$ 
do not vanish.
\footnote{In view of Remark~\ref{RemBVOddSection}, 
the fact that the extension~$\bs^{\dagger}$ remains undefined at the
zero locus of all these coefficients makes no harm.}
\end{convent}

\begin{example}
If 
\begin{equation}\label{EqReverseLoop}
a^i=\sum\nolimits_{k\in\BBZ}(\text{loop})^k,
\quad\text{ then }\quad
a_i^{\dagger}=\sum\nolimits_{k\in\BBZ}(\text{loop})^{-k},
\end{equation}
that is, all the reiterations of a closed path are walked backwards.
\end{example}

\begin{rem}\label{RemConventionTrueReverse}
Convention~\eqref{EqSetACrossAtX} means that, whenever each component~$s^i$ of the map~$\bs$ is just a single word, the 
respective dual~$a^{\dagger}_i$ becomes the weighted \emph{re}ver\-se --~and true \emph{in}ver\-se~-- of the path~$a^i(\bx,\bvx^{\pm1})$.
\end{rem}



\begin{rem}
When cyclic words $(a)$ are evaluated using~\eqref{EqSuperSection} at points~$\bx\in M^n$, each resulting cyclic word from the algebra~$\cX(\bvx^{\,\pm1})$ acquires an overall coefficient (which is supposed to be a smooth function on~$M^n$). The associativity of multiplication~$\cdot$ of the coefficients~$f_J(\bx)$ is used here. Note however that the commutativity of~$\cdot$ can be relaxed, yet if so, the result of evaluation $(a){\bigr|}_{(\bs,\bs^\dagger)}(\bx)\in\cX(\bvx^{\,\pm1})$ would depend on a position of the lock~$\binfty$ between letters of the word~$a$, see~\eqref{EqManyValues} and Remark~\ref{RemMergeModels} on p.~\pageref{RemMergeModels}.
\end{rem}

\subsection{Elementary (non)commutative variations}\label{SecElemNCVar} 
The precedence $\vec{e}_1\prec\ldots\prec\vec{e}_m\prec\vec{e}^{{}\,\dagger,1}\prec\ldots\prec\vec{e}^{{}\,\dagger,m}$
of the basic vectors for virtual shifts endows the Cartesian sum 
$T_{\ba}\Span(a^1,\dots,a^m)\oplus 
T_{\ba^{\dagger}}\Span(a^{\dagger}_1,\dots,a^{\dagger}_m)$
of the 
{dual} spaces with an orientation; it fixes the signs in all the structures of (non)commutative symplectic
geometry. The signs show up through the two couplings
$T_{\ba}V\times T_{\ba^{\dagger}}V^{\dagger}\to\Bbbk$ and $T_{\ba^{\dagger}}V^{\dagger}\times T_{\ba}V\to\Bbbk$
(which we denote by $\langle\,,\,\rangle$ in both cases, making no confusion; for the sequential order is essential).
Namely, we have that
\begin{equation}\label{EqTwoCouplings}
\bigl\langle\underrightarrow{\vec{e}_i,\vec{e}^{{}\,\dagger,j}}\bigr\rangle=\boldsymbol{\delta}_i{}^j
\quad\text{ and }\quad
\bigl\langle\underrightarrow{\vec{e}^{{}\,\dagger,j},\vec{e}_i}\bigr\rangle=-\boldsymbol{\delta}_i{}^j,
\end{equation}
where $\boldsymbol{\delta}_i{}^j$ is the Kronecker symbol that equals unit iff $i=j$ and which is set equal to zero
otherwise, see~\cite[\S2.2]{gvbv}.

Note that the virtual deformations $\delta\ba=\delta a^i(\bx)\cdot\vec{e}_i(\bx)$ and
$\delta\ba^{\dagger}=\delta a^{\dagger}_j(\bx)\cdot\vec{e}^{{}\,\dagger,j}(\bx)$ can be dependent on $\bx\in M^n$~---
and they should be such. By construction, each pair $(\delta\ba$,\ $\delta\ba^\dagger)$ of virtual shifts for the generators~$a^i$ and~$a^\dagger_i$ is a map belonging to the space $\Map\bigl(M^n\to T_{(\ba,\ba^\dagger)}\text{span}_\Bbbk(\ba;\ba^\dagger)\bigr)$.
We let the shifts be \emph{normalised} at all internal points $\bx\in\supp(\delta a^i)\subseteq M^n$
by the constraint
$$\delta a^i(\bx)\cdot\delta a_i^{\dagger}(\bx)\equiv1.\qquad
\lefteqn{\text{(no summation!)}}$$
This is why the couplings of virtual deformations are invisible in the 
ready\/-\/to\/-\/use formulae. Indeed, it is enough to know the signs
\begin{subequations}\label{EqNormalizeShifts}
\begin{align}
\left.\langle\delta a^i(\bx)\cdot\underrightarrow{\smash{\stackrel{\text{first}}{\vec{e}_i(\bx)}},\ %
\smash{\stackrel{\text{second}}{\vec{e}^{{}\,\dagger,i}(\bby)}}}\cdot\delta a_i^{\dagger}(\bby)\rangle\right|_{\bx=\bby}&=+1\label{EqNormalizeShiftsPlus}\\
\intertext{and}
\left.\langle\delta a^{\dagger}_i(\bby)\cdot\underrightarrow{\smash{\stackrel{\text{first}}{\vec{e}^{{}\,\dagger,i}(\bby)}},\ %
\smash{\stackrel{\text{second}}{\vec{e}_i(\bx)}}}\cdot\delta a^i(\bx)
\rangle\right|_{\bx=\bby}&=-1\label{EqNormalizeShiftsMinus},
\end{align}
\end{subequations}
at all the internal points $\bx$ of the support $\supp(\delta a^i)$, 
see
~\cite{gvbv,sqs13,dq15} for illustrations.%
\footnote{The usefullness of carrying the coefficients $\delta\ba(\,\cdot\,)$ and $\delta\ba^{\dagger}(\,\cdot\,)$ all way long 
is revealed in the geometry of \emph{iterated} variations;
let us also remember that we shall not always indicate the reference 
of frames $\vec{e}_i(\,\cdot\,)$ and $\vec{e}^{{}\,\dagger,i}(\,\cdot\,)$ to points of the substrate manifold~$M^n$. However, 
the fact that such reference 
is not impossible is crucial for the consistency of the formalism.}

\subsection{Parity\/-\/odd neighbours $\boldb=\Pi(\ba^{\dagger})$}\label{SecNeighbour}
From now on, let the set\/-\/up be $\BBZ_2$-\/graded by the function $|\,\cdot\,|$ that takes values in $\BBZ$ and determines the parity $(-)^{|\,\cdot\,|}$. 
All the objects which have been considered in the preceding sections were parity-even, of proper grading~$0$.
Let us relay the parity of symbols $a_i^{\dagger}$ by postulating that the new parity-odd variables carry the grading $+1$
(or \emph{minus} one, or any other (un)conventional odd integer number). To keep track of the reversed parity, let us denote these generators by $\boldb=(b_1,\,\dots\,,b_m)$ so that $\Pi\colon a_i^{\dagger}\rightleftarrows b_i$.

In the cyclic world, 
the concept of $\BBZ_2$-\/grading works as follows:\footnote{\label{FootAllTurnAround}%
Let $F$ be a homogeneous word of grading~$|F|$, written by using the $\BBZ_2$-\/graded alphabet. A full turn $F\mapsto \mathfrak{t}^{\lambda(F)}\,(F)$ along the orientation on the circle that carries the cyclic word~$F$ of length~$\lambda(F)$ yields the sign factor $(-)^{|F|\cdot(|F|-1)}=(+)$; the equality is valid because the product of two consecutive integers standing in the exponent is always even.
This argument shows also 
that, for a cyclic word to be rotated from a given configuration (determined by the position of the lock $\pmb{\infty}$ in between the word's letters) to another one, a choice to direct that rotation 
(coun\-ter)\/clock\-wise does not matter. Indeed, every clockwise rotation can be realised via one full turn clockwise (that would leave no effect by the above) followed by the appropriate shift backwards, in the counterclockwise direction.}
\begin{equation}\label{EqGradedShift}
\gotht\,(\gamma_1\circ\ldots\circ\gamma_{\lambda})=
(-)^{|\gamma_1\circ\ldots\circ\gamma_{\lambda-1}|\cdot|\gamma_{\lambda}|}
\gamma_{\lambda}\circ\gamma_1\circ\ldots\circ\gamma_{\lambda-1}.
\end{equation}
We denote by~$\cA^{(0|1)}$ the graded commutative unital 
non\/-\/associative algebra of cyclic words written in the alphabet
$1
$,\ $a^1$,\ $\dots$,\ $a^m$,\ $b_1$,\ $\dots$,~$b_m$.
By introducing the notation~$\cAZO$ we stress that the su\-per\-di\-men\-sion, equal to $(m|m)$, is positive in both the parity\/-\/even and odd components of the generators space~$\text{span}_\Bbbk(\ba;\boldb)$.


\begin{rem}[``$(abab)=0$\,?"]\label{RemABAB}
The idea that cyclic words acquire and accumulate the extra sign factors, whenever a parity\/-\/odd symbol overtakes the rest of the word, 
creates the following subtlety.

Set $m=1$ for definition and, omitting the symbols~$\circ$ of associative multiplication, first consider the cyclic word $(abaab)$. The identical, parity\/-\/odd letters~$b$ contained in it
can be distinguished nevertheless: one of them is \emph{followed} by $aa$ but preceded only by~$a$, whereas the other is \emph{preceded} by~$aa$ and followed by just a single copy of letter~$a$; we have that $(abaab)\sim-(aabab)$.

On the other hand, the cyclic word $(abab)$ does not contain any mechanism to distinguish between the two parity\/-\/odd entries~$b$, yet $(\underline{ab}ab)\sim-(ab\underline{ab})$ by construction. In fact, this word is synonymic to zero in the algebra
of cyclic words which are written in the parity\/-\/extended alphabet.\footnote{%
Analogous notions of \emph{zero} non\/-\/oriented graphs equipped with edge ordering and of \emph{zero} oriented graphs with an ordering of outgoing edges at every vertex are known from~\cite{Ascona96} and~\cite{Ascona96,KontsevichFormality}, respectively (cf.\ \cite{f16,cpp,JNMP2017} for illustrations).%
}
Let us be aware of the existence of this class of synonyms for zero; 
the calculus of iterated variations which we presently develop 
is indifferent to these synonyms 
existence. 
\end{rem}



\begin{model}[
The BV-\/geometry]\label{ModelBV}
We take the algebra $\cAZO$ as fibre%
\footnote{Note that the parity reversion~$\Pi$ does not modify the topology of spaces, whence conventions~\eqref{EqNormalizeShifts} remain valid for the virtual variations
$\delta\boldb=\delta b_i(\bx)\cdot\vec{e}^{{}\,\dagger,i}(\bx)$.
Note also that the presence of grading
does not modify our earlier 
convention~\eqref{EqSetACrossAtX} for the evaluation of symbols~--- as soon as a calculation governed by such graded arithmetic rule is over.} 
in the noncommutative superbundle~$\bpiNCZO$ over the sheaf~$\MnnC$.
This picture 
is summarised in Fig.~\ref{FigSetup}$(a)$, in which one easily recognises the noncommutative generalisation of the
classical Batalin\/--\/Vilkovisky geometry (see Fig.~\ref{FigSetup}$(b)$).%
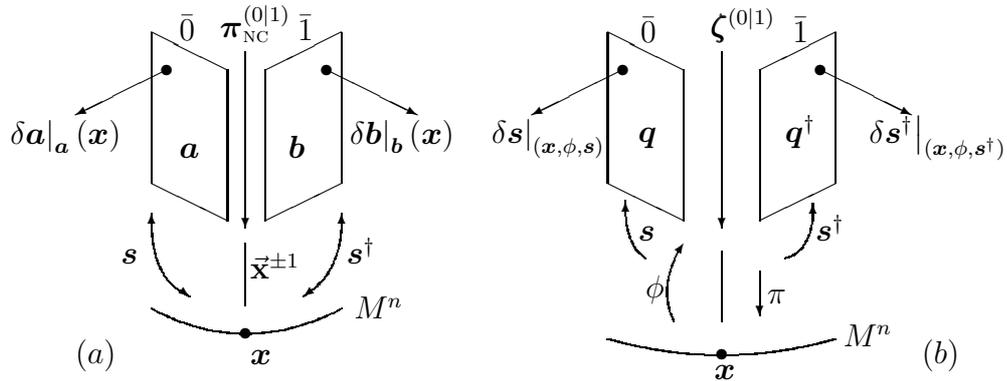
\begin{figure}[htb]
\begin{center}
\unitlength=1mm
\linethickness{0.4pt}
\begin{picture}(116.00,45.00)
\put(15.00,25.00){\line(2,-1){10.00}}
\put(25.00,20.00){\line(0,1){20.00}}
\put(25.00,40.00){\line(-2,1){10.00}}
\put(15.00,45.00){\line(0,-1){20.00}}
\put(30.00,20.00){\line(2,1){10.00}}
\put(40.00,25.00){\line(0,1){20.00}}
\put(40.00,45.00){\line(-2,-1){10.00}}
\put(30.00,40.00){\line(0,-1){20.00}}
\put(17.00,40.00){\circle*{1.33}}
\put(38.00,40.00){\circle*{1.33}}
\put(17.00,40.00){\vector(-2,-1){12.00}}
\put(38.00,40.00){\vector(2,-1){12.00}}
\put(27.33,8.67){\line(0,1){8.33}}
\bezier{112}(15.00,8.33)(27.33,1.67)(40.00,8.33)
\put(27.33,5.00){\circle*{1.33}}
\bezier{52}(35.00,10.00)(40.00,13.00)(40.00,20.00)
\put(35,10){\vector(-4,-3){0}}
\bezier{52}(20.00,10.00)(15.00,12.00)(15.00,20.00)
\put(20,10){\vector(4,-3){0}}
\put(15.00,19.00){\vector(0,1){2.00}}
\put(40.00,19.00){\vector(0,1){2.00}}
\put(18.67,44.00){\makebox(0,0)[lb]{$\bar0$}}
\put(34.25,44.00){\makebox(0,0)[lb]{$\bar1$}}
\put(18.67,28.33){\makebox(0,0)[lb]{$\ba$}}
\put(33.00,28.33){\makebox(0,0)[lb]{$\boldb$}}
\put(-4.00,29){\makebox(0,0)[lb]{$\left.\delta\ba\right|_{\ba}(\bx)$}}
\put(40.67,29){\makebox(0,0)[lb]{$\left.\delta\boldb\right|_{\boldb}(\bx)$}}
\put(11.00,14.33){\makebox(0,0)[lb]{$\bs$}}
\put(28.00,12.33){\makebox(0,0)[lb]{$\bvx^{\pm1}$}}
\put(40.67,14.33){\makebox(0,0)[lb]{$\bs^{\dagger}$}}
\put(41.67,7.00){\makebox(0,0)[lb]{$M^n$}}
\put(28.00,1.00){\makebox(0,0)[lb]{$\bx$}}
\put(5.00,0.00){\makebox(0,0)[lb]{$(a)$}}
\put(24,43.67){\makebox(0,0)[lb]{$\bpiNCZO$}}
\put(27.33,42.33){\vector(0,-1){23.33}}
\put(75.00,25.00){\line(2,-1){10.00}}
\put(85.00,20.00){\line(0,1){20.00}}
\put(85.00,40.00){\line(-2,1){10.00}}
\put(75.00,45.00){\line(0,-1){20.00}}
\put(95.00,40.00){\line(0,-1){20.00}}
\put(95.00,20.00){\line(2,1){10.00}}
\put(105.00,25.00){\line(0,1){20.00}}
\put(105.00,45.00){\line(-2,-1){10.00}}
\put(90.00,42.33){\vector(0,-1){23.33}}
\put(77.00,40.00){\vector(-2,-1){12.00}}
\put(103.00,40.00){\vector(2,-1){12.00}}
\put(103.00,40.00){\circle*{1.33}}
\put(77.00,40.00){\circle*{1.33}}
%
\put(90,6.67){\line(0,1){9.33}}   
\bezier{124}(75.00,4.33)(90.00,0.00)(105.00,4.33)
\put(90.00,2.33){\circle*{1.33}}
\put(95.00,13.33){\vector(0,-1){6.00}}
\bezier{40}(83.67,6.67)(81.67,11.00)(84.33,16.00)
\put(83.00,14){\vector(2,3){2}}
\bezier{32}(80.00,15.00)(77.33,16.67)(77.33,22.00)
\put(77.33,21.00){\vector(0,1){1.67}}
\bezier{36}(98.33,15.00)(102.00,16.00)(102.00,21.67)
\put(102.00,20.00){\vector(0,1){2.33}}
\put(88.67,43.67){\makebox(0,0)[lb]{$\pmb{\zeta}^{(0|1)}$}}
\put(79.33,43.67){\makebox(0,0)[lb]{$\bar0$}}
\put(99.33,43.67){\makebox(0,0)[lb]{$\bar1$}}
\put(98.33,29.67){\makebox(0,0)[lb]{$\bq^{\dagger}$}}
\put(79.00,29.67){\makebox(0,0)[lb]{$\bq$}}
   \put(59.33,28.33){\makebox(0,0)[lb]{$\left.\delta\bs\right|_{(\bx,\phi,\bs)}$}}
\put(109.00,28){\makebox(0,0)[lb]{$\left.\delta\bs^{\dagger}\right|_{(\bx,\phi,\bs^{\dagger})}$}}
\put(79.00,17.33){\makebox(0,0)[lb]{$\bs$}}
\put(102.33,17.33){\makebox(0,0)[lb]{$\bs^{\dagger}$}}
\put(96.00,9.00){\makebox(0,0)[lb]{$\pi$}}
\put(80.00,9.00){\makebox(0,0)[lb]{$\phi$}}
\put(106.00,3.33){\makebox(0,0)[lb]{$M^n$}}
\put(89.00,-1.00){\makebox(0,0)[lb]{$\bx$}}
\put(116.33,0.00){\makebox(0,0)[lb]{$(b)$}}
\end{picture}
\caption{The elementary displacements~$\bvx^{\pm1}$ in a tiling of~$M^n$ versus the gauge connection fields~$\phi$ over the space\/-\/time~$M^n$;
the canonical duality of diagonal variations for the opposite\/-\/parity halves of the alphabet versus the opposite\/-\/parity field\/-\/antifield and ghost\/-\/antighost pairs.}\label{FigSetup}
\end{center}
\end{figure}
The r\^ole of physical fields~$\phi$ as sections of their bundle~$\pi$ is now played by the primitive displacements~$\bvx^{\pm1}$ in granulated space. 
The fibre algebra generated by the symbols~$a^i$ and~$b_i$ was known to us before as the Whitney sum of parity\/-\/even and odd components in the 
Batalin\/--\/Vilkovisky superbundle~$\boldsymbol{\zeta}^{(0|1)}$, pulled back --\,by the projection~$\pi$\,-- over the total space of the bundle of physical fields. The symbols~$\ba$ and $\boldb=\Pi(\ba^{\dagger})$ of opposite parities form the
noncommutative analogue of the BV-\/zoo $\bq$, $\bq^{\dagger}$ inhabited by the (anti)\/fields and (anti)\/ghosts. The r\^ole of
the BV-\/bundle sections is granted 
to the two maps~$\bs$ and~$\bs^{\dagger}$.%
\footnote{We recall from~\cite{gvbv} that the normalised variations $\delta\bs$ and $\delta\bs^{\dagger}$ were the dual
components in sections of the tangent bundle $\BBT\boldsymbol{\zeta}^{(0|1)}$;
the vectors $\delta\bs\left(\bx,\phi(\bx),\bs(\bx,\phi(\bx))\right)$ and
$\delta\bs^{\dagger}\left(\bx,\phi(\bx),\bs^{\dagger}(\bx,\phi(\bx))\right)$
were attached at points of graphs of sections for the BV-superbundle induced over $\pi$. The construction of these test
shifts was laborious 
indeed in the graded\/-\/commutative world. On the other hand, the noncommutative target spaces contain nothing else but the
basic letters $\ba$ and $\boldb$ that undergo the virtual deformations, so that the picture is simplified considerably.}
\end{model}

\begin{rem}
It will readily be seen that both the Batalin\/--\/Vilkovisky Laplacian of the integral functionals given by \emph{zero} words --\,or the Schouten bracket taken for zero word functionals with any other cyclic\/-\/word functional\,-- vanish identically.
\end{rem}

\begin{rem}\label{RemMergeModels}
Models~\ref{ModelMatr} and~\ref{ModelStar}, as well as Models~\ref{ModelBV} and~\ref{ModelStar} can be combined. 
For instance (for the latter pair), the (quasi)\/crystal tiling of an affine manifold~$M^n$ yields the alphabet~$\bvx^{\,\pm1}$ and 
concatenation~$\circ$ in the algebra of formal paths that show up in~\eqref{EqSuperSection}, whereas a given Poisson structure on that manifold~$M^n$ yields the associative $\star$-\/product which is used to multiply the coefficients~$f_J(\bx;\hbar)\in C^\infty(M^n)((\hbar))$ occurring in~\eqref{EqSuperSection}.

However, it is the graded\/-\/commutative model over the sheaves~$\MnnC$ of algebras of walks along a tiling of~$M^n$ which will be the default set\/-\/up in the further study.
\end{rem}


\subsection{The ring of noncommutative local functionals}\label{SecLocalFunctionals}
Let us proceed from \emph{functions} on the space $J^\infty\bigl(\bpiNCZO 
\bigr)$ of jets of sections~\eqref{EqSuperSection} to the notion of \emph{functionals} that take the evaluation mappings~$(\bs$,\ $\bs^\dagger)$
to formal cyclic words\footnote{Such cyclic words are \emph{formal} because (\textit{i}) they could encode no realisable paths along the edges of the graph and (\textit{ii}), although ``cyclic'' by construction, each homogeneous component of such words could not encode a \emph{closed} walk, even if it did specify \emph{some} walk along the edges.}
written in the alphabet~$\bvx^{\pm1}$ of edges in the adjacency graph for a given crystal tiling of the substrate manifold~$M^n$.

\begin{convent}\label{RemFatInfty}
On the infinite jet space $J^\infty\bigl(\bpiNCZO 
\bigr)$,
every 
cyclic word~$(f)$ is a sum of its homogeneous components, each weighted by the 
coefficients that (can) depend on points~$\bx$ of the substrate manifold~$M^n$. For the sake of definition, let us assume that every such coefficient is $C^\infty$-\/smooth on~$M^n$; their asymptotic behaviour must also be specified in advance so that the integration by parts makes sense.
Specifically, if the manifold~$M^n$ is closed, then there is nothing to discuss: the empty boundary carries no boundary terms.
However, should there be one, $\dd M^n\neq\varnothing$, or should the manifold~$M^n$ be non\/-\/compact (e.g., let $M^n=\BBR^n$ with the standard Euclidean topology), then we postulate that the coefficients decay rapidly towards the boundary~$\dd M^n$ or spatial infinity, respectively.

Likewise, we suppose that the supports $\supp\delta a^i$ of the $C^\infty(M^n)$-\/smooth infinitesimal variations $\delta a^i(\cdot)\cdot\vec{e}_i(\cdot)\colon M^n\to T_\ba\spanOp(a^1,\ldots,a^m)$ are compact
and $\supp\delta a^i \cap \partial M^n=\varnothing$.
\end{convent}


The \emph{volume element} $\dvol(\bx)$ on~$M^n$ in the construction of integral functionals over the jet space $J^\infty\bigl(\bpiNCZO 
\bigr)$ is another piece of external data.

\begin{convent}\label{RemVolumeElement}
We suppose that a volume element~$\dvol(\bx)$ is given at all points~$\bx\in M^n$ (possibly, in a way that depends on the tiling at hand). Also, we technically assume in this text 
that the volume element~$\dvol(\bx)$ may not depend on a choice of the mappings~$(\bs,\bs^\dagger)$~--- that is, in a sense, on a configuration of noncommutative ``fields'' over the granulation~$\MnnC$ of the physical space~$M^n$.

One could think that the volume element $\dvol(\cdot)$ is placed in the locks of cyclic words; this idea is practical because, whenever any such word is unlocked, it is converted at once into a singular linear integral operator supported on the diagonal;
the volume element then disappears, giving way to the attachment points' congruence mechanism through the locality of couplings~\eqref{EqTwoCouplings} in~\eqref{EqNormalizeShifts}.
\end{convent}

\begin{convent}
From now on we restrict the study to the class of functionals such that densities of the integral functionals 
$F=\int f\bigl(\bx$,\ 
$[\ba]$,\ $[\boldb]\bigr)\circ\dvol(\bx)$
do not depend explicitly on the edge alphabet~$\bvx^{\,\pm1}$ (but can do so implicitly through a differential dependence of densities on~$\ba$ or~$\boldb$, which are evaluated at the jets $j^\infty_{\bx}(\bs,\bs^\dagger)$ of sections~\eqref{EqSuperSection} for~$\bpiNCZO$. 
(We recall that such vertical subtheory makes the full theory in Models~\ref{ModelMatr} and~\ref{ModelStar}, cf.\ footnote~\ref{FootWhyNoEdgesInDensity} on p.~\pageref{FootWhyNoEdgesInDensity}.)
\end{convent}

\begin{notation}
The vector space of such integral functionals 
will be denoted by~$\bar{H}^n\bigl(\bpiNCZO 
\bigr)$. 
\end{notation}

Integral functionals $F_1$,\ $\ldots$,\ $F_\ell\in\bar{H}^n\bigl(\bpiNCZO 
\bigr)$ are the building blocks in the \emph{local} functionals such as $F_1\times\ldots\times F_\ell\in{\bar{H}^{n^{\otimes\ell}}}\bigl(\bpiNCZO 
\bigr)$.

\begin{define}\label{DefProduct}
Let $F_1=\int f_1\bigl(\bx_1,[\ba],[\boldb]\bigr)\circ\dvol(\bx_1)$ and
$F_2=\int f_2\bigl(\bx_2,[\ba],[\boldb]\bigr)\circ\dvol(\bx_2)$ be two linear integral functionals the densities of which do not depend explicitly on any letters from the edge alphabet~$\bvx^{\pm1}$. The \emph{product}
\[
F_1\times F_2 = \iint (f_1){\bigr|}_{(\bx_1,[\ba],[\boldb])} \times
(f_2){\bigr|}_{(\bx_2,[\ba],[\boldb])} \circ\dvol(\bx_1)\cdot\dvol(\bx_2)
\in{{\bar{H}^{n^{\otimes 2}}}}\bigl(\bpiNCZO 
\bigr)
\]
is the horizontal cohomology class of linear integral functionals over $\bigl({M^n}^{\otimes 2}$,\ $\dvol(\,\cdot\,)^{\otimes 2}\bigr)$ such that their densities are equivalent to the product~$(f_1)\times(f_2)$ in~$\cA^{(0|1)}$.

Setting ${\bar{H}^{n^{\otimes 0}}}\bigl(\bpiNCZO 
\bigr)$ equal to $\Bbbk\cdot(\BOne)$ by definition, we extend the bi\/-\/linear operation~$\times$ recursively from pairs of integral functionals to the multiplication of products of any nonnegative number of functionals.
Because the operation~$\times$ is not associative, there are the respective Catalan number ways to arrange the multiplications in $F_1\times\ldots\times F_\ell$ by inserting the $\ell-1$ balanced pairs of parentheses. We let the \emph{default ordering} be lexicographic: $(\cdots(F_1\times F_2)\times\ldots\times F_{\ell-1})\times F_\ell$. 
\end{define}

\begin{cor}
The multiplication~$\times$ of local functionals over~$J^\infty\bigl(\bpiNCZO 
\bigr)$ is graded\/-\/commutative: $F\times G=(-)^{|F|\cdot|G|}G\times F$ for~$F$ and~$G$ homogeneous.
\end{cor}

\begin{notation}
Denote by 
\begin{equation}\label{EqLocal}
\overline{\mathfrak{M}}^n\bigl(\bpiNCZO 
\bigr) =
\bigoplus_{\ell\geqslant0}{\bar{H}^{n^{\otimes\ell}}}\bigl(\bpiNCZO 
\bigr)
\end{equation}
the $\BBZ_2$-\/graded commutative non\/-\/associative unital ring of local functionals in the noncommutative set\/-\/up under study.
\end{notation}

To define the value of a local functional~$F$ at a section~$(\bs,\bs^\dagger)$, first let us consider 
the class of \emph{integral} functionals such as 
$F=\int f\bigl(\bx$,\ 
$[\ba]$,\ $[\boldb]\bigr)\circ\dvol(\bx)$, where the cyclic word~$(f\circ\dvol(\bx))$ marks an 
equivalence class modulo integrations by parts (no boundary terms!\,). 

\begin{define}
The value of such integral functional at a given mapping~$(\bs,\bs^\dagger)$ is
\begin{equation}\label{EqFs}
F(\bs,\bs^\dagger)\eqdef\int_{M^n} f\bigl(\bx, 
\jet^\infty_{\bx}(\bs),\jet^\infty_{\bx}(\bs^\dagger)\bigr)\circ\dvol(\bx)\in\cX(\bvx^{\pm1});
\end{equation}
the integral makes sense due to our earlier assumptions on the global choice of alphabet~$\bvx^{\pm1}$ on the entire~$M^n$ (that is, the tiling $M^n=\bigcup_\alpha \overline{\Delta}_\alpha$ is not \emph{quasi}crystal) and on the class of functional coefficients depending on~$\bx$, so that the (im)\/proper integral converges. 

The evaluation of products $F_1\times\ldots\times F_\ell$ of functionals at a given mapping~$(\bs,\bs^\dagger)$ goes as follows; 
without loss of generality suppose~$\ell=2$.
First, double $(\bs,\bs^\dagger)\mapsto(\bs,\bs^\dagger)^{\otimes2}$ for the 
$\ell=2$ copies of the substrate manifold~$M^n$, and then integrate over
${M^n}^{\otimes2}$ in the element of~${\bar{H}^{n^{\otimes 2}}}\bigl(\bpiNCZO
\bigr)$.
\end{define}

\begin{rem}\label{RemTwoMult}
Through the evaluation procedure, local functionals
keep track of the fibre algebra~$\cAZO$ of cyclic words (even though neither the letters $a^i$ nor $b_j$ show up in the functionals' values that belong to the functionals value space~$\cX(\bvx^{\pm1})$ of cyclic words written in the edge alphabet~$\bvx^{\pm1}$).

Indeed, we recall from Remark~\ref{RemNotAlgHom} that generally speaking,
\[
\bigl(F_1\mathop{\stackrel{\cA^{(0|1)}}{\times}}F_2\bigr)\,(\bs,\bs^\dagger)
\neq
F_1(\bs,\bs^\dagger) \mathop{\stackrel{\cX(\bvx^{\pm1})}{\times}}
F_2(\bs,\bs^\dagger).
\]
Moreover, 
although the multiplication~$\times$ in~$\cAZO$ is~$\BBZ_2$-\/graded commutative, that grading is lost in the course of functionals' evaluation at the mappings~$(\bs,\bs^\dagger)$;
the multiplication~$\times$ in the non\/-\/graded algebra~$\cX(\bvx^{\pm1})$ is just commutative.\\
\centerline{\rule{1in}{0.7pt}}
\end{rem}

\noindent%
In the remaining part of this chapter we reveal the structure of differential (shifted-) graded Lie algebra
--\,more specifically
, the BV algebra\,-- on the $\BBZ_2$-\/graded commutative non\/-\/associative unital ring $\overline{\mathfrak{M}}^n\bigl(\bpiNCZO 
\bigr)$ of local functionals.
First we introduce some notation.
Let us recall that the generators~$a^i$ and~$b_i$ are evaluated at 
sections $(\bs,\bs^\dagger)$,
whereas the generator virtual shifts~$(\delta\ba,\delta\boldb)$
are taken from the space $\Map\bigl(M^n\to T_{(\ba,\boldb)}\text{span}_\Bbbk(\ba;\boldb)\bigr)$.
To permit the iteration of variations, one has to deal with the space of local functionals such that densities of their integral building blocks can contain not only the generators but also their shifts 
(see footnote~\ref{FootTpi} on p.~\pageref{FootTpi} below). 

\begin{notation}
In order to avoid an agglomeration of formulae, let us denote by
\[
\overline{\mathfrak{N}}^n\bigl(\BBT\bpiNCZO 
\bigr)=
\bigoplus_{\ell\geqslant0}{\bar{H}^{n^{\otimes\ell}}}\bigl(\BBT\bpiNCZO 
\bigr)
\]
the vector space of such local functionals over the jet space~$J^\infty\bigl(\bpiNCZO 
\bigr)$.
\end{notation}


\begin{rem}\label{RemWordWasFirst}
The multiplication of 
functionals, as part of the construction of space $\overline{\mathfrak{M}}^n\bigl(\bpiNCZO
\bigr)$, is provided by Definition~\ref{DefProduct}. The BV~Laplacian~$\Delta$ (see p.~\pageref{DefLaplace} below and Definition~\ref{DefSchouten} on p.~\pageref{DefSchouten}) is a local variational operator 
on the space of local functionals, hence every argument of~$\Delta$ is encoded by a cyclic word. This means that \emph{first} such argument is formed (if necessary, by using the structure~$\times$ of algebra~$\cAZO$ whenever that input object is a product of several integral functionals; parentheses would specify the consecutive order of non\/-\/associative multiplications). Secondly, the BV~algebra's differential operations~$\Delta$ or~$\lshad\,,\,\rshad$ rework the input into an element of~$\overline{\mathfrak{M}}^n\bigl(\bpiNCZO
\bigr)$. In particular, at no moment are any intermediate objects from~$\overline{\mathfrak{N}}^n\bigl(\BBT\bpiNCZO 
\bigr)$ multiplied anew by using the structure~$\times$ for~\eqref{EqLocal}.

For example, identity~\eqref{1a} on p.~\pageref{1a} below 
frames an application of the differential structure~$\lshad\,,\,\rshad$ to the functional $F\times(G\times H)\in\overline{\mathfrak{M}}^n\bigl(\bpiNCZO 
\bigr)$, referred to at least three copies of the underlying manifold~$M^n$. The same ordering --\,multiplication, then variation over~$M^n$\,-- applies to both sides of identity~\eqref{EqLapSchouten} where the BV~Laplacian~$\Delta$ works on the product~$F\times G$ twice (in particular, via~$\lshad\,,\,\rshad$ to which the operator~$\Delta$ is parent).

Therefore, let us remember that it is the ring $\overline{\mathfrak{M}}^n\bigl(\bpiNCZO
\bigr)$
but not the larger vector space $\overline{\mathfrak{N}}^n\bigl(\BBT\bpiNCZO
\bigr)$ on which the BV~al\-ge\-b\-ra structure is well defined. The reduction from $\overline{\mathfrak{N}}^n\bigl(\BBT\bpiNCZO 
\bigr)$ to $\overline{\mathfrak{M}}^n\bigl(\bpiNCZO 
\bigr)$ amounts to a perfect matching and then coupling of the (co)\/vectors~$\vec{e}_i$ and~$\vec{e}^{\,\dagger,i}$ in all pairs of canonically dual components~$\delta\ba$ and~$\delta\boldb$ of the variations. In the next section we recall the geometric mechanism of integration by parts; the way \emph{how} the couplings are reconfigured itself is the algorithmic definition of the BV~al\-ge\-b\-ra structure (see section~\ref{SecLaplacian}).
\end{rem}

\subsection{Elements of the geometric theory of variations}\label{SecElementsGVBV}
The Gel'fand framework of singular integral distributions is known, e.g., from~\cite{GelfandShilov}. In our case, the space $\overline{\mathfrak{N}}^n\bigl(\BBT\bpiNCZO
\bigr)$ of local functionals over the tangent superbundle~$\BBT\bpiNCZO$
is the space of basic objects on which the variations act by singular linear integral operators. 

For consistency, let us outline 
key ideas in the geometry of iterated variations (introduced in~\cite{gvbv,sqs13} and illustrated in~\cite{dq15,prg15});
%
they are as follows.
\begin{itemize}
\item
The unlinking of a cyclic word, together with an intention to paste the open string of symbols contained in it into another word as an uninterrupted fragment, converts the (procedure of) insertion of that string into a singular linear integral operator supported on the diagonal.
\item
Such operators are singular because the restriction to the diagonal over points in copies of the substrate manifold~$M^n$ is ensured by 
ordered couplings~\eqref{EqTwoCouplings} which are \emph{not defined} off the diagonal $\bx=\bby$ in~\eqref{EqNormalizeShifts}.
\item
The definitions of the Batalin\/--\/Vilkovisky Laplacian~$\Delta$ and variational Schouten bracket $\lshad\,,\,\rshad$ are 
operational, that is, every such definition is an algorithm for the on\/-\/the\/-\/diagonal reconfiguration of the couplings.
\item
The objects that are usually viewed in the calculus of variations as differential forms are either the volume element $\dvol(\bx)$ on the substrate manifold~$M^n$
or the dual bases $\vec e_i$, $\vec e^{{}\,\dagger,i}$ in the tangent spaces attached at the point~$(\ba,\bb)$ of the fibre 
algebra (this is what its alphabet was doubled for). 
The \emph{orientation} uniquely determines the signs
of couplings~\eqref{EqTwoCouplings} by ordering the tangent vectors. This explains why such differential $1$-\/forms anticommute.
\end{itemize}
\begin{convent}
In the course of virtual variation of the symbols $a^i_{\sigma}$ and $b_{j,\tau}$ by using%
\footnote{It is readily seen that the congruence of multi-indexes $\sigma$ in $(\dd/\dd\bx)^{\sigma}$ and $a^i_{\sigma}$
(as well as in the partial derivative $\vec{\dd}/\dd a^i_{\sigma}$, see~\eqref{EqShiftAOp} below) refers to the definition of vector as an
equivalence class of curves passing through a point.}
\begin{equation}\label{EqShiftsJetVariables}
{(\delta a^i)\Bigl(\frac{\overleftarrow{\dd}}{\dd\bx}\Bigr)^{\sigma}(\bx)\cdot\vec e_i(\bx)}
\quad\text{ and }\quad
{(\delta b_j)\Bigl(\frac{\overleftarrow{\dd}}{\dd\bx}\Bigr)^{\tau}(\bx)\cdot\vec e^{{}\,\dagger,j}(\bx)},
\end{equation}
the responses of integral functionals are always expanded with respect to the 
dual bases~$\vec e^{{}\,\dagger,i}$ and~$\vec e_i$. 
For instance, we obtain the singular linear integral operators
\begin{subequations}\label{EqShiftABOp}
\begin{align}
\overrightarrow{\delta\ba}(\cdot)&=\int_{M^n}\Id\bby\sum_{i=1}^m\sum_{|\sigma|\geqslant0}(\delta a^i)
\Bigl(\frac{\overleftarrow{\dd}}{\dd\bby}\Bigr)^{\sigma}(\bby)\cdot
\bigl\langle
\underrightarrow{\stackrel{\text{first}}{\vec e_i(\bby)},\stackrel{\text{second}}{\vec e^{{}\,\dagger,i}(\cdot)}}
\bigr\rangle
\frac{\overrightarrow{\dd}}{\dd a^i_{\sigma}}\label{EqShiftAOp}
\\
\intertext{and}
\overrightarrow{\delta\bb}(\cdot)&=\int_{M^n}\Id\bz\sum_{j=1}^m\sum_{|\tau|\geqslant0}(\delta b_j)
\Bigl(\frac{\overleftarrow{\dd}}{\dd\bz}\Bigr)^{\tau}(\bz)\cdot
\bigl\langle
\underrightarrow{\stackrel{\text{first}}{(-\vec e^{{}\,\dagger,j})}(\bz),\stackrel{\text{second}}{\vec e_i(\cdot)}}
\bigr\rangle
\frac{\overrightarrow{\dd}}{\dd b_{j,\tau}}\ .\label{EqShiftBOp}
\end{align}
\end{subequations}
This convention will be illustrated in the sequel.
\end{convent}
\begin{itemize}
\item
Given by its own singular integral operator, each variation brings a new copy of the integration domain~$M^n$ into the
picture. In consequence, all the intermediate objects 
$\text{Obj}\in\overline{\mathfrak{N}}^n\bigl(\BBT\bpiNCZO 
\bigr)$
that emerge in the course of calculations do retain a kind of memory
of the way how they were obtained from the input data.%
\footnote{\label{FootTpi}%
In the (graded-)commutative language of bundles this means that their products
$\pmb{\zeta}^{(0|1)}\times\BBT\pmb{\zeta}^{(0|1)}\times\ldots\times\BBT\pmb{\zeta}^{(0|1)}$,
standing over 
$M^n\times M^n\times\ldots\times M^n$,
are taken, but not their Whitney sums
$\pmb{\zeta}^{(0|1)}\times_{M^n}\BBT\pmb{\zeta}^{(0|1)}\times_{M^n}\ldots\times_{M^n}\BBT\pmb{\zeta}^{(0|1)}$
are fibred over a single copy of the base manifold $M^n$
.}\label{FootWhitney}
That is, no calculation can be interrupted along the way.
\end{itemize}

\begin{lemma}[Integration by parts]
From the powers of partial derivatives $\bigl(\overleftarrow{\dd}/\dd\bby\bigr)^\sigma$ that act on the test shifts in~\eqref{EqShiftABOp} one obtains, 
due to the locality of couplings $\langle{\cdot},{\cdot}\rangle$, 
the powers of minus total derivatives $\bigl(-\overrightarrow{\Id}/\Id\bx\bigr)^\sigma$ that act on densities of integral functionals.
\end{lemma}

\begin{proof}[Explanation \textup{(see~\cite{gvbv})}]
Consider a point~$\bby$ of the affine manifold~$M^n$ and denote by $\bby+\delta\bby\in M^n$ a near\/-\/by point with coordinates~$y^i+\delta y^i$, here and immediately below $1\leqslant i,\alpha\leqslant n$; the notation $\lim_{\delta\bby\to 0}$ makes obvious sense. For the sake of brevity, put $\sigma\mathrel{{:}{=}}\{x^\alpha\}$. We have that, due to the absence of boundary terms and then by definition (by Newton\/--\/Leibniz and in the last line, by S.\,Lie),
\begin{multline*}
\int\Id\bby \bigl\langle(\delta a^i)\frac{\overleftarrow{\dd}}{\dd y^\alpha}(\bby)\cdot
\underrightarrow{\vec{e}_i(\bby),\vec{e}^{\,\dagger,i}(\bx)}\,
\Bigl(\frac{\overrightarrow{\dd}}{\dd a^i_{x^\alpha}} f(\bx,[\ba],[\boldb])\Bigr)\Bigr|_{\jet^\infty_\bx(\bolds,\bolds^\dagger)}\bigr\rangle={}
\\
{}=
\int\Id\bby\ \delta a^i(\bby)\bigl(-\frac{\overrightarrow{\dd}}{\dd y^\alpha}\bigr)
\bigl\langle\underrightarrow{\vec{e}_i(\bby),\vec{e}^{\,\dagger,i}(\bx)}\,
\Bigl(\frac{\overrightarrow{\dd}}{\dd a^i_{x^\alpha}} f(\bx,[\ba],[\boldb])\Bigr)\Bigr|_{\jet^\infty_\bx(\bolds,\bolds^\dagger)}\bigr\rangle
\\
{}
\stackrel{\text{def}}{=} -\int\Id\bby\ \delta a^i(\bby)\,\lim\limits_{\delta y^\alpha\to+0}\frac{1}{\delta y^\alpha} 
\left\{
\begin{aligned}
\bigl\langle & \underbrace{
\vec{e}_i(\bby+\delta y^\alpha),\vec{e}^{\,\dagger,i}(\bx)
}_{\text{$+1$ if $\bx=\bby+\delta y^\alpha $}}\,
\frac{\overrightarrow{\dd}}{\dd a^i_{x^\alpha}} f(\bx,[\ba],[\boldb])\bigr|_{\jet^\infty_\bx(\bolds,\bolds^\dagger)}\bigr\rangle
\\
&
-\bigl\langle\underbrace{\vec{e}_i(\bby),\vec{e}^{\,\dagger,i}(\bx)}_{\text{$+1$ if $\bx=\bby$}}\,
\frac{\overrightarrow{\dd}}{\dd a^i_{x^\alpha}} f(\bx,[\ba],[\boldb])\bigr|_{\jet^\infty_\bx(\bolds,\bolds^\dagger)}\bigr\rangle
\end{aligned}
\right\}
\\
{}
\stackrel{\text{def}}{=} \int\Id\bby\ \delta a^i(\bby)
\bigl\langle\underrightarrow{\vec{e}_i(\bby),\vec{e}^{\,\dagger,i}(\bx)}\,
\bigl(-\frac{\overrightarrow{\dd}}{\dd x^\alpha}\bigr)
\Bigl( \frac{\overrightarrow{\dd}}{\dd a^i_{x^\alpha}} f(\bx,[\ba],[\boldb])\bigr|_{\jet^\infty_\bx(\bolds,\bolds^\dagger)} \Bigr) \bigr\rangle
\\
{}
\stackrel{\text{def}}{=} \int\Id\bby\ \delta a^i(\bby)
\langle\underrightarrow{\vec{e}_i(\bby),\vec{e}^{\,\dagger,i}(\bx)}\rangle\cdot
\Bigl(\bigl(-\frac{\overrightarrow{\Id}}{\Id x^\alpha}\bigr)
\frac{\overrightarrow{\dd}}{\dd a^i_{x^\alpha}} f(\bx,[\ba],[\boldb])\Bigr)\Bigr|_{\jet^\infty_\bx(\bolds,\bolds^\dagger)}.
\end{multline*}
For multi\/-\/indexes~$\sigma$ longer than~$\{x^\alpha\}$ the powers of partial derivatives $(\overleftarrow{\dd}/\dd\bby)^\sigma$ are processed by repeated integrations by parts; this yields the powers of minus total derivatives $(-\overrightarrow{\Id}/\Id\bx)^\sigma$. 
In the course of derivation of densities with respect to not~$a^i_\sigma$ but~$b_{j,\tau}$ and so, in the course of using the other of two (co)\/vectors couplings, all reasonings are still performed in exactly the same way.
\end{proof}
\begin{convent}
In every calculation, the integrations by parts are performed last, prior only to the reconfigurations of couplings and
their evaluation by using~\eqref{EqNormalizeShifts}. For instance, the derivative
$(\overleftarrow{\dd}/\dd\bby)^{\sigma}$ in 
formula~\eqref{EqShiftAOp} channels through $\vec e_i(\bby)$ and
$\vec e^{{}\,\dagger,i}(\bx)$ on the diagonal $\bby=\bx$ (which is the locus 
where the coupling is defined);
the derivative thus becomes $(-\vec{\Id}/\Id\bx)^{\sigma}$ that falls on (a derivative of) the argument's density at
$\bx\in M^n$.

This principle 
makes the variations (graded-)\/permutable.
\end{convent}

\begin{notation}
To keep track 
where the total derivatives would come from after integration by parts and to emphasize that such integrations are performed at the end of a calculation, we embrace the (powers of) minus the total derivatives by using the delimiters $\vphantom{\bigl|}^{\lceil}\ldots\vphantom{\bigr|}^{\rceil}$. 
Likewise, in the notation for those total derivatives we preserve the base variables from singular linear integral operators. (We remember that couplings~\eqref{EqNormalizeShifts} wright the diagonal, hence the above convention refers to notation only.) In these terms, operators~\eqref{EqShiftABOp} can be realised by using the formulas
\begin{align*}
\overrightarrow{\delta\ba}(\cdot)&=\int_{M^n}\Id\bby\sum_{i=1}^m\sum_{|\sigma|\geqslant0} \delta a^i (\bby)\cdot
\bigl\langle
\underrightarrow{\stackrel{\text{first}}{\vec e_i(\bby)},\stackrel{\text{second}}{\vec e^{{}\,\dagger,i}(\cdot)}}
\bigr\rangle
\,\vphantom{\Bigl|}^{\lceil}\Bigl(-\frac{\overrightarrow{\Id}}{\Id\bby}\Bigr)^{\sigma}\vphantom{\Bigr|}^{\rceil}
\frac{\overrightarrow{\dd}}{\dd a^i_{\sigma}}
\\
\intertext{and}
\overrightarrow{\delta\bb}(\cdot)&=\int_{M^n}\Id\bz\sum_{j=1}^m\sum_{|\tau|\geqslant0} \delta b_j (\bz)\cdot
\bigl\langle
\underrightarrow{\stackrel{\text{first}}{(-\vec e^{{}\,\dagger,j})}(\bz),\stackrel{\text{second}}{\vec e_i(\cdot)}}
\bigr\rangle
\,\vphantom{\Bigl|}^{\lceil}\Bigl(-\frac{\overrightarrow{\Id}}{\Id\bz}\Bigr)^{\tau}\vphantom{\Bigr|}^{\rceil}
\frac{\overrightarrow{\dd}}{\dd b_{j,\tau}}\ ,
\end{align*}
respectively.
\end{notation}
\begin{itemize}
\item
By construction, \emph{iterated} variations of a functional over a copy of~$M^n$ 
never spread from it to the fragments of other functionals in
any composite object during multiple integrations by parts over~$M^n$
(e.g., see~\cite{sqs13,prg15}).
\end{itemize}
Summarising, the BV~calculus of iterated variations relies heavily on a reference of each object to the copy of manifold~$M^n$ over which that object was defined; the locality of couplings~\eqref{EqNormalizeShifts} provides a 
restriction to the diagonal over all such copies at the end of the day. 
The association with 
own bases is the mechanism that discriminates between the fibre
letters from different words in the input. Indeed, 
integrations by parts over the words' substrates~$M^n$ act by total derivatives only on the letters from the respective words.

We refer to~\cite{gvbv,sqs13,dq15,prg15} for more details and illustrations of these guiding principles.%

\subsection{How the Batalin\/--\/Vilkovisky Laplacian determines the Schouten bracket}\label{SecLaplacian}
Now we are ready to outline the construction of parity-odd BV~
La\-pla\-ci\-an~$\Delta$. On the space of local functionals over
the jet space $J^{\infty}(\bpiNCZO 
)$, it is the parent structure for the noncommutative 
variational Schouten bracket~$\lshad\,,\,\rshad$. We establish the main properties of these structures, recalling further the relations between
them.

\begin{define}\label{DefLaplace}
The Batalin\/--\/Vilkovisky Laplacian is the reconfiguration --\,shown in 
Fig.~\ref{ReconfigCouplingsBV}\,--
of (co)vec\-tor couplings in the second variation 
$\overrightarrow{\delta\ba}(\overrightarrow{\delta\bb}(\cdot))$
of a local functional on the jet space $J^{\infty}(\bpiNCZO 
)$.
\end{define}
\begin{figure}[h]
\begin{center}
\unitlength=1mm
\special{em:linewidth 0.4pt}
\linethickness{0.4pt}
\begin{picture}(75.00,30.00)
\put(11.67,20){\makebox(0,0)[rc]{$\underbrace{\phantom{MMMll}}_{\overrightarrow{\delta\ba}}$}}
\put(12.67,24.67){\makebox(0,0)[rc]{$\langle^{(1)}\Mars\ {}^{(2)}\Venus\rangle$}}
\put(27.67,15.67
){\makebox(0,0)[cc]{$\langle^{(3)}\Venus\ {}^{(4)}\Mars\rangle$}}
\put(27.67,20.5){\makebox(0,0)[cc]{$\overbrace{\phantom{MMMlll}}^{\overrightarrow{\delta\boldb}}$}}
\put(48.00,20.17
){\makebox(0,0)[lc]{{\Large$\mapsto$}}}
\put(60.00,24.67){\makebox(0,0)[lc]{$\langle^{(1)}\Mars\,\qquad{}^{(3)}\Venus\rangle$}}
\put(70.00,15.67){\makebox(0,0)[lc]{$\langle^{(2)}\Venus\,\qquad{}^{(4)}\Mars\rangle$}}
\end{picture}
\vspace{-10mm}
\caption{The on\/-\/the\/-\/diagonal reconfiguration of couplings is the ope\-r\-a\-tio\-nal definition of BV \/Laplacian~$\Delta$; the variations are normalised by~\eqref{EqNormalizeShifts}.}\label{ReconfigCouplingsBV}
\end{center}
\end{figure}

\begin{proof}[The analytic construction of BV Laplacian~$\Delta$]
First, let us consider an \emph{integral} functional 
$F=\int f\bigl(\bx,[\ba],[\boldb]\bigr)\circ\dvol(\bx)
\in\bar{H}^n(\bpiNCZO 
)$. Let
$\delta a^{i_1}(\bby_1)\cdot\vec e_{i_1}(\bby_1)$ and $\delta b_{i_2}(\bby_2)\cdot\vec e^{{}\,\dagger,i_2}(\bby_2)$
be a pair of test shifts of the parity-even and odd letters in the fibre alphabet; assume 
normalization~\eqref{EqNormalizeShifts}. Construct the second variation%
\footnote{Summation over the (multi)\/indices~$i_\alpha$,\ $\sigma$, $\tau$ or the like is implicit in this formula and in what follows.} 
\begin{multline*}
\overrightarrow{\delta\ba}(\overrightarrow{\delta\bb}(F))=\iint_{M^n}\Id\bby_1\,\Id\bby_2\int
\Bigl\{
(\delta a^{i_1})\Bigl(\frac{\overleftarrow\dd}{\dd\bby_1}\Bigr)^{\sigma_1}(\bby_1)\cdot
\left\langle\vec e_{i_1}(\bby_1)\Bigl|\vec e^{{}\,\dagger,i_1}(\bx)\right\rangle
\frac{\overrightarrow{\dd}}{\dd a^{i_1}_{\sigma_1}}\;\circ\\
\circ(\delta b_{i_2})\Bigl(\frac{\overleftarrow\dd}{\dd\bby_2}\Bigr)^{\sigma_2}(\bby_2)\cdot
\left\langle(-\vec e^{{}\dagger,i_2})(\bby_2)\Bigl|\vec e_{i_2}(\bx)\right\rangle
\frac{\overrightarrow{\dd}}{\dd b_{i_2,\sigma_2}}f(\bx,
[\ba],[\bb])\Bigr\}
\,\dvol(\bx).
\end{multline*}
At the end of a reasoning (of which the object $\Delta F$ could be only a small piece), the integrations by parts carry the
derivatives off the virtual test shifts, which yields
\begin{multline*}
\iint_{M^n}\Id\bby_1\,\Id\bby_2\int
\Bigl\{\delta a^{i_1}(\bby_1)\cdot\left\langle\vec e_{i_1}(\bby_1)\bigl|\vec e^{{}\,\dagger.i_1}(\bx)\right\rangle\cdot\\
\cdot\delta b_{i_2}(\bby_2)\cdot\left\langle(-\vec e^{{}\,\dagger,i_2})(\bby_2)\bigl|\vec e_{i_2}(\bx)\right\rangle
\vphantom{\Bigl|}^{\lceil}\Bigl(-\frac{\Id}{\Id\bx}\Bigr)^{\sigma_1\cup\sigma_2}\vphantom{\Bigr|}^{\rceil}
\frac{\overrightarrow{\dd}^2}{\dd a^{i_1}_{\sigma_1}\dd b_{i_2,\sigma_2}}
f(\bx,
[\ba],[\bb])\Bigr\}\,\dvol(\bx).
\end{multline*}
Finally, the two pairs of couplings are reconfigured according to the scenario in Fig.~\ref{ReconfigCouplingsBV}, which gives
the action of operator
\begin{equation*}
\iint_{M^n}\Id\bby_1\,\Id\bby_2
\left\{
\begin{matrix}
\langle\delta a^{i_1}(\bby_1)\vec e_{i_1}(\bby_1)|&&|\delta b_{i_2}(\bby_2)\cdot(-\vec e^{{}\,\dagger,i_2})(\bby_2)\rangle\\
&\langle\vec e^{{}\,\dagger.i_1}(\bx)|&&|\vec e_{i_2}(\bx)\rangle
\end{matrix}
\right\}
\end{equation*}
on the basic (co)vectors over $\bx\in M^n$. The couplings wright 
the diagonal $i_1=i_2$ in the summation over the indexes. 
Normalization~\eqref{EqNormalizeShifts} and the couplings 
values~\eqref{EqTwoCouplings} make 
each line in the formula above equal to $-1$; their product equals unit.
\end{proof}

\begin{cor}
In particular, this gives us the integrand of~$\Delta F$ whenever this object \emph{is} the endpoint of a reasoning;
namely, we obtain
\begin{equation*}
\sum_{i=1}^m\sum_{\substack{|\sigma_1|\geqslant0\\|\sigma_2|\geqslant0}}
\left(-\frac{\overrightarrow{\Id}}{\Id\bx}\right)^{\sigma_1\cup\sigma_2}
\left(\frac{\overrightarrow{\dd}^2}{\dd a^i_{\sigma_1}\dd b_{i,\sigma_2}}f\right)
(\bx,
[\ba],[\bb]).
\end{equation*}
We emphasize that, should the object $\Delta F$ itself be a constituent element of a larger expression, other partial
derivatives $\overrightarrow{\dd}/\dd a^{j_1}_{\tau_1}$ or 
$\overrightarrow{\dd}/\dd b_{j_2,\tau_2}$ could accumulate
at the given density $f$ of the functional $F$, whereas all the powers of minus the total derivatives would still gather outside those higher\/-\/order partial derivatives.
\end{cor}

\begin{lemma}\label{LBVOperatorDifferential}
The linear operator
\[
\Delta\colon \bar{H}^{n(1+k)}\bigl(\BBT\bpiNCZO 
\bigr)\longrightarrow\bar{H}^{n(2+k)}\bigl(\BBT\bpiNCZO 
\bigr)
\]
is a differential for every~$k\geqslant0$.
\end{lemma}

\begin{proof}
The idea is as follows: if two normalised variations are interchanged 
in an integral fun\-c\-tio\-nal within the image of~$\Delta^2$, this 
yields an indistinguishable result of opposite sign.%
\footnote{It is readily seen that this mechanism establishes the property~$\Delta^2=0$ of the BV~Laplacian to be a differential whenever acting on any local, i.\,e.\ not only integral functional.
Indeed, within the definition of~$\Delta$, both the derivations with respect to the generators work by the graded Leibniz rule along the argument's cyclic word; it is the integrations by parts over the manifold~$M^n$ which keep track of a possible composite structure $H=F_1\times\ldots\times F_\ell$ of that cyclic word, should it be made from integral functionals~$F_1$,\ $\ldots$,\ $F_\ell\in\bar{H}^{n(1+k)}\bigl(\BBT\bpiNCZO 
\bigr)$.}

Namely, let $\delta\bolds_1=(\delta a^i_1,\delta b_{1,i})$ and 
$\delta\bolds_2=(\delta a^j_2,\delta b_{2,j})$ be two normalised 
shifts of the generators~$\ba$ and~$\boldb$, 
and let 
$H=\int h(\bx,[\ba],[\boldb])\,\dvol(\bx)$ be an integral functional over~$J^\infty\bigl(\bpiNCZO 
\bigr)$. 
(It suffices to consider the minimal 
picture $H\in\bar{H}^n\bigl(\bpiNCZO 
\bigr)$ 
without 
any variations already built into~$H$.)
By definition, we have that
\begin{multline*}
\Delta(\Delta H)(\bolds,\bolds^\dagger)=\int_M\Id\bz_1\int_M\Id\bz_2\int_M\Id\bby_1\int_M\Id\bby_2\int_M\dvol(\bx)\cdot\\
\cdot
\Biggl\{
\left\langle
(\delta a^{\alpha}_1)
\left(\tfrac{\overleftarrow{\dd}}{\dd\bz_1}\right)^{\sigma_1}(\bz_1)\,
\overbrace{\vec{e}_{\alpha}(\bz_1),(-\vec{e}^{{}\,\dagger \alpha})(\bz_2)}^{-1}\,
(\delta b_{1,\alpha})
\left(\tfrac{\overleftarrow{\dd}}{\dd\bz_2}\right)^{\sigma_2}(\bz_2)
\right\rangle
\underbrace{\left\langle
\vec{e}^{{}\,\dagger \alpha}(\bx),\vec{e}_{\alpha}(\bx)\right\rangle}_{-1}\\
{}\quad\left\langle
(\delta a^{\beta}_2)
\left(\tfrac{\overleftarrow{\dd}}{\dd\bby_1}\right)^{\tau_1}(\bby_1)\,
\overbrace{\vec{e}_{\beta}(\bby_1),(-\vec{e}^{{}\,\dagger \beta})(\bby_2)}^{-1}\,
(\delta b_{2,\beta})
\left(\tfrac{\overleftarrow{\dd}}{\dd\bby_2}\right)^{\tau_2}(\bby_2)
\right\rangle
\underbrace{\left\langle
\vec{e}^{{}\,\dagger \beta}(\bx),\vec{e}_{\beta}(\bx)\right\rangle}_{-1}\\
\frac{\overrightarrow{\dd}}{\dd a^{\alpha}_{\sigma_1}}
\frac{\overrightarrow{\dd}}{\dd b_{\alpha,\sigma_2}}
\frac{\overrightarrow{\dd}}{\dd a^{\beta}_{\tau_1}}
\frac{\overrightarrow{\dd}}{\dd b_{\beta,\tau_2}}
h(\bx,[\ba],[\boldb])
\left.\Biggr\}\right|_{\jet^{\infty}_{\bx}(\bolds,\bolds^\dagger)}.
\end{multline*}
By exchanging the integrand's upper two lines and then relabelling $\alpha\rightleftarrows\beta$, $\sigma\rightleftarrows\tau$ so that $\delta a_1^{\alpha}\rightleftarrows\delta a_2^{\beta}$ and
$\delta b_{1,\alpha}\rightleftarrows\delta b_{2,\beta}$, and by swapping the reference
$\bby\rightleftarrows\bz$ to copies of the base manifold $M^n$, we almost recover the initial expression (which should
be the case), yet the order in which the parity\/-\/odd partial derivatives follow is inverse,
\[
\frac{\overrightarrow{\dd}}{\dd b_{\alpha,\sigma_2}}
\circ\frac{\overrightarrow{\dd}}{\dd b_{\beta,\tau_2}}
\longmapsto
\frac{\overrightarrow{\dd}}{\dd b_{\beta,\tau_2}}
\circ\frac{\overrightarrow{\dd}}{\dd b_{\alpha,\sigma_2}}
=-\frac{\overrightarrow{\dd}}{\dd b_{\alpha,\sigma_2}}
\circ\frac{\overrightarrow{\dd}}{\dd b_{\beta,\tau_2}}.
\]
Therefore the integrand of functional~$\Delta^2(H)$ vanishes, which proves the assertion.
\end{proof}

\begin{rem}[The geometric realisations of $\Delta$]\label{RemV1V2}
There are at least two schemes to algorithmically define the BV Laplacian~$\Delta$ in the noncommutative set\/-\/up: on the basis of a minimal model, which we denote by {\footnotesize$)($} in Fig.~\ref{FigV1V2} below, a larger construction~$\asymp$ can be built. Still both options reproduce the same structure~$\Delta$ whenever the alphabet $a^i$,\ $b_i$ is proclaimed graded\/-\/commutative.%
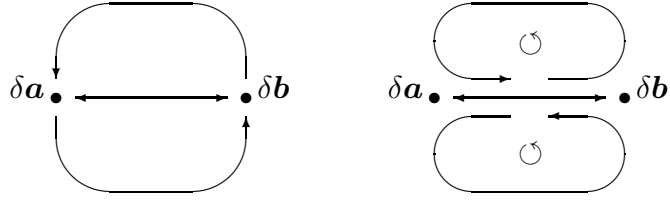
\begin{figure}[htb]
\[
{\unitlength=1mm
\begin{picture}(25,25)(0,-12.5)
\put(0,0){\circle*{1.5}}
\put(25,0){\circle*{1.5}}
\put(-1.5,0){\llap{$\delta\ba$}}
\put(26.5,0){$\delta\boldb$}
\put(2.5,0){\line(1,0){20}}
\put(2.5,0){\vector(-1,0){0}}
\put(22.5,0){\vector(1,0){0}}
\put(12.5,2.5){\oval(25,20)[t]}
\put(0,2.5){\vector(0,-1){0}}
\put(12.5,-2.5){\oval(25,20)[b]}
\put(25,-2.5){\vector(0,1){0}}
\end{picture}}
\qquad\qquad\qquad
{\unitlength=1mm
\begin{picture}(25,25)(0,-12.5)
\put(0,0){\circle*{1.5}}
\put(25,0){\circle*{1.5}}
\put(-1.5,0){\llap{$\delta\ba$}}
\put(26.5,0){$\delta\boldb$}
\put(2.5,0){\line(1,0){20}}
\put(2.5,0){\vector(-1,0){0}}
\put(22.5,0){\vector(1,0){0}}
\put(12.5,7.5){\oval(25,10)[t]}
\put(15,7.5){\oval(20,10)[rb]}
\put(10,7.5){\oval(20,10)[lb]}
\put(10,2.5){\vector(1,0){0}}
%
\put(12.5,-7.5){\oval(25,10)[b]}
\put(15,-7.5){\oval(20,10)[rt]}
\put(10,-7.5){\oval(20,10)[lt]}
\put(15,-2.5){\vector(-1,0){0}}
\put(11,6){$\circlearrowleft$}
\put(11,-8.5){$\circlearrowleft$}
\end{picture}}
\]
\caption{There remains only one cyclic word within the minimal scheme~{\footnotesize$)($} yet there appears a product of two cyclic (sub)\/words if the scheme~$\asymp$ is adopted.}\label{FigV1V2}
\end{figure} %
The minimal option~{\footnotesize$)($} suggests an orientation\/-\/preserving attachment~${\downarrow}{\uparrow}$ of the respective pairs of loose ends in the argument
$F=\int f(\bx,[\ba],[\boldb])\circ\dvol(\bx)$ of~$\Delta$. 
Specifically, for a cyclic word $(f)=w(\bx)\cdot(c_1\circ\ldots\circ c_\lambda)$ written using letters~$c_\alpha$ from the alphabet $a^i_\sigma$,\ $b_{j,\tau}$ (and weighted by a smooth coefficient~$w$ depending on~$\bx\in M^n$), the BV~Laplacian yields the sum (in a term portrayed here, without loss of generality w.r.t.\ the sequential order of the letters~$a$ and~$b$), $\Delta(F)={}$
\begin{multline*}
\sum_{i=1}^m \sum_{\substack{|\sigma|\geqslant0\\ |\tau|\geqslant0}} 
\int \bigl(-\frac{\Id}{\Id\bx}\bigr)^{\sigma\cup\tau}
\frac{\vec{\dd}^2}{\dd a^i_\sigma \dd b_{i,\tau}} w(\bx)\cdot
\bigl(c_1\dots c_{\mu-1} a^{i_0}_{\sigma_0} c_{\mu+1}\dots c_{\nu_1} b_{i_0,\tau_0} c_{\nu+1}\dots c_\lambda\bigr)\cdot\dvol(\bx)\\
\stackrel{ )( }{=} \sum_{i_0=1}^m \sum_{\substack{|\sigma_0|\geqslant0\\ |\tau_0|\geqslant0}}
\pm\int \bigl(-\frac{\Id}{\Id\bx}\bigr)^{\sigma\cup\tau}
w(\bx)\cdot\bigl(c_1\dots c_{\mu-1} \circ c_{\mu+1}\dots c_{\nu_1} \circ c_{\nu+1}\dots c_\lambda\bigr)\cdot\dvol(\bx).
\end{multline*}
The sign~$\pm$ in front of each integral is $(-)^{|\dd/\dd b|\cdot|c_1\dots c_{\nu-1}| }$.
The two pairs of adjacent loose ends, namely, $c_{\mu-1}\circ\widehat{a^{i_0}_{\sigma_0}}\circ c_{\mu+1}$ and $c_{\nu-1}\circ\widehat{b_{i_0,\tau_0}}\circ c_{\nu+1}$, link to
$c_{\mu-1}\circ c_{\mu+1}$ and $c_{\nu-1}\circ c_{\nu+1}$ respectively, so that the integrand of every term in~$\Delta(F)$ is a cyclic word in which two letters were erased but the cyclic ordering of all the remaining letters is preserved.

Conversely, according to the second scheme, which we denote by~$\asymp$, the cyclic word~$(f)$ is disrupted at both $a^{i_0}_{\sigma_0}$ and $b_{i_0,\tau_0}$; next, either of the strings of adjacent symbols $c_{\mu+1}\dots c_{\nu-1}$ and $c_{\nu+1}\dots c_\lambda c_1\dots c_{\mu-1}$ is rolled into a cyclic word. In this way, the integrand of every term in such realisation of~$\Delta(F)$ is the word which itself is the product of two cyclic words:
$\Delta(F)\stackrel{\asymp}{=}{}$
\[
\sum_{i=1}^m \sum_{\substack{|\sigma|\geqslant0\\ |\tau|\geqslant0}} \pm\int \bigl(-\frac{\Id}{\Id\bx}\bigr)^{\sigma\cup\tau} w(\bx)\cdot\langle
\text{sign}\rangle\cdot\bigl(c_{\nu+1}\dots c_\lambda c_1\dots c_{\mu-1}\bigr)\times\bigl(
c_{\mu+1}\dots c_{\nu-1}\bigr)\cdot\dvol(\bx).
\]
In this formula, the overall sign $\pm=(-)^{|\dd/\dd b|\cdot|c_1\dots c_{\nu-1}| }$ in front of the integral is the same as before. Yet now there are three ways to define the sign from the linking:
\begin{enumerate}
\item $\langle\text{sign}\rangle \mathrel{{:}{=}} (-)^{|c_{\nu+1}\dots c_\lambda|\cdot(
|c_1\dots c_{\nu-1}|+1)}$, i.e.\ $b_{i,\tau}$~is involved,
\item $\langle\text{sign}\rangle \mathrel{{:}{=}} (-)^{|c_{\nu+1}\dots c_\lambda|\cdot
|c_1\dots c_{\nu-1}|}$, or
\item $\langle\text{sign}\rangle \mathrel{{:}{=}} (-)^{|c_{\nu+1}\dots c_\lambda|\cdot
|c_1\dots c_{\mu-1}|}$ regardless of the other subword.
\end{enumerate}
It is the second variant for which a unique term $\bigl(c_1\dots c_{\mu-1}\circ c_{\mu+1}\dots c_{\nu-1}\circ c_{\nu+1}\dots c_\lambda\bigr)$ from the scheme~{\footnotesize$)($} is reproduced --\,with proper sign\,-- when the product $\bigl(c_{\nu+1}\dots c_\lambda c_1\dots  c_{\mu-1}\bigr)\times\bigl(c_{\mu+1}\dots c_{\nu-1}\bigr)$ of two cyclic words is expanded.
Let us remember however that for any choice of that sign within the scheme~$\asymp$, other cyclic words can appear. Obviously, such would be the terms in which the original consecutive order of letters along~$(f)$ is broken by the graded extension of formula~\eqref{EqDefMult} on p.~\pageref{EqDefMult}.
Indeed, under the multiplication~$\times$ the content of the first co\/-\/multiple is pasted in between \emph{every} pair of letters in the other co\/-\/multiple (and \textit{vice versa}).

We use the minimal scheme~{\footnotesize$)($} throughtour this text in the realisations of Definition~\ref{DefLaplace}. This choice is motivated in Remark~\ref{RemConventionTrueReverse}: whenever a pair of mutually inverse paths~$s^i$ and~$s^\dagger_i$ is skipped out from a given closed contour~$(f){\bigr|}_{\bs,\bs^\dagger}$, the remaining disjoint parts are linked, orientation preserved, to a new closed contour.

On the other hand, the larger scheme~$\asymp$ is reminiscent of the matrix integral methods from string theory~\cite{Bar2007,BarCRM2010}. Let $\mathsf{n}\gg 1$ and consider the algebra~$\Mat(\mathsf{n}\times\mathsf{n},\Bbbk)$ of square matrices with $\Bbbk$-\/valued entries. Recall that the $\Bbbk$-\/valued trace~$\tr$ of a product of such matrices is insensitive to cyclic permutations of comultiples. Let $\dot{\ba}=\bQ$ be a polynomial (in~$\ba$) vector field on the space of generators~$\ba$ of the matrix algebra. The divergence~$\text{div}\,\bQ$ is quadratic with respect to the traces of cyclic subwords that are formed by sub\/-\/strings of letters in the polynomial coefficients of~$\bQ$: it is readily seen that $\text{div}\,\bQ=\sum \tr(\circlearrowleft)\stackrel{\Bbbk}{\cdot}\tr(\circlearrowleft)$. Now in a larger setting, suppose that the vector field~$\bQ=\bQ^f=\{f,\cdot\}$ itself is obtained using the parity\/-\/odd symplectic form~$\Id\ba\wedge\Id\boldb$ for the double alphabet~$\ba$,\ $\boldb$, that is, $\bQ^f$~is produced by applying the skew gradient to a given Hamiltonian~$f$. Then the calculation of~$\text{div}\,(\grad f)$ goes in parallel with the construction of BV~Laplacian~$\Delta$ within the scheme~$\asymp$. Still let us remember that in such framework of~\cite{BarCRM2010}, it is the $\Bbbk$-\/valued traces which are multiplied using the operation~$\cdot$ in the ground field~$\Bbbk$, but not the cyclic words themselves (which can be multiplied using~$\times$ in the $\Bbbk$-\/algebra~$\cAZO$). In this respect the matrix integral formalism differs from our present study.
\end{rem}

\begin{define}[$\Delta(F\times G)$]\label{DefSchouten}
Let $F$ and~$G$ be integral functionals on $J^{\infty}(\bpiNCZO 
)$
and assume $F$~homogeneous. 
Applied to the product $F\times G$ of two integral functionals (see Definition~\ref{DefProduct} on p.~\pageref{DefProduct} and Definition~\ref{DefLaplace} on p.~\pageref{DefLaplace}), the 
BV~La\-p\-la\-cian~$\Delta$
is the parent structure for the (non)\/commutative variational \emph{Schouten bracket} $\lshad\,,\,\rshad$, or \emph{antibracket},
\begin{equation}\label{1b=DeviationDerivation}
\Delta(F\times G)\stackrel{\text{def}}{=}\Delta(F)\times G+(-)^{|F|}\lshad F,G\rshad+(-)^{|F|}F\times\Delta G.
\end{equation}
In other words, the bracket $\lshad\,,\,\rshad$ measures the deviation for~$\Delta$ from being a graded derivation.

The definition of $\Delta$~acting on products~$F\times G$ of (homogeneous) local functionals~$F$ and~$G$ over~$J^\infty\bigl(\bpiNCZO 
\bigr)$ is recursive; it extends by linearity to the entire space of local functionals.
\end{define}

\begin{cor}\label{PropSchoutenSkew}
The (non)\/commutative variational Schouten bracket is shifted\/-\/graded skew\/-\/symmetric:
\begin{equation}\label{EqSchoutenSkew}
\lshad F,G\rshad=-(-)^{(|F|-1)\cdot(|G|-1)}\lshad G,F\rshad
\end{equation}
for any homogeneous 
local functionals~$F$ and~$G$ over~$J^\infty\bigl(\bpiNCZO 
\bigr)$.
\end{cor}



\begin{proof}[The analytic construction of the Schouten bracket~$\lshad{\,},{\,}\rshad$]
By the graded Leibniz rule, we have that
\begin{multline*}
\bigl(\overrightarrow{\delta\ba}
\circ\overrightarrow{\delta\boldb}\bigr)(F\times G)=
\overrightarrow{\delta\ba}\Bigl(\overrightarrow{\delta\boldb}(F)\times G 
+ (-)^{|F|}F\times\overrightarrow{\delta\boldb}(G)\Bigr)
={}\\
{}
=\bigl(\overrightarrow{\delta\ba}\circ\overrightarrow{\delta\boldb}\bigr)(F)\times G
+ (-)^{|F|}\overrightarrow{\delta\ba}(F)\times\overrightarrow{\delta\boldb}(G)
+ \overrightarrow{\delta\boldb}(F)\times\overrightarrow{\delta\ba}(G)
+ (-)^{|F|}F\times\bigl(\overrightarrow{\delta\ba}\circ\overrightarrow{\delta\boldb}\bigr)(G).
\end{multline*}
Using the 
lemma
$\overrightarrow{\dd}/\dd\boldb(F)=(-)^{|F|-1}(F)\overleftarrow{\dd}/\dd\boldb$, let us reverse the direction in which the operators~$\overrightarrow{\delta\ba}$ and~$\overrightarrow{\delta\boldb}$ act on~$F$ in the second and third terms of the formula above; 
this yields\footnote{%
Further processing of the first and last terms in the formula at hand --\,that is, the on\/-\/the\/-\/diagonal reconfigurations of couplings and integrations by parts\,-- is 
analogous to the algorithm for dealing with the second and third terms, see Definition~\ref{DefLaplace}. 
The result is~\eqref{1b=DeviationDerivation}.}
\[
{}
=\bigl(\overrightarrow{\delta\ba}\circ\overrightarrow{\delta\boldb}\bigr)(F)\times G
+ (-)^{|F|}\Bigl(
(F)\overleftarrow{\delta\ba}\times\overrightarrow{\delta\boldb}(G)
- (F)\overleftarrow{\delta\boldb}\times\overrightarrow{\delta\ba}(G)
\Bigr)
+ (-)^{|F|}F\times\bigl(\overrightarrow{\delta\ba}\circ\overrightarrow{\delta\boldb}\bigr)(G).
\]
Let us have a closer look at the difference of second and third terms: for integral functionals~$F$ and~$G$ it is
\begin{align*}
{}&\iint\Id\bby_1\dvol(\bx_1)\,\bigl(f(\bx_1,[\ba],[\boldb])\bigr)
\frac{\overleftarrow{\dd}}{\dd a^{i_1}_{\sigma_1}}\,
\Bigl\langle\underleftarrow{\stackrel{\text{second}}{\vec{e}^{\,\dagger,i_1}(\bx_1)},\stackrel{\text{first}}{\vec{e}_{i_1}(\bby_1)}} \cdot
\Bigl(\frac{\overrightarrow{\dd}}{\dd\bby_1}\Bigr)^{\sigma_1}
(\delta a^{i_1})(\bby_1)\Bigr\rangle\times{}
\\
{}&{}\quad{}\times\iint\Id\bby_2
\Bigl\langle (\delta b_{i_2})\Bigl(\frac{\overleftarrow{\dd}}{\dd\bby_2}\Bigr)^{\sigma_2}(\bby_2)\cdot
\underrightarrow{\stackrel{\text{first}}{(-\vec{e}^{\,\dagger,i_2})(\bby_2)},\stackrel{\text{second}}{\vec{e}_{i_2}(\bx_2)}}
\Bigr\rangle\,\frac{\overrightarrow{\dd}}{\dd b_{i_2,\sigma_2}}\bigl(g(\bx_2,[\ba],[\boldb])\bigr)\,\dvol(\bx_2)-{}
\\
{}&-\iint\Id\bby_1\dvol(\bx_1)\,\bigl(f(\bx_1,[\ba],[\boldb])\bigr)
\frac{\overleftarrow{\dd}}{\dd b_{i_1,\sigma_1}}\,
\Bigl\langle\underleftarrow{\stackrel{\text{second}}{\vec{e}_{i_1}(\bx_1)},\stackrel{\text{first}}{(-\vec{e}^{\,\dagger,i_1})(\bby_1)}} \cdot
\Bigl(\frac{\overrightarrow{\dd}}{\dd\bby_1}\Bigr)^{\sigma_1}
(\delta b_{i_1})(\bby_1)\Bigr\rangle\times{}
\\
{}&{}\quad{}\times\iint\Id\bby_2
\Bigl\langle (\delta a^{i_2})\Bigl(\frac{\overleftarrow{\dd}}{\dd\bby_2}\Bigr)^{\sigma_2}(\bby_2)\cdot
\underrightarrow{\stackrel{\text{first}}{\vec{e}_{i_2}(\bby_2)},\stackrel{\text{second}}{\vec{e}^{\,\dagger,i_2}(\bx_2)}}
\Bigr\rangle\,\frac{\overrightarrow{\dd}}{\dd a^{i_2}_{\sigma_2}}\bigl(g(\bx_2,[\bq],[\bq^\dagger])\bigr)\,\dvol(\bx_2),
\end{align*}
where the (co)\/vectors marked `{second}' replace the respective letters in the already\/-\/built product of cyclic words~$f$ and~$g$; let us remember that in the construction of~$\Delta(F\times G)$, the multtiplication~$\times$ is performed \textit{ab initio} 
and let us bear in mind that the (co)\/vec\-tors belonging to shifts~\eqref{EqShiftsJetVariables}, here marked `first', do \emph{not} become parts of that cyclic word.
The conversion of two pairs of variations in 
$(F)\overleftarrow{\delta\ba}\times\overrightarrow{\delta\boldb}(G)
-(F)\overleftarrow{\delta\boldb}\times\overrightarrow{\delta\ba}(G)$
into one integral object --\,via integrations by parts on the diagonal $\bx_1=\bby_1=\bby_2=\bx_2$ through many consecutive reconfigurations of the couplings\,-- determines the 
functional\footnote{%
The remaining volume element can be either~$\dvol(\bx_1)$ or~$\dvol(\bx_2)$; 
its final location is prescribed by either the right\/-\/to\/-\/left or left\/-\/to\/-\/right (which is the case here) direction of couplings in the output.
From~\eqref{EqTwoCouplings} 
it is clear that a simultaneous swap ``first~$\rightleftarrows$ second'' in a pair of couplings would give the extra factor $(-1)\cdot(-1)=+1$, so that expression's overall sign does not change.}
\begin{align}
&\iiiint\Id\bx_1\Id\bby_1\Id\bby_2\,\dvol(\bx_2)\,
\Bigl\langle \delta a^{i_1}(\bby_1)\cdot
\underline{\underline{\stackrel{\text{first}}{\vec{e}_{i_1}(\bby_1)}}},
\underline{\underline{\stackrel{\text{second}}{(-\vec{e}^{\,\dagger,i_2})}(\bby_2)}} \cdot \delta b_{i_2}(\bby_2) \Bigr\rangle \cdot{}\notag\\
&\bigl(f(\bx_1,[\ba],[\boldb])\bigr)
\underbrace{\frac{\overleftarrow{\dd}}{\dd a^{i_1}_{\sigma_1}}
\vphantom{\Bigl{|}}^{\lceil}
\Bigl(-\frac{\overleftarrow{\Id}}{\Id\bby_1}\Bigr)^{\sigma_1}
\vphantom{\Bigr{|}}^{\rceil}
\circ\Bigl\langle\underline{\stackrel{\text{first}}{\vec{e}^{\,\dagger,i_1}(\bx_1)}}\Bigr|}
\times
\underbrace{\Bigl|\underline{\stackrel{\text{second}}{\vec{e}_{i_2}(\bx_2)}}
\Bigr\rangle\circ
\vphantom{\Bigl{|}}^{\lceil}
\Bigl(-\frac{\overrightarrow{\Id}}{\Id\bby_2}\Bigr)^{\sigma_2}
\vphantom{\Bigr{|}}^{\rceil}
\frac{\overrightarrow{\dd}}{\dd b_{i_2,\sigma_2}}}\bigl(g(\bx_2,[\ba],[\boldb])\bigr)-{}\notag
\\
{}&-
\iiiint\Id\bx_1\Id\bby_1\Id\bby_2\,\dvol(\bx_2)\,
\Bigl\langle 
\delta b_{i_1}(\bby_1)
\underline{\underline{\stackrel{\text{first}}{(-\vec{e}^{\,\dagger,i_1})}(\bby_1)}},
\underline{\underline{\stackrel{\text{second}}{\vec{e}_{i_2}(\bby_2)}}}
\cdot \delta a^{i_2}(\bby_2)\Bigr\rangle \cdot{}\notag
\\
&\bigl(f(\bx_1,[\ba],[\boldb])\bigr)
\underbrace{\frac{\overleftarrow{\dd}}{\dd b_{i_1,\sigma_1}}
\vphantom{\Bigl{|}}^{\lceil}
\Bigl(-\frac{\overleftarrow{\Id}}{\Id\bby_1}\Bigr)^{\sigma_1}
\vphantom{\Bigr{|}}^{\rceil}
\circ\Bigl\langle\underline{\stackrel{\text{first}}{\vec{e}_{i_1}(\bx_1)}}\Bigr|}
\times
\underbrace{\Bigl|
\underline{\stackrel{\text{second}}{\vec{e}^{\,\dagger,i_2}(\bx_2)}}
\Bigr\rangle\circ
\vphantom{\Bigl{|}}^{\lceil}
\Bigl(-\frac{\overrightarrow{\Id}}{\Id\bby_2}\Bigr)^{\sigma_2}
\vphantom{\Bigr{|}}^{\rceil}
\frac{\overrightarrow{\dd}}{\dd a^{i_2}_{\sigma_2}}}\bigl(g(\bx_2,[\ba],[\boldb])\bigr).\label{EqUniqueDef}
\end{align}
Evaluating both couplings in the minuend, we obtain $(-1)\cdot(-1)=+1$;
likewise, in the subtrahend we have that $(+1)\cdot(+1)=+1$; at every value of the indexes, the respective shift components contribute with $\delta a^\bullet\cdot\delta a^\dagger_\bullet=1$.
We emphasize that the expression~$\lshad F,G\rshad$, which has been constructed by following the couplings' 
re\/-\/attachment mechanism, itself can serve as a constituent part of a larger object. Because the reconfigurations of couplings and integrations by parts occur prior only to the restriction of output to the jet~$\jet^\infty_{\bx}(\bolds,\bolds^\dagger)$ at the diagonal~$\bx_1=\bx_2=\bby_1=\bby_2\mathrel{{=}{:}}\bx\in M^n$ for the section $(\bolds,\bolds^\dagger) 
$, 
this means that 
other partial derivatives
can freely overtake the horizontal derivatives along the base~$M^n$.
This is 
why the total derivatives were embraced by using 
$\vphantom{\bigl{|}}^\lceil\ldots\vphantom{\bigr{|}}^\rceil$ and why the shifts' own base variables~$\bby_i$ were used in~\eqref{EqUniqueDef} instead of the variables~$\bx_i$ from the functionals' densities.
\end{proof}


\begin{rem}
In effect, the only minus sign making the \emph{difference} of two terms is determined by the precedence $\ba\prec\boldb$ versus succedence $\boldb\succ\ba$, that is, by the sequential order in which the parity\/-\/even and odd partial derivatives 
are distributed between the ordered pair $F\prec G$ of input objects.
\end{rem}

\begin{cor}
Suppose 
that the Schouten bracket of integral functionals~$F$ and~$G$ is the endpoint of a calculation, that is, the reasoning stops there and the object~$\lshad F,G \rshad\colon
\Gamma(\bpiNCZO)\to\cX(\bvx^{\pm1})$ 
is used only for its evaluation at mappings~$(\bolds,\bolds^\dagger)$
but it is \emph{not} contained 
in any larger formula. 
Should this be known in advance, then one re\/-\/derives the familiar 
expression,
\begin{multline}
\lshad F,G\rshad =
\int\Bigl\{\Bigl( (f) \underbrace{\frac{\overleftarrow{\dd}}{\dd a^i_\sigma}
\vphantom{\Bigl{|}}^\lceil\Bigl(-\frac{\overleftarrow{\Id}}{\Id\bx}\Bigr)^\sigma\vphantom{\Bigr{|}}^\rceil}
\times
\underbrace{
\vphantom{\Bigl{|}}^\lceil\Bigl(-\frac{\overrightarrow{\Id}}{\Id\bx}\Bigr)^\tau\vphantom{\Bigr{|}}^\rceil
\frac{\overrightarrow{\dd}}{\dd b_{i,\tau}}} (g)\\
{}-
(f) \underbrace{\frac{\overleftarrow{\dd}}{\dd b_{i,\sigma}}
\vphantom{\Bigl{|}}^\lceil\Bigl(-\frac{\overleftarrow{\Id}}{\Id\bx}\Bigr)^\sigma\vphantom{\Bigr{|}}^\rceil}
\times
\underbrace{
\vphantom{\Bigl{|}}^\lceil\Bigl(-\frac{\overrightarrow{\Id}}{\Id\bx}\Bigr)^\tau\vphantom{\Bigr{|}}^\rceil
\frac{\overrightarrow{\dd}}{\dd a^i_\tau}} (g)
\Bigr) \bigl(\bx,[\ba],[\boldb]\bigr)\Bigr\}\,\dvol(\bx),\label{EqFamiliar}
\end{multline}
where, we remember, the multiplication~$f\times g$ is performed \textit{ab initio} to contruct the object~$F\times G$ over~$M^n\times M^n$; the underbraced operators then proceed by the four Leibniz rules along the two comultiples, either of which is built into the product but exists over the respective copy of underlying manifold~$M^n$.
\end{cor}

\begin{rem}[The geometric realisation of $\lshad\,,\,\rshad$]\label{RemGeometricSchouten}
The geometric construction of every term in the noncommutative variational Schouten bracket of integral functionals goes as follows. Without loss of generality suppose that 
either of the arguments~$F$ and~$G$ 
consists of just one 
cyclic word (otherwise, proceed by linearity).

For consistency let us first recall the geometric mechanism of left multiplication~$(F\times)\,G$.
\begin{figure}[htb]
\begin{center}
\unitlength=1mm
\linethickness{0.4pt}
\begin{picture}(120.00,18.67)(0,3)
\put(15.00,10.00){\circle{14.00}}
\put(15.00,10.00){\makebox(0,0)[cc]{{\large$\circlearrowleft$}}}
\put(15.00,17.00){\circle*{1.33}}
\put(15.00,18.67){\makebox(0,0)[cb]{$\pmb{\infty}_F$}}
\put(-2.00,12.00){\makebox(0,0)[lc]{$F={}$}}
\put(20.67,15.33){\makebox(0,0)[lb]{$r$}}
\put(21.33,10.00){\line(1,0){1.33}}
\put(17.00,15.67){\makebox(0,0)[cc]{$\cdot$}}
\put(19.33,14.33){\makebox(0,0)[cc]{$\cdot$}}
\put(20.33,12.33){\makebox(0,0)[cc]{$\cdot$}}
\put(25.00,10.00){\makebox(0,0)[lc]{{\Huge$\mapsto$}}}
\put(45.00,10.00){\circle{14.00}}
\put(45.00,10.00){\makebox(0,0)[cc]{{\large$\circlearrowleft$}}}
\put(45.00,18.33){\makebox(0,0)[cb]{$\pmb{\infty}_F$}}
\put(37.33,10.00){\line(1,0){1.33}}
\put(38.33,15.00){\makebox(0,0)[lb]{$r$}}
\put(39.33,12.67){\makebox(0,0)[cc]{$\cdot$}}
\put(40.67,14.33){\makebox(0,0)[cc]{$\cdot$}}
\put(42.67,15.67){\makebox(0,0)[cc]{$\cdot$}}
\put(45.00,17.00){\circle*{1.33}}
\put(75.00,10.00){\circle{14.00}}
\put(75.00,10.00){\makebox(0,0)[cc]{{\large$\circlearrowleft$}}}
\put(75.00,17.00){\circle*{1.33}}
\put(81.33,10.00){\line(1,0){1.33}}
\put(81.67,14.67){\makebox(0,0)[lb]{$p$}}
\put(75.00,18.67){\makebox(0,0)[cb]{$\pmb{\infty}_G$}}
\put(46.67,17){\vector(-1,0){1.33}}
\put(43.33,17){\vector(1,0){1.33}}
\put(73.33,17){\vector(1,0){1.33}}
\put(76.67,17){\vector(-1,0){1.33}}
\put(77.00,15.33){\makebox(0,0)[cc]{{\tiny$+$}}}
\put(79.00,13.67){\makebox(0,0)[cc]{{\tiny$+$}}}
\put(80.33,11.33){\makebox(0,0)[cc]{{\tiny$+$}}}
\put(95.00,10.00){\makebox(0,0)[rc]{{\Huge$\gets$}}}
\put(94.4,12){\line(0,-1){2.7}}
\put(110.00,10.00){\circle{14.00}}
\put(110.00,10.00){\makebox(0,0)[cc]{{\large$\circlearrowleft$}}}
\put(110.00,17.00){\circle*{1.33}}
\put(110.00,18.33){\makebox(0,0)[cb]{$\pmb{\infty}_G$}}
\put(120.00,10.00){\makebox(0,0)[lc]{${}=G$}}
\put(103,14.33){\makebox(0,0)[lb]{$p$}}
\put(104.33,11.33){\makebox(0,0)[cc]{{\tiny$+$}}}
\put(106.33,13.33){\makebox(0,0)[cc]{{\tiny$+$}}}
\put(108.00,14.67){\makebox(0,0)[cc]{{\tiny$+$}}}
\put(102.33,10.00){\line(1,0){1.33}}
\end{picture}
\end{center}
\end{figure}
Namely, by using~$\gotht^{r_a+r_b}$ rotate the necklace~$F$ counterclockwise until $r_a\geqslant0$ parity\/-\/even and $r_b\geqslant0$ parity-odd symbols would have passed through the lock~$\pmb{\infty}_F$\;; when 
the $(r_a+r_b+1)$th symbol approaches $\pmb{\infty}_F$, open that lock.
Likewise, using~$\gotht^{-(p_a+p_b)}$ rotate the ring~$G$ clockwise and, as soon as $p_a\geqslant0$ parity\/-\/even and $p_b\geqslant0$ parity\/-\/odd symbols would have passed through~$\pmb{\infty}_G$, unlock~$G$ 
just before its $(p_a+p_b+1)$th symbol. 
Place the loose ends of the two open words next to each other, preserving the orientation of two strings of
symbols,
\begin{figure}[htb]
\begin{center}
\unitlength=1mm
\linethickness{0.4pt}
\begin{picture}(76.33,37.00)
\put(40.5,12){\vector(1,0){33}}
\put(40.5,12){\vector(-1,0){34}}
\put(5.33,12){\circle*{1.33}}
\put(75,12){\circle*{1.33}}
\put(3,12){\llap{$\delta\ba$ or $\delta\boldb$}}
\put(77,12){$\delta\boldb$ or $\delta\ba$}
\put(23,13){{\small $+$ or $-$, respectively}}
\put(20.00,30.00){\circle*{1.33}}
\bezier{144}(20.00,30.00)(5.00,35.00)(5.00,15.00)
\bezier{124}(5.00,15.00)(5.00,0.00)(20.00,5.00)
\bezier{68}(20.00,5.00)(25.00,10.00)(35.00,10.00)
\bezier{60}(20.00,30.00)(23.67,25.00)(35.00,25.00)
\put(33,25.00){\vector(1,0){3}}
\put(33,10){\vector(1,0){3}}
\put(60.00,30.00){\circle*{1.33}}
\bezier{68}(60.00,30.00)(55.33,25.00)(45.00,25.00)
\bezier{156}(60.00,30.00)(75.00,37.00)(75.00,15.00)
\bezier{124}(75.00,15.00)(75.00,0.00)(60.00,5.00)
\bezier{56}(60.00,5.00)(54.33,10)(45,10)
\put(48,10){\vector(-1,0){3}}
\put(48,25){\vector(-1,0){3}}
\put(14.33,19.33){\makebox(0,0)[cc]{{\large$\circlearrowleft$}}}
\put(65.00,19.33){\makebox(0,0)[cc]{{\large$\circlearrowleft$}}}
\put(20.00,31.67){\makebox(0,0)[lb]{$\pmb{\infty}_F$}}
\put(60,31.67){\makebox(0,0)[rb]{$\pmb{\infty}_G$}}
\put(25.33,27){\makebox(0,0)[lb]{$r$}}
\put(52.33,26.33){\makebox(0,0)[lb]{$p$}}
\put(35.67,17.33){\makebox(0,0)[lb]{$\langle\,{}\ ,\ \rangle$}}
\put(1.00,20.33){\makebox(0,0)[lb]{$F$}}
\put(76.33,20.33){\makebox(0,0)[lb]{$G$}}
\put(40.00,27.00){\circle*{1.33}}
\put(40.00,28.67){\makebox(0,0)[cb]{$\pmb{\infty}$}}
\put(20.67,27.00){\makebox(0,0)[cc]{$\cdot$}}
\put(23.67,25.33){\makebox(0,0)[cc]{$\cdot$}}
\put(27.00,24.00){\makebox(0,0)[cc]{$\cdot$}}
\put(52.67,23.67){\makebox(0,0)[cc]{{\tiny$+$}}}
\put(55.67,24.67){\makebox(0,0)[cc]{{\tiny$+$}}}
\put(58.00,25.67){\makebox(0,0)[cc]{{\tiny$+$}}}
\end{picture}
\end{center}
\end{figure}
and join the facing ends of the two strings, forming the new cyclic word
that inherits 
the orientation.
\footnote{Note that by the above construction, the symbols from~$F$ preserve their consecutive order when forming a sub\/-\/string in the cyclic word~$F\times G$, as well as the symbols from~$G$~do.}

From the old markers~$\pmb{\infty}_F$ and~$\pmb{\infty}_G$ where the reading of cyclic words~$F$ and~$G$ started, in opposite directions issue the derivations~$\dd/\dd a^i_\sigma$ and~$\dd/\dd b_{i,\tau}$ of opposite parities. Let one of them work against the orientation~$\circlearrowleft$, i.\,e.\ %
\emph{clockwise} over~$F$ and let the other act \emph{counterclockwise}, i.\,e.\ along the orientation on~$G$. (Each obeying the Leibniz rule, either of those derivations of course also reworks the $r_a+r_b$ --\,resp., $p_a+p_b$\,-- symbols which are found in the string of~$F$ --\,resp., in~$G$ with its~$\pmb{\infty}_G$\,-- behind the lock~$\pmb{\infty}_F$ with respect to the orientation of cyclic words. The calculation of grading and parity then involves negative integer numbers.) The antecedence $\left.\dd/\dd a^i_\sigma\right|_{F} \prec
\left.\dd/\dd b_{i,\tau}\right|_{G}$ yields the plus sign, whereas the opposite sequential order of~$F$ vs~$G$ yields the minus sign in front of the corresponding term in the Leibniz rule expansions.\footnote{One easily recognises the sign convention from~\eqref{EqTwoCouplings} in the antecedence of derivations.}
In every such term we integrate by parts in order to shake~$|\sigma|$ and~$|\tau|$ derivatives off the arguments~$a^i_\sigma$ and~$b_{i,\tau}$ of two derivations. Recall that the emerging powers of minus the total derivatives now act in~$F\times G$ over~$M^n\times M^n$ \emph{only} on the sub\/-\/strings from the words~$F$ or~$G$ where the symbol~$a^i_\sigma$ and~$b_{i,\tau}$ initially belonged to, see~\eqref{EqUniqueDef}.

Finally, rotate the letters around the new word counterclockwise so that the old location of $\pmb{\infty}_G$ in between the symbols of~$G$ or after to the last symbol of~$G$ reaches the new linking $\pmb{\infty}_{\lshad F,G\rshad}$ of strings, nearest to $\pmb{\infty}_G$ in the positive direction. 
The terminal configuration is displayed here;\label{PicFinalInSchouten} 
\begin{figure}[h]
\begin{center}
\unitlength=2mm
\linethickness{0.4pt}
\begin{picture}(27.00,19.00)(6,3)
\bezier{100}(13.00,10.00)(13.5,16.5)(20.00,17.00)
\bezier{100}(27.00,10.00)(26.5,16.5)(20.00,17.00)
\bezier{100}(13.00,10.00)(13.5,3.5)(20.00,3.00)
\bezier{100}(27.00,10.00)(26.5,3.5)(20.00,3.00)
\put(20,17){\circle*{0.7}}
\put(13.00,10.00){\circle*{0.7}}
\put(14.9,5.3){\vector(1,-1){0.2}}
\put(14.9,14.7){\vector(1,1){0.2}}
\put(16,14.8){\makebox(0,0)[lc]{{\tiny$+$}}}
\put(17,15.5){\makebox(0,0)[lc]{{\tiny$+$}}}
\put(18,16){\makebox(0,0)[lc]{{\tiny$+$}}}
\put(13.4,11){\makebox(0,0)[lc]{$\cdot$}}
\put(13.8,12){\makebox(0,0)[lc]{$\cdot$}}
\put(14.2,13){\makebox(0,0)[lc]{$\cdot$}}
\bezier{30}(14.6,4.6)(15,5)(15.5,5.5)
\bezier{30}(14.6,15.4)(15,15)(15.5,14.5)
\put(20.00,10.00){\makebox(0,0)[cc]{{\large$\circlearrowleft$}}}
\put(9.00,9.00){\makebox(0,0)[lb]{$\pmb{\infty}_F$}}
\put(18.00,18.00){\makebox(0,0)[lb]{$\pmb{\infty}_{\lshad F,G\rshad}=\pmb{\infty}_G$}}
\put(17.67,13.67){\makebox(0,0)[cc]{{\scriptsize$G$}}}
\put(25,9){\makebox(0,0)[cc]{{\scriptsize$G$}}}
\put(22,14.7){\makebox(0,0)[cc]{{\scriptsize$G$}}}
\put(16.5,16.8){\makebox(0,0)[cc]{$p$}}
\put(13,13){\makebox(0,0)[cc]{$r$}}
\put(15,7.2){\makebox(0,0)[cc]{{\scriptsize$F$}}}
\put(15.5,12.2){\makebox(0,0)[cc]{{\scriptsize$F$}}}
\put(17.33,5.00){\makebox(0,0)[cc]{{\scriptsize$G$}}}
\end{picture}
\end{center}
\end{figure}
it carries $|F|+|G|-1$~parity\/-\/odd symbols, it preserves the orientation of both the input words~$F$ and~$G$, and it carries the sign factor determined by the ordered coupling of (co)\/vectors.
\end{rem}

\begin{cor}\label{CorSchoutenAsDerivation}
For a given homogeneous integral functional $F\in\bar{H}^n(\bpiNCZO 
)$ of grading~$|F|$, the operator $\lshad F,\,\cdot\,\rshad$ proceeds over letters of its cyclic\/-\/word(s) argument by the graded Leibniz rule (and by linearity); 
this operator's proper grading $|\lshad F,\,\cdot\,\rshad|$ is~$|F|-1$.
\end{cor}

\begin{cor}\label{CorSchoutenOnProduct}
The bi\/-\/linear (non)commutative variational Schouten bracket $\lshad\,,\,\rshad$ itself is 
a shifted\/-\/graded derivation of the product~$\times$ in the algebra of local functionals:
\begin{equation}\label{1a}
\lshad F,G\times H\rshad=\lshad F,G\rshad\times H+(-)^{(|F|-1)\cdot|G|}G\times\lshad F,H\rshad,
\end{equation}
where $F$ and~$G$ are assumed homogeneous and where both terms in the right\/-\/hand side are understood as applications of~$\lshad F,{\cdot}\rshad$ to the cyclic word~$G\times H$ within the BV~Laplacian action $\Delta\bigl(F\times(G\times H)\bigr)$ on the non\/-\/associative product of three comultiples.
\end{cor}

\begin{proof}
It is 
clear that the terms in $\lshad F,G\times H\rshad$ are
grouped in two parts: those in which the 
parity\/-\/odd 
derivations $\overrightarrow{\dd}/\dd b_{i,\tau}$ act on~$G$ and those for~$H$; the former do not contribute with any extra sign factors whereas the latter do --- in a way which depends on the parity~$\GH{G}$. This means that $\lshad F,G\times H\rshad=\lshad F,G\rshad\times H+\ldots$ in terms of~$\lshad F,\cdot\rshad$ acting on the product~$G\times H$. 
Proceeding by linearity if necessary, suppose also that $H$~is also homogeneous.
To grasp the sign in front of the term which has been omitted, let us swap the graded multiples~$G$ and~$H$. We have that
$G\times H=(-)^{\GH{G}\cdot\GH{H}}H\times G$, whence $\lshad F,G\times H\rshad=(-)^{\GH{G}\cdot\GH{H}}\lshad F,H\rshad\times G+\cdots$ in terms of $\lshad F,\cdot\rshad$~acting on the product~$H\times G$.
By recalling that the grading~$\GH{\lshad F,H\rshad}$ of the respective class of substrings in~$\lshad F,G\times H\rshad$ equals~$\GH{F}+\GH{H}-1$, we conclude that
\[
\lshad F,G\times H\rshad=\lshad F,G\rshad\times H+(-)^{\GH{G}\cdot\GH{H}}(-)^{(\GH{F}+\GH{H}-1)\cdot\GH{G}}G\times \lshad F,H\rshad,
\]
which yields formula~\eqref{1a}.
\end{proof}

\begin{rem}
Shifted\/-\/graded skew\/-\/symmetry~\eqref{EqSchoutenSkew} of the noncommutative variational Schouten bra\-cket for homogeneous \emph{local} functionals $F,G\in\overline{\mathfrak{M}}^n\bigl(\bpiNCZO 
\bigr)$ 
can now be 
re\/-\/derived, from Corollaries
~\ref{PropSchoutenSkew} and
~\ref{CorSchoutenOnProduct}, by induction on the respective numbers $\ell'$,\ $\ell''$ of building blocks in the arguments~$F$ and~$G$.
\end{rem}

\begin{theor}\label{Theor1}
Let $F$,\ $G$,\ and $H$ be homogeneous integral functionals on $J^{\infty}(
\bpiNCZO 
)$ so that their
gradings are $|F|,\ |G|$,\ and~$|H|$ respectively. Then each of the following three 
tautologically equivalent statements is valid\textup{:}
\begin{itemize}
\item[(i)] 
The noncommutative variational Schouten bracket satisfies the shifted\/-\/graded Jacobi identity
\begin{multline*}
(-)^{(|F|-1)\cdot(|H|-1)}\lshad F,\lshad G,H\rshad\rshad+
(-)^{(|F|-1)\cdot(|G|-1)}\lshad G,\lshad H,F\rshad\rshad+{}\\
{}+
(-)^{(|G|-1)\cdot(|H|-1)}\lshad H,\lshad F,G\rshad\rshad=0.
\end{multline*}
\item[(ii)]
The Jacobi identity for the bracket $\lshad\,,\,\rshad$ is the graded Leibniz rule for the operator
$\lshad F,\,\cdot\,\rshad$ acting on $\lshad G,H\rshad$, namely,
\begin{equation}\label{Jacobi4Schouten}
\lshad F,\lshad G,H\rshad\rshad=\lshad\lshad F,G\rshad,H\rshad+(-)^{(|F|-1)\cdot(|G|-1)}\lshad G,\lshad F,H\rshad\rshad.
\end{equation}
\item[(iii)]
The 
graded commutator of operators $\lshad F,\,\cdot\,\rshad$ and $\lshad G,\,\cdot\,\rshad$ is equal to the ope\-ra\-tor $\lshad\lshad F,G\rshad\,\cdot\,\rshad$, that is,
\begin{equation}\label{JacobiCommutator}
\lshad F,\lshad G,\,\cdot\,\rshad\rshad(H)-(-)^{(|F|-1)\cdot(|G|-1)}\lshad G,\lshad F,\,\cdot\,\rshad\rshad(H)=
\lshad\lshad F,G\rshad,\,\cdot\,\rshad(H).
\end{equation}
\end{itemize}
The arrangement of parentheses in~\eqref{JacobiCommutator} is~$(F\times G)\times H$\textup{;}
both the other variants~\textup{(}i\/--\/ii\textup{)} are obtained from~\eqref{JacobiCommutator} using multiplication by sign factors.%
\footnote{\label{FootWhyFGbyHinJacobi}%
Each reading of the Jacobi identity for~$\lshad\,,\,\rshad$ is valid regardless of the sequential order of multiplications in $F\times G\times H$ after a reduction to the graded\/-\/commutative set\/-\/up. From the first paragraph in the proof below it is seen why the parentheses configuration is $(F\times G)\times H$ in the non\/-\/associative setting. In the meantime, we conclude that the Jacobi identity for~$\lshad\,,\,\rshad$ renders the fact that the commutator of adjoint actions is the adjoint action by the bracket, cf.~\cite{Galli10}.}%
\end{theor}
\noindent%
Proven immediately below for the case of integral building blocks from
$\bar{H}^n(\bpiNCZO 
)$, assertion~(\textit{iii}) of Theorem~\ref{Theor1} is then extended by induction to the space $\overline{\mathfrak{M}}^n\bigl(\bpiNCZO 
\bigr)$ of noncommutative local functionals.

\begin{proof}
Consider the consecutive action of operators $\lshad F\,\cdot\,\rshad$ and $\lshad G\,\cdot\,\rshad$ of gradings $|F|-1$ and $|G|-1$, respectively, on an integral functional~$H$. Each operator proceeds over letters in every cyclic word of~$H$ by the graded Leibniz rule. It is readily seen that by taking the 
graded \emph{difference} of the two applications, as it stands in the left\/-\/hand side of~\eqref{JacobiCommutator}, we cancel all the terms in which the strings
of symbols from~$F$ and~$G$ are pasted into~$H$ not hitting each other (that is, rather staying next to each other or
becoming separated by the argument's own letters). Therefore, both sides of~\eqref{JacobiCommutator} contain the
\emph{second} variation of~$F$ or~$G$ but only the \emph{first} variation of~$H$.

Note further that all the integrals by parts always involve 
only the letters that belong to (what remains of) the functional
which is varied, see section~\ref{SecElementsGVBV}. Consequently, both sides of~\eqref{JacobiCommutator} contain the
same configurations of powers of total derivatives that fall on the letters from the second or first, first or second, and first variations of $F$, $G$, and~$H$, respectively. This shows that it is sufficient to inspect the matching of signs~--- as they occur in the
left- and right\/-\/hand side of~\eqref{JacobiCommutator}~--- in front of the insertions of symbols from $F$ into $G$, and
\textit{vice versa}. Without loss of generality, let us suppose that each of the functionals $F$ and $G$ consist of just
a single cyclic word.

Every term in $\lshad G,\,\cdot\,\rshad(H)$ is obtained from the cyclic words

\medskip

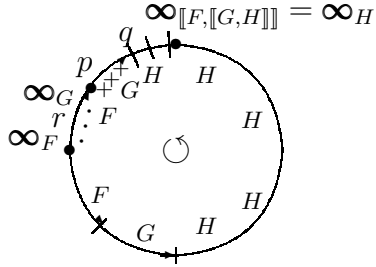
\begin{center}
\unitlength=1mm
\linethickness{0.4pt}
\begin{picture}(67.00,23.33)
\put(15.00,15.00){\circle{14.00}}
\put(15.00,22.00){\circle*{1.33}}
\put(13.67,23.33){\makebox(0,0)[lb]{$\pmb{\infty}_G$}}
\put(15,15){\makebox(0,0)[cc]{{\large$\circlearrowleft$}}}
\put(-2.00,13.33){\makebox(0,0)[lb]{$G={}$}}
\put(30.00,13.33){\makebox(0,0)[lb]{and}}
\put(40.00,13.33){\makebox(0,0)[lb]{$H={}$}}
\put(60.00,15){\circle{14.00}}
\put(60.00,22){\circle*{1.33}}
\put(57.67,23.33){\makebox(0,0)[lb]{$\pmb{\infty}_H$}}
\put(60,15){\makebox(0,0)[cc]{{\large$\circlearrowleft$}}}
\put(52.00,19.33){\makebox(0,0)[lb]{$q$}}
\put(22.00,19.33){\makebox(0,0)[lb]{$p$}}
\put(21.00,15.00){\line(1,0){2.00}}
\put(52.00,15){\line(1,0){2.00}}
\bezier{20}(52.33,17.00)(53.00,16.67)(54.00,16.33)
\bezier{20}(53,19.00)(53.67,18.6)(54.5,18.25)
\bezier{20}(54.33,20.43)(55,19.77)(55.67,19.43)
\bezier{20}(56.33,22)(56.67,21.33)(57,20.67)
\put(16.33,19.67){\makebox(0,0)[lb]{{\tiny$+$}}}
\put(18.33,17.67){\makebox(0,0)[lb]{{\tiny$+$}}}
\put(19.67,16.00){\makebox(0,0)[lb]{{\tiny$+$}}}
\end{picture}
\end{center}

\vspace{-25pt}

\noindent
as follows (see Remark~\ref{RemGeometricSchouten}). First, the ring $G$ is rotated counterclockwise, transporting $p$ odd
symbols through $\pmb{\infty}_G$, which gives the sign $(-)^{p\cdot(|G|-1)}$, and then $G$ is unlocked at $\pmb{\infty}_G$.
At the same time, $H$~is rotated clockwise and unlocked as soon as $q$~odd letters would have passed the lock
$\pmb{\infty}_H$. 
The word obtained from~$G$ 
is pasted, orientation preserved, into the similarly
obtained fragments of~$H$; the loose ends of the two strings are joined, making a new circle.  
Contracting one pair of variations $(\delta\ba,\delta\bb)$ destroys one parity-odd symbol in either~$G$
or~$H$. 
Finally, the $q$~parity\/-\/odd letters of~$H$ are pushed counterclockwise~--- so many of them that the old~$\pmb{\infty}_H$ coincides with $\pmb{\infty}_{\lshad G,H\rshad}$, placed at the moment of linking at the concatenation of strings' loose ends nearest to~$\pmb{\infty}_H$ in positive direction. 
The sign factor which is gained when the lock of $H$ is restored on its proper place
equals $(-)^{q\cdot(|G|-1)}$; the \emph{minus one} in the exponent counts the parity-odd letter destroyed by the coupling. The resulting
necklace --~a term in $\lshad G,H\rshad$~-- looks like this:

\begin{center}
\unitlength=2mm
\linethickness{0.4pt}
\begin{picture}(27.00,19.00)(6,0)
\bezier{100}(13.00,10.00)(13.5,16.5)(20.00,17.00)
\bezier{100}(13.01,10.00)(13.21,12.5)(14.01,13.67)
\bezier{100}(27.00,10.00)(26.5,16.5)(20.00,17.00)
\bezier{100}(13.00,10.00)(13.5,3.5)(20.00,3.00)
\bezier{100}(27.00,10.00)(26.5,3.5)(20.00,3.00)
\put(13.00,10.00){\circle*{1}}
\put(20.00,17.00){\circle*{1}}
\put(15.00,5.33){\vector(1,-1){1.33}}
\put(18.00,3.33){\vector(-3,1){1.67}}
\put(20.00,10.00){\makebox(0,0)[cc]{{\large$\circlearrowleft$}}}
\put(8.00,10.00){\makebox(0,0)[lb]{$\pmb{\infty}_G$}}
\put(20.00,18.00){\makebox(0,0)[cb]{$\pmb{\infty}_{\lshad G,H\rshad}$}}
\put(13.67,12.67){\vector(1,2){0.67}}
\bezier{28}(13.7,14.5)(14.4,14.2)(15.1,13.9)
\bezier{28}(14.20,15.5)(14.87,15)(15.53,14.6)
\bezier{28}(15.33,16.33)(15.87,15.67)(16.5,15)
\bezier{28}(16.8,17)(17.33,16)(17.5,15.54)
\bezier{28}(18.00,17.5)(18.25,16.8)(18.5,16.1)
\put(17.67,13.67){\makebox(0,0)[cc]{{\scriptsize$H$}}}
\put(25,9){\makebox(0,0)[cc]{{\scriptsize$H$}}}
\put(15.33,18.00){\makebox(0,0)[cc]{$q$}}
\put(12.00,12.33){\makebox(0,0)[cc]{$p$}}
\put(14,12.50){\makebox(0,0)[lc]{{\tiny$+$}}}
\put(13.5,10.83){\makebox(0,0)[lc]{{\tiny$+$}}}
\put(15.33,9.00){\makebox(0,0)[cc]{{\scriptsize$G$}}}
\put(15.33,6.00){\makebox(0,0)[cc]{{\scriptsize$G$}}}
\put(21.00,4.67){\makebox(0,0)[cc]{{\scriptsize$H$}}}
\end{picture}
\end{center}

\noindent
The total sign accumulated up to this moment is 
$(-)^{p\cdot(|G|-1)}\cdot(-)^{q\cdot(|G|-1)}$. Now the operator
$\lshad F,\,\cdot\,\rshad$ approaches that ring from the left. 
Arguing as above, we rotate the cyclic word

\begin{center}
\unitlength=1mm
\linethickness{0.4pt}
\begin{picture}(27.67,18.00)(1,0)
\put(20.00,10.00){\circle{14.00}}
\put(20.00,17.00){\circle*{1.00}}
\put(20.00,10.00){\makebox(0,0)[cc]{{\large$\circlearrowleft$}}}
\put(20.00,18.00){\makebox(0,0)[cb]{$\pmb{\infty}_F$}}
\put(22,15.67){\makebox(0,0)[lc]{$\cdot$}}
\put(23.85,14.33){\makebox(0,0)[lc]{$\cdot$}}
\put(24.8,12.75){\makebox(0,0)[lc]{$\cdot$}}
\put(26,15){\makebox(0,0)[lc]{$r$}}
\put(26.33,10.00){\line(1,0){1.33}}
\put(3.00,10.00){\makebox(0,0)[lc]{$F={}$}}
\end{picture}
\end{center}

\noindent
counterclockwise, letting $r$ parity-odd symbols pass through $\pmb{\infty}_F$ (this yields $(-)^{r\cdot(|F|-1)}$). Having
unlocked that ring at $\pmb{\infty}_F$, we carry this term in $\lshad F,\,\cdot\,\rshad$ of grading $|F|-1$ along the
$p+q$ parity-odd symbols in the pre-fabricated linking of $G$ and $H$. By the time the loose ends of
$\lshad F,\,\cdot\,\rshad$ reach the former location of $\pmb{\infty}_G$ in $G$, the sign factor
$(-)^{(p+q)\cdot(|F|-1)}$ is accumulated, and the configuration is this:

\begin{center}
\unitlength=2mm
\linethickness{0.4pt}
\begin{picture}(27.00,19.00)(6,1.5)
\bezier{100}(13.00,10.00)(13.5,16.5)(20.00,17.00)
\bezier{100}(13.01,10.00)(13.21,12.5)(14.01,13.67)
\bezier{100}(27.00,10.00)(26.5,16.5)(20.00,17.00)
\bezier{100}(13.00,10.00)(13.5,3.5)(20.00,3.00)
\bezier{100}(27.00,10.00)(26.5,3.5)(20.00,3.00)
\put(13.00,10.00){\circle*{.7}}
\put(20.00,17.00){\circle*{.7}}
\put(14.4,14.1){\circle*{.7}}
\put(14.3,6){\vector(1,-1){1}}
\put(19,3){\vector(1,0){1}}
\put(15.50,15){\vector(1,1){1.33}}
\put(20.00,10.00){\makebox(0,0)[cc]{{\large$\circlearrowleft$}}}
\put(9.00,10.00){\makebox(0,0)[lb]{$\pmb{\infty}_F$}}
\put(10.00,13.00){\makebox(0,0)[lb]{$\pmb{\infty}_G$}}
\put(18.00,18.00){\makebox(0,0)[lb]{$\pmb{\infty}_{\lshad F,\lshad G,H\rshad\rshad}=\pmb{\infty}_H$}}
\put(13.67,12.67){\vector(1,2){0.67}}
\put(20,2.5){\line(0,1){1}}
\bezier{28}(14.5,4.5)(15,5)(15.4,5.4)
\bezier{28}(16.8,17)(17.33,16)(17.5,15.54)
\bezier{28}(18.00,17.5)(18.25,16.8)(18.5,16.1)
\bezier{28}(19.50,17.8)(19.55,17)(19.6,16.2)
\put(18.5,15){\makebox(0,0)[cc]{{\scriptsize$H$}}}
\put(22,15){\makebox(0,0)[cc]{{\scriptsize$H$}}}
\put(22,5){\makebox(0,0)[cc]{{\scriptsize$H$}}}
\put(17,14){\makebox(0,0)[cc]{{\scriptsize$G$}}}
\put(15.5,12.3){\makebox(0,0)[cc]{{\scriptsize$F$}}}
\put(25,12){\makebox(0,0)[cc]{{\scriptsize$H$}}}
\put(14,15.5){\makebox(0,0)[cc]{$p$}}
\put(16.7,17.50){\makebox(0,0)[cc]{$q$}}
\put(12.3,12){\makebox(0,0)[cc]{$r$}}
\put(14,12.50){\makebox(0,0)[lc]{{$\cdot$}}}
\put(13.75,11.50){\makebox(0,0)[lc]{{$\cdot$}}}
\put(13.5,10.50){\makebox(0,0)[lc]{{$\cdot$}}}
\put(14.7,14){\makebox(0,0)[lc]{{\tiny$+$}}}
\put(15.2,14.7){\makebox(0,0)[lc]{{\tiny$+$}}}
\put(15.8,15.4){\makebox(0,0)[lc]{{\tiny$+$}}}
\put(15,7.00){\makebox(0,0)[cc]{{\scriptsize$F$}}}
\put(18,4.50){\makebox(0,0)[cc]{{\scriptsize$G$}}}
\put(25.00,6.67){\makebox(0,0)[cc]{{\scriptsize$H$}}}
\end{picture}
\end{center}

\noindent
By having realised the scenario which the first term in the left\/-\/hand side of~\eqref{JacobiCommutator} provides, we obtain the overall sign
\begin{equation
}\label{EqSignLHS}
(-)^{r\cdot(|F|-1)}\cdot(-)^{p\cdot(|G|-1)}\cdot(-)^{q\cdot(|G|-1)}\cdot(-)^{(p+q)\cdot(|F|-1)}=
(-)^{r\cdot(|F|-1)}\cdot(-)^{(p+q)\cdot(|F|+|G|-2)}.
\end{equation
}
Moreover, now it is clear what the extra sign contribution to the formula above would
be, should the insertion of the unlocked~$F$ start later --~with respect to the cyclic order~-- than the starting point $\pmb{\infty}_G$ of the 
turned\/-\/and\/-\/unlocked cyclic word~$G$.

On the other hand, let us calculate the overall sign factor of the very same geometric configuration in the right-hand side
of~\eqref{JacobiCommutator}. So, we first produce the respective term in $\lshad F,G\rshad$. Let us recall from the above
that the word

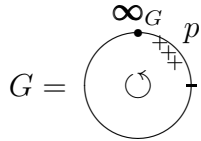
\begin{center}
\unitlength=1mm
\linethickness{0.4pt}
\begin{picture}(27.67,18.00)(2,2)
\put(20.00,10.00){\circle{14.00}}
\put(20.00,17.00){\circle*{1.00}}
\put(20.00,10.00){\makebox(0,0)[cc]{{\large$\circlearrowleft$}}}
\put(20.00,18.00){\makebox(0,0)[cb]{$\pmb{\infty}_G$}}
\put(21.75,15.67){\makebox(0,0)[lc]{{\tiny$+$}}}
\put(22.85,14.33){\makebox(0,0)[lc]{{\tiny$+$}}}
\put(23.7,13){\makebox(0,0)[lc]{{\tiny$+$}}}
\put(26,17){\makebox(0,0)[lc]{$p$}}
\put(26.33,10.00){\line(1,0){1.33}}
\put(3.00,10.00){\makebox(0,0)[lc]{$G={}$}}
\end{picture}
\end{center}

\noindent
is unlocked straight after $\pmb{\infty}_G$, but

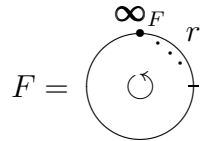
\begin{center}
\unitlength=1mm
\linethickness{0.4pt}
\begin{picture}(27.67,18.00)(0,2)
\put(20.00,10.00){\circle{14.00}}
\put(20.00,17.00){\circle*{1.00}}
\put(20.00,10.00){\makebox(0,0)[cc]{{\large$\circlearrowleft$}}}
\put(20.00,18.00){\makebox(0,0)[cb]{$\pmb{\infty}_F$}}
\put(21.75,15.67){\makebox(0,0)[lc]{$\cdot$}}
\put(23.35,14.33){\makebox(0,0)[lc]{$\cdot$}}
\put(24.7,13){\makebox(0,0)[lc]{$\cdot$}}
\put(26,17){\makebox(0,0)[lc]{$r$}}
\put(26.33,10.00){\line(1,0){1.33}}
\put(3.00,10.00){\makebox(0,0)[lc]{$F={}$}}
\end{picture}
\end{center}

\noindent
is first rotated counterclockwise by $r$ parity-odd slots; this yields the sign $(-)^{r\cdot(|F|-1)}$ and gives the word

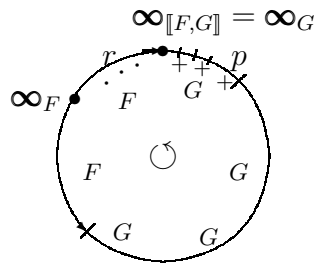
\begin{center}
\unitlength=2mm
\linethickness{0.4pt}
\begin{picture}(27.00,19.00)(6,1.5)
\bezier{100}(13.00,10.00)(13.5,16.5)(20.00,17.00)
\bezier{100}(27.00,10.00)(26.5,16.5)(20.00,17.00)
\bezier{100}(13.00,10.00)(13.5,3.5)(20.00,3.00)
\bezier{100}(27.00,10.00)(26.5,3.5)(20.00,3.00)
\put(14.2,13.8){\circle*{0.7}}
\put(20.00,17.00){\circle*{0.7}}
\put(14.8,5.2){\vector(1,-1){0.2}}
\put(18.67,17){\vector(1,0){1.33}}
\put(16,14.8){\makebox(0,0)[lc]{$\cdot$}}
\put(17,15.5){\makebox(0,0)[lc]{$\cdot$}}
\put(18,16){\makebox(0,0)[lc]{$\cdot$}}
\bezier{30}(14.6,4.6)(15,5)(15.5,5.5)
\bezier{30}(25.4,15.4)(25,15)(24.5,14.5)
\bezier{30}(21.2,17.2)(21.05,16.9)(21,16.6)
\bezier{30}(22.2,17)(22.1,16.7)(22,16.4)
\bezier{30}(23.2,16.6)(23.05,16.35)(22.9,16.1)
\put(20.00,10.00){\makebox(0,0)[cc]{{\large$\circlearrowleft$}}}
\put(10.00,13.00){\makebox(0,0)[lb]{$\pmb{\infty}_F$}}
\put(18.00,18.00){\makebox(0,0)[lb]{$\pmb{\infty}_{\lshad F,G\rshad}=\pmb{\infty}_G$}}
\put(17.67,13.67){\makebox(0,0)[cc]{{\scriptsize$F$}}}
\put(25,9){\makebox(0,0)[cc]{{\scriptsize$G$}}}
\put(22,14.4){\makebox(0,0)[cc]{{\scriptsize$G$}}}
\put(16.5,16.5){\makebox(0,0)[cc]{$r$}}
\put(25,16.3){\makebox(0,0)[cc]{$p$}}
\put(20.5,16.1){\makebox(0,0)[lc]{{\tiny$+$}}}
\put(22,15.8){\makebox(0,0)[lc]{{\tiny$+$}}}
\put(23.5,14.9){\makebox(0,0)[lc]{{\tiny$+$}}}
\put(15.33,9.00){\makebox(0,0)[cc]{{\scriptsize$F$}}}
\put(17.33,5.00){\makebox(0,0)[cc]{{\scriptsize$G$}}}
\put(23.00,4.67){\makebox(0,0)[cc]{{\scriptsize$G$}}}
\end{picture}
\end{center}

\noindent
It contains $|F|+|G|-1$ parity\/-\/odd letters; let us use it in the action of 
$\lshad\lshad F,G\rshad,\,\cdot\,\rshad$ on~$H$. By rotating the word to-paste counterclockwise by $p$ parity\/-\/odd symbols,
we gain the sign $(-)^{p\cdot(|F|+|G|-2)}$\;; proceeding by the Leibniz rule over $q$ parity-odd letters in $H$, we obtain
another sign factor $(-)^{q\cdot(|F|+|G|-2)}$. In total, the overall sign that occurs in the right\/-\/hand side
of~\eqref{JacobiCommutator} for the configuration that we knew before is
$$(-)^{r\cdot(|F|-1)}\cdot(-)^{p\cdot(|F|+|G|-2)}\cdot(-)^{q\cdot(|F|+|G|-2)}.$$
This is exactly~\eqref{EqSignLHS}.

To process --\,in both sides of~\eqref{JacobiCommutator}\,-- the configurations in which the symbols from~$G$ are pasted in between the letters of~$F$, and those are already
installed in~$H$, let us first swap~$F$ and~$G$.
By Corollary~\ref{PropSchoutenSkew}, the right\/-\/hand side of~\eqref{JacobiCommutator} becomes
$$-(-)^{(|F|-1)\cdot(|G|-1)}\lshad\lshad G,F\rshad,\,\cdot\,\rshad(H).$$
Second, multiply both sides of~\eqref{JacobiCommutator} by the sign factor $-(-)^{(|F|-1)\cdot(|G|-1)}$; this gives
$$-(-)^{(|F|-1)\cdot(|G|-1)}\lshad F,\lshad G,\,\cdot\,\rshad\rshad(H)+\lshad G,\lshad F,\,\cdot\,\rshad\rshad(H)
\quad\text{ versus }\quad\lshad\lshad G,F\rshad,\,\cdot\,\rshad(H).$$
Finally, relabel $F\rightleftarrows G$ back; by having thus recovered both sides of~\eqref{JacobiCommutator} in its
authentic form, we convert the configurations to\/-\/consider into those which we did cope with. The proof is complete.
\end{proof}

\begin{lemma}\label{LemmaBaseLapSchouten}
Let $F\in\bar{H}^{n(1+k)}\bigl(\BBT\bpiNCZO 
\bigr)$ and 
$G\in\bar{H}^{n(1+\ell)}\bigl(\BBT\bpiNCZO 
\bigr)$ be two 
integral functionals\textup{(}here $k,\ell\geqslant0$\textup{),} and assume $F$~homogeneous. Then 
\begin{equation}\label{EqLapSchouten}
\Delta\bigl(\schouten{F,G}\bigr) = \schouten{\Delta F,G} + (-)^{\GH{F}-1}\schouten{F,\Delta G}.
\end{equation}
\end{lemma}
\noindent%
This claim will be extended to all elements of 
the algebra of local functionals over $J^\infty\bigl(\bpiNCZO 
\bigr)$;
the inductive proof of Theorem~\ref{ThLaplaceOnSchouten} on p.~\pageref{ThLaplaceOnSchouten} is based on this lemma.

\begin{proof}
The key idea is that the structures $\Delta$ and $\lshad\,,\,\rshad$ yield equivalence classes of integral functionals which, after integration by parts at the end of the day, are \emph{independent} of a choice of the built-in test shifts normalized by~\eqref{EqNormalizeShifts}. Consequently, the composite structure $\Delta(\lshad{\cdot},{\cdot}\rshad)$ does not change under swapping $\delta a_1^{\alpha}\rightleftarrows\delta a_2^{\beta}$, $\delta b_{1,\alpha}\rightleftarrows\delta b_{2,\beta}$ of the respective variations~$\delta\bolds_1$ and~$\delta\bolds_2$ in~$\Delta$ and in~$\lshad\,,\,\rshad$. Hence the terms which are skew\/-\/symmetric under such exchange necessarily vanish (cf.\ the proof of Lemma~\ref{LBVOperatorDifferential} on p.~\pageref{LBVOperatorDifferential} above).

For the sake of clarity, let us assume that $F=\int f\bigl(\bx_1,[\ba],[\boldb]\bigr)\,\dvol(\bx_1)$ and
$G=\int g\bigl(\bx_2,[\ba],[\boldb]\bigr)\,\dvol(\bx_2)$ are 
building blocks from the cohomology group
$\bar{H}^n(\bpiNCZO 
)$; this simplification is legitimate because new variations which come from~$\Delta$ and
$\lshad\,,\,\rshad$ do not interfere with any other test shifts if those are already absorbed by the 
densities~$f$ and~$g$.
Suppose that $\delta\bolds_1$ and $\delta\bolds_2$ are two normalized variations of the generators~$a^i$ and~$b_i$.
By definition, we have that\footnote{To keep track of their origin, we preserve the notation for base variables~$\bby_\mu$ and~$\bz_\nu$ in the minus total derivatives acting at the end of the day on densities of the functionals~$F$ and~$G$.%
}
\begin{align*}
&
\overrightarrow{\delta\ba}\,\overrightarrow{\delta\boldb}\,%
\left(\lshad F,G\rshad\right)
=\int_M\Id\bz_1\int_M\Id\bz_2\int_M\Id\bby_1\int_M\Id\bby_2\int_M\Id\bx_1 
\int_M\dvol(\bx_2) 
\\
&\Biggl\{
\bigl\langle
(\delta a_2^{j_1})\Bigl(\tfrac{\overleftarrow{\dd}}{\dd\bz_1}\Bigr)^{\tau_1}(\bz_1)\cdot
\underrightarrow{\vec{e}_{j_1}(\bz_1),\vec{e}^{{}\,\dagger,j_1}(\cdot)}\bigr\rangle\,
\frac{\overrightarrow{\dd}}{\dd a^{j_1}_{\tau_1}} \circ
\bigl\langle
(\delta b_{2,j_2})\Bigl(\tfrac{\overleftarrow{\dd}}{\dd\bz_2}\Bigr)^{\tau_2}(\bz_2)\cdot
\underrightarrow{(-\vec{e}^{{}\,\dagger,j_2}(\bz_2),\vec{e}_{j_2}(\cdot)}\bigr\rangle\,
\frac{\overrightarrow{\dd}}{\dd b_{j_2,\tau_2}} 
\Biggr\}
\end{align*}
\begin{align*}
&\quad\smash{\Biggl[}
\bigl\langle \delta a_1^{i_1}(\bby_1)\,\underrightarrow{\vec{e}_{i_1}(\bby_1),
(-\vec{e}^{{}\,\dagger,i_2})(\bby_2)}\,\delta b_{1,i_2}(\bby_2)\bigr\rangle\cdot{}
\\
&\bigl(f(\bx_1,[\ba],[\boldb])\bigr){\frac{\overleftarrow{\dd}}{\dd a^{i_1}_{\sigma_1}}}
\vphantom{\Bigl|}^\lceil\Bigl(-\tfrac{\overleftarrow{\Id}}{\Id\bby_1}\Bigr)^{\sigma_1}\vphantom{\Bigr|}^\rceil
\,
\underline{\langle\vec{e}^{{}\,\dagger,i_1}(\bx_1)|}\times
\underline{|\vec{e}_{i_2}(\bx_2)\rangle}\,
\vphantom{\Bigl|}^\lceil\Bigl(-\tfrac{\overrightarrow{\Id}}{\Id\bby_2}\Bigr)^{\sigma_2}\vphantom{\Bigr|}^\rceil
\frac{\overrightarrow{\dd}}{\dd b_{i_2,\sigma_2}}\bigl(g(\bx_2,[\ba],[\boldb])\bigr)
-{}
\\
&{}-
\bigl\langle \delta b_{1,i_1}(\bby_1)\,\underrightarrow{(-\vec{e}^{\dagger,i_1})(\bby_1),
\vec{e}_{i_2}(\bby_2)}\,\delta a_1^{i_2}(\bby_2)\bigr\rangle\cdot{}
\\
&\bigl(f(\bx_1,[\ba],[\boldb])\bigr){\frac{\overleftarrow{\dd}}{\dd b_{i_1,\sigma_1}}}
\vphantom{\Bigl|}^\lceil\Bigl(-\tfrac{\overleftarrow{\Id}}{\Id\bby_1}\Bigr)^{\sigma_1}\vphantom{\Bigr|}^\rceil
\,
\underline{\langle\vec{e}_{i_1}(\bx_1)|}\times
\underline{|\vec{e}^{\,\dagger,i_2}(\bx_2)\rangle}\,
\vphantom{\Bigl|}^\lceil\Bigl(-\tfrac{\overrightarrow{\Id}}{\Id\bby_2}\Bigr)^{\sigma_2}\vphantom{\Bigr|}^\rceil
\frac{\overrightarrow{\dd}}{\dd a^{i_2}_{\sigma_2}}\bigl(g(\bx_2,[\ba],[\boldb])\bigr)
\Biggr]
\end{align*}
The partial derivatives 
$\overrightarrow{\dd}/\dd a^{j_1}_{\tau_1}\circ\overrightarrow{\dd}/\dd b_{j_2,\tau_2}$
are distributed between the arguments~$f$ and~$g$ by the graded Leibniz rule. 
Whenever \emph{none} of the two operators overtakes the density of~$F$,
the reconfiguration yields $\lshad\Delta F,G\rshad$. 
Likewise, if \emph{both} derivatives indexed by~$j$
overtake~$F$ and then also overtake an old derivative that fell on~$g$, 
we obtain $(-)^{\GH{F}-1}\lshad F,\Delta G\rshad$,
which is the second term in the right\/-\/hand side of~\eqref{EqLapSchouten}.

We claim that the remaining four terms cancel out by virtue of independence --\,of both~$\Delta$ and $\lshad\,,\,\rshad$\,-- of a choice of normalized virtual shifts. 

The two mixed 
terms can informally be visualised using
\[
\frac{\overrightarrow{\delta}}{\delta\boldb}(f)\frac{\overleftarrow{\delta}}{\delta\ba}\times
\frac{\overrightarrow{\delta}}{\delta\ba}\frac{\overrightarrow{\delta}}{\delta\boldb}(g)
\pm
\frac{\overrightarrow{\delta}}{\delta\ba}(f)\frac{\overleftarrow{\delta}}{\delta\boldb}\times
\frac{\overrightarrow{\delta}}{\delta\boldb}\frac{\overrightarrow{\delta}}{\delta\ba}(g).
\]
They contribute to the integrand with 
the difference of equal terms,
\begin{multline*}
\bigl\langle\delta a_1^{i_1}(\bby_1)\,\underrightarrow{\vec{e}_{i_1}(\bby_1),
 (-\vec{e}^{\,\dagger,i_2}(\bby_2)}\,\delta b_{1,i_2}(\bby_2)\bigr\rangle \cdot
\bigl\langle\delta b_{2,j_2}(\bz_2)\,\underrightarrow{(-\vec{e}^{\,\dagger,j_2})(\bz_2),
 \vec{e}_{j_1}(\bz_1)}\,\delta a_2^{j_1}(\bz_1)\bigr\rangle \cdot{}
 \\
\bigl(f(\bx_1,[\ba],[\boldb])\bigr)
\left\{
\begin{aligned}
\frac{\overleftarrow{\dd}}{\dd a^{i_1}_{\sigma_1}} \vphantom{\Bigl|}^{\lceil}
\Bigl(-\frac{\overleftarrow{\Id}}{\Id\bby_1}\Bigr)^{\sigma_1}\vphantom{\Bigr|}^{\rceil}
\langle\underline{\vec{e}^{\,\dagger,i_1}(\bx_1)}| &\times
 |\underline{\vec{e}_{i_2}(\bx_2)}\rangle\,
\vphantom{\Bigl|}^{\lceil}\Bigl(-\frac{\overrightarrow{\Id}}{\Id\bby_2}\Bigr)^{\sigma_2}\vphantom{\Bigr|}^{\rceil} \frac{\overrightarrow{\dd}}{\dd b_{i_2,\sigma_2}} 
\\
(-)^{\GH{F}-1}
\frac{\overleftarrow{\dd}}{\dd b_{j_2,\tau_2}} \vphantom{\Bigl|}^{\lceil}
\Bigl(-\frac{\overleftarrow{\Id}}{\Id\bz_2}\Bigr)^{\tau_2}\vphantom{\Bigr|}^{\rceil}
\langle\underline{\underline{\vec{e}_{j_2}(\bx_1)}}| &\times
 |\underline{\underline{\vec{e}^{\,\dagger,j_1}(\bx_2)}}\rangle\,
\vphantom{\Bigl|}^{\lceil}\Bigl(-\frac{\overrightarrow{\Id}}{\Id\bz_1}\Bigr)^{\tau_1}\vphantom{\Bigr|}^{\rceil} \frac{\overrightarrow{\dd}}{\dd a^{j_1}_{\tau_1}} 
\end{aligned}
\right\}
\bigl(g(\bx_2,[\ba],\boldb])\bigr)
\\
-
\bigl\langle\delta b_{1,i_1}(\bby_1)\,\underrightarrow{(-\vec{e}^{\,\dagger,i_1})(\bby_1),
\vec{e}_{i_2}(\bby_2)}\,\delta a_1^{i_2}(\bby_2)\bigr\rangle \cdot
\bigl\langle\delta a_2^{j_1}(\bz_1)\,\underrightarrow{\vec{e}_{j_1}(\bz_1),
(-\vec{e}^{\,\dagger,j_2})(\bz_2)}\,\delta b_{2,j_2}(\bz_2)\bigr\rangle\cdot{}
\\
\bigl(f(\bx_1,[\ba],\boldb])\bigr)
\left\{
\begin{aligned}
\frac{\overleftarrow{\dd}}{\dd b_{i_1,\sigma_1}} \vphantom{\Bigl|}^{\lceil}
\Bigl(-\frac{\overleftarrow{\Id}}{\Id\bby_1}\Bigr)^{\sigma_1}\vphantom{\Bigr|}^{\rceil}
\langle\underline{\vec{e}_{i_1}(\bx_1)}| &\times
 |\underline{\vec{e}^{\,\dagger,i_2}(\bx_2)}\rangle\,
\vphantom{\Bigl|}^{\lceil}\Bigl(-\frac{\overrightarrow{\Id}}{\Id\bby_2}\Bigr)^{\sigma_2}\vphantom{\Bigr|}^{\rceil} \frac{\overrightarrow{\dd}}{\dd a^{i_2}_{\sigma_2}} 
\\
\frac{\overleftarrow{\dd}}{\dd a^{j_1}_{\tau_1}} \vphantom{\Bigl|}^{\lceil}
\Bigl(-\frac{\overleftarrow{\Id}}{\Id\bz_1}\Bigr)^{\tau_1}\vphantom{\Bigr|}^{\rceil}
\langle\underline{\underline{\vec{e}^{\,\dagger,j_1}(\bx_1)}}| &\times
 |\underline{\underline{\vec{e}_{j_2}(\bx_2)}}\rangle\,
\vphantom{\Bigl|}^{\lceil}\Bigl(-\frac{\overrightarrow{\Id}}{\Id\bz_2}\Bigr)^{\tau_2}\vphantom{\Bigr|}^{\rceil} (-)^{\GH{F}-1}\frac{\overrightarrow{\dd}}{\dd b_{j_2,\tau_2}} 
\end{aligned}
\right\}
\bigl(g(\bx_2,[\ba],\boldb])\bigr),
\end{multline*}
which yields zero after summation over all the (multi)\/indices.

To prove that each of the remaining two terms,\footnote{Note that
${\overrightarrow{\delta}}/{\delta\boldb}\bigl((f){\overleftarrow{\delta}}/{\delta\boldb}\bigr)=
\bigl({\overrightarrow{\delta}}/{\delta\boldb}(f)\bigr)\overleftarrow{\delta}/{\delta\boldb}$.}
\[
\frac{\overrightarrow{\delta}}{\delta\boldb}(f)\frac{\overleftarrow{\delta}}{\delta\boldb}\times
\frac{\overrightarrow{\delta}}{\delta\ba}\frac{\overrightarrow{\delta}}{\delta\ba}(g)
\qquad\text{and}\qquad
\frac{\overrightarrow{\delta}}{\delta\ba}(f)\frac{\overleftarrow{\delta}}{\delta\ba}\times
\frac{\overrightarrow{\delta}}{\delta\boldb}\frac{\overrightarrow{\delta}}{\delta\boldb}(g),
\]
cancels by itself, let us 
inspect its behaviour 
under a swap $\delta\bolds_1\rightleftarrows\delta\bolds_2$ of coefficients in the normalized test shifts.\footnote{This mechanism has already been implemented in the short proof of Lemma~\ref{LBVOperatorDifferential}, see p.~\pageref{LBVOperatorDifferential}.}

Namely, the integrand of 
the third term is~$(-)^{\GH{F}-1}$ times
\begin{multline*}
\bigl\langle\delta b_{2,j_2}(\bz_2)\,\underrightarrow{(-\vec{e}^{\,\dagger,j_2})(\bz_2),
\vec{e}_{j_1}(\bz_1)}\,\delta a_2^{j_1}(\bz_1)\bigr\rangle\cdot
\bigl\langle\delta b_{1,i_2}(\bby_2)\,\underrightarrow{(-\vec{e}^{\,\dagger,i_2})(\bby_2),
\vec{e}_{i_1}(\bby_1)}\,\delta a_1^{i_1}(\bby_1)\bigr\rangle\cdot{}
\\
\bigl(f(\bx_1,[\ba],\boldb])\bigr)
\left\{
\begin{aligned}
\frac{\overleftarrow{\dd}}{\dd b_{j_2,\tau_2}}\vphantom{\Bigl|}^{\lceil}\Bigl(-\frac{\overleftarrow{\Id}}{\Id\bz_2}\Bigr)^{\tau_2}\vphantom{\Bigr|}^{\rceil}\,\langle\underline{\vec{e}_{j_2}(\bx_1)}|
&\times
|\underline{\vec{e}^{\,\dagger,j_1}(\bx_2)}\rangle\,\vphantom{\Bigl|}^{\lceil}\Bigl(-\frac{\overrightarrow{\Id}}{\Id\bz_1}\Bigr)^{\tau_1}\vphantom{\Bigr|}^{\rceil}\frac{\overrightarrow{\dd}}{\dd a^{j_1}_{\tau_1}}
\\
\frac{\overleftarrow{\dd}}{\dd b_{i_2,\sigma_2}}\vphantom{\Bigl|}^{\lceil}\Bigl(-\frac{\overleftarrow{\Id}}{\Id\bby_2}\Bigr)^{\sigma_2}\vphantom{\Bigr|}^{\rceil}\,\langle\underline{\underline{\vec{e}_{i_2}(\bx_1)}}|
&\times
|\underline{\underline{\vec{e}^{\,\dagger,i_1}(\bx_2)}}\rangle\,\vphantom{\Bigl|}^{\lceil}\Bigl(-\frac{\overrightarrow{\Id}}{\Id\bby_1}\Bigr)^{\sigma_1}\vphantom{\Bigr|}^{\rceil}\frac{\overrightarrow{\dd}}{\dd a^{i_1}_{\sigma_1}}
\end{aligned}
\right\}
\bigl(g(\bx_2,[\ba],\boldb])\bigr).
\end{multline*}
By construction, the lower\/-\/line derivations --\,from~$\lshad\,,\,\rshad$\,-- act first on~$f$ and~$g$, and then the (graded-)\/derivations from the upper line --\,from~$\Delta$\,-- work on the respective arguments.
Now let the (multi)\/indexes be relabelled as above: $i\rightleftarrows j$, $\sigma\rightleftarrows\tau$, and
$\smash{\delta a^{i}_1\rightleftarrows\delta a^{j}_2}$, 
$\delta b_{1,i}\rightleftarrows\delta b_{2,i}$ on top of $\bby\rightleftarrows\bz$.
On the one hand, no relabelling of summation indices would affect any sum.
On the other hand, such relabelling swaps the two lines between~$f$ and~$g$, producing the minus sign factor due to the interchange of two parity\/-\/odd derivatives that fall on the first argument~$f$.
Consequently, the entire sum vanishes.

The only remaining term is processed analogously; the same relabelling of (multi)\/in\-di\-ces swaps the parity\/-\/odd derivations that act on the second argument~$g$. 
Equal to minus itself, the fourth term vanishes. 
This concludes the proof.
\end{proof}

\begin{theor}
\label{ThLaplaceOnSchouten}
Let $F$ and $G$ be two noncommutative local functionals over the infinite jet space $J^\infty\bigl(\bpiNCZO 
\bigr)$\textup{;} suppose $F$~is homogeneous. 
The Batalin\/--\/Vilkovisky Laplacian~$\Delta$ satisfies the relation
\begin{equation}\tag{\ref{EqLapSchouten}}\label{EqZimes}
\Delta\bigl(\lshad F,G\rshad\bigr)=
\lshad\Delta F,G\rshad+(-)^{|F|-1}\lshad F,\Delta G\rshad.
\end{equation}
\end{theor}
\noindent%
In other words, the operator~$\Delta$ is a graded derivation of the noncommutative variational Schouten bracket~$\lshad\,,\,\rshad$.

\begin{proof}
We prove this by induction over the number of building blocks in each argument 
of the Schouten bracket in the left hand side of~\eqref{EqZimes
}. To assert the claim in full, one reduces the set\/-\/up to integral functionals~$F$, swaps the arguments $F\rightleftarrows G$ of the Schouten bracket~$\lshad\,,\,\rshad$ by using formula~\eqref{EqSchoutenSkew}, and repeats the reasoning.\footnote{This is essential because Jacobi identity~\eqref{JacobiCommutator}, which will be used explicitly in the proof below, requires the arrangement of parentheses $(({\cdot}\times{\cdot})\times{\cdot})$ but not $({\cdot}\times({\cdot}\times{\cdot}))$ in the course of multiplication of the three functionals~$F$,\ $G$,\ and~$H$ in~\eqref{EqLeibnizFirst}.}

If~$F$ and~$G$ both belong to $\bar{H}^*\bigl(\BBT\bpiNCZO %
\bigr)$, then Lemma~\ref{LemmaBaseLapSchouten} states the assertion, which is the base of induction. 
To make an inductive step, without loss of generality let us assume that the first 
argument of $\schouten{\,,\,}$ in~\eqref{EqZimes
} is a product of \emph{two} elements 
from~$\bar{\mathfrak{N}}^n\bigl(\BBT\bpiNCZO 
\bigr)$, 
each of them containing fewer 
multiples from~$\bar{H}^*\bigl(\BBT\bpiNCZO %
\bigl)$ 
than the product. Denote such factors by~$F$ and~$G$ and suppose for definition that either of them, as well as the second argument~$H$ of the Schouten bracket, is homogeneous. 
Using Corollaries~\ref{PropSchoutenSkew} and~\ref{CorSchoutenOnProduct}, we expand
$(F\times G)\,\overleftarrow{\lshad{\cdot},H\rshad}$ and deduce that
\begin{equation}\label{EqLeibnizFirst}
\schouten{F\times G, H} = F\times\schouten{G,H} + (-)^{\GH{G}\cdot(\GH{H}-1)}\schouten{F,H}\times G.
\end{equation}
Therefore, recalling Definition~\ref{DefSchouten} of the Schouten bracket, 
we have that
\begin{multline*}
\Delta(\lshad F\times G, H\rshad)
=\Delta F\times\schouten{G,H} + (-)^{\GH{F}}\schouten{F,\schouten{G,H}} + (-)^{\GH{F}} F\times\Delta\bigl(\schouten{G,H}\bigr)
  \\
+(-)^{\GH{G}\cdot(\GH{H}-1)}\,\Bigl\{\Delta\bigl(\schouten{F,H}\bigr)\times G
  +(-)^{\GH{F}+\GH{H}-1}\schouten{\schouten{F,H},G}
  +(-)^{\GH{F}+\GH{H}-1}\schouten{F,H}\times\Delta G\Bigr\}.
\end{multline*}
Using the inductive hypothesis in the third 
and fourth 
terms of the right\/-\/hand side in the above formula, 
we continue the equality and obtain
\begin{align}
\Delta F\times\schouten{G,H} + (-)^{\GH{F}}&\schouten{F,\schouten{G,H}} 
 + (-)^{\GH{F}} \Bigl\{ F\times\schouten{\Delta G,H} + (-)^{\GH{G}-1} F\times\schouten{G,\Delta H}\Bigr\}\notag
\\
{}+(-)^{\GH{G}\cdot(\GH{H}-1)}\,\Bigl\{ &\schouten{\Delta F,H}\times G + (-)^{\GH{F}-1}\schouten{F,\Delta H}\times G \notag
\\
{}&\quad{} + (-)^{\GH{F}+\GH{H}-1}\,\Bigl[ \schouten{\schouten{F,H},G} + \schouten{F,H}\times\Delta G \Bigr] \Bigr\}.
\label{EqIndStep}
\end{align}
On the other hand, let us expand the right\/-\/hand side of~\eqref{EqZimes}, which now is
\[
\schouten{\Delta(F\times G),H} + (-)^{\GH{F}+\GH{G}-1}\schouten{F\times G,\Delta H};
\]
we recall the definition of~$\lshad\,,\,\rshad$ and we then use~\eqref{EqLeibnizFirst}.
We obtain
\begin{multline}
\schouten{\Delta F\times G+(-)^{\GH{F}}\schouten{F,G}+(-)^{\GH{F}} F\times\Delta G,H}
 + (-)^{\GH{F}+\GH{G}-1}\schouten{F\times G,\Delta H} 
\\
{}={}(-)^{\GH{G}\cdot(\GH{H}-1)}\schouten{\Delta F,H}\times G + \Delta F\times\schouten{G,H}
+(-)^{\GH{F}}\schouten{\schouten{F,G},H} 
  \\
{}\qquad +(-)^{\GH{F}}\Bigl\{ (-)^{(\GH{G}-1)\cdot(\GH{H}-1)} \schouten{F,H}\times\Delta G + F\times\schouten{\Delta G,H} \Bigr\} 
  \\
{}\qquad +(-)^{\GH{F}+\GH{G}-1}\Bigl\{(-)^{\GH{G}\cdot\GH{H}}\schouten{F,\Delta H}\times G + F\times\schouten{G,\Delta H} \Bigr\}.\label{EqClaimExpand}
\end{multline}
Comparing~\eqref{EqClaimExpand} with~\eqref{EqIndStep}, which was derived from the inductive hypothesis, 
we see that all the terms match except for
\[
(-)^{\GH{F}}\,\Bigl\{\schouten{F,\schouten{G,H}} - (-)^{(\GH{F}-1)\cdot(\GH{G}-1)}\schouten{G,\schouten{F,H}}
\Bigr\}
\]
from~\eqref{EqIndStep} versus
\[
(-)^{\GH{F}}\,\schouten{\schouten{F,G},H}
\]
from~\eqref{EqClaimExpand}. These three terms constitute $(-)^{\GH{F}}$ times the left-{} vs right\/-\/hand sides of Jacobi identity~\eqref{JacobiCommutator} for the noncommutative variational Schouten bracket. 
The balance of~\eqref{EqIndStep} and~\eqref{EqClaimExpand} completes the inductive step and concludes the proof.
\end{proof}

\begin{theor}
\label{ThBVDifferential}
The Batalin\/--\/Vilkovisky Laplacian~$\Delta$ is a differential on the space of local functionals over 
$J^{\infty}(\bpiNCZO 
)$\textup{,} \[\Delta^2=0.\]
\end{theor}
\noindent%
Summarising, the space $\overline{\mathfrak{M}}^n\bigl(\bpiNCZO 
\bigr)$ of cyclic word\/-\/valued local functionals is a (non)\/asso\-ci\-a\-tive graded\/-\/commutative BV~algebra.

\begin{proof}
We prove Theorem~\ref{ThBVDifferential} by induction over the number of building blocks 
from $\bar{H}^* 
\bigl(\BBT\bpiNCZO 
\bigr)$ 
in the argument $H \in \overline{\mathfrak{N}}^n 
\bigl(\BBT\bpiNCZO 
\bigr)$ of~$\Delta^2$. 
If $H \in \overline{H}^* 
\bigl(\bpiNCZO 
\bigr)$ itself is an integral functional, 
then by Lemma~\ref{LBVOperatorDifferential}
there remains nothing to prove. 
Suppose now that $H = F\times G$ for some $F,G \in \overline{\mathfrak{N}}^n 
\bigl(\BBT\bpiNCZO 
\bigr)$ and assume that the functional~$F$ is homogeneous. 
Then Definition~\ref{DefSchouten} 
yields that
\begin{align*}
\Delta^2&(F\times G) = \Delta\left(\Delta F\times G + (-)^{\GH{F}}\schouten{F,G} + 
 (-)^{\GH{F}}F\times\Delta G\right).
\intertext{Using Definition~\ref{DefSchouten} 
again and also Theorem~\ref{ThLaplaceOnSchouten}, 
we continue the equality:}
={}& \Delta^2F\times G + (-)^{\GH{\Delta F}}\schouten{\Delta F,G} 
+ (-)^{\GH{\Delta F}}\Delta F\times\Delta G \\
{}&{}+(-)^{\GH{F}}\schouten{\Delta F,G} + (-)^{\GH{F}}(-)^{\GH{F}-1}\schouten{F,\Delta G}
\\
{}&{}+(-)^{\GH{F}}\Delta F\times\Delta G 
+(-)^{\GH{F}}(-)^{\GH{F}}\schouten{F,\Delta G} + (-)^{\GH{F}}(-)^{\GH{F}}F\times\Delta^2G.
\end{align*}
By the inductive hypothesis, the first and last terms in the above formula vanish; 
taking into account that $\GH{\Delta F} = \GH{F}-1$ in~$\BBZ_2$, the terms with $\Delta F\times\Delta G$ 
cancel against each other, as do the terms containing $\schouten{\Delta F,G}$ and~$\schouten{F,\Delta G}$. 
The proof is complete.
\end{proof}


\begin{rem}
In the BV~context, the non\/-\/associativity of the algebra of cyclic words is a property still not a burden. To establish that the BV~Laplacian~$\Delta$ is a differential on the algebra of local functionals, we de facto proved that for any three such functionals~$F$, $G$, and~$H$ one has that $\Delta^2(F\times G\times H)=0$. In view of the non\/-\/associativity of the product~$\times$, the parentheses were arranged in the lexicographic order~$(F\times G)\times H$. This was essential for a verification of Jacobi identity~\eqref{Jacobi4Schouten}, see footnote~\ref{FootWhyFGbyHinJacobi} on p.~\pageref{FootWhyFGbyHinJacobi}. Yet because the multiplication~$\times$ is graded\/-\/commutative so that $F\times(G\times H)=(-)^{|F|\cdot(|G|+|H|)} (G\times H)\times F$, the arrangement $(\cdot\times(\cdot\times\cdot))$ is transformed into~$((\cdot\times\cdot)\times\cdot)$, which was considered before. Now relabelling the \textsl{arbitrary} functionals via $F\gets G\gets H\gets F$, we deduce that the non\/-\/associativity of operation~$\times$ in the argument of $\Delta^2(F\times G\times H)$ is not restrictive.\footnote{In the weight factor~$\exp\bigl(\tfrac{\boldsymbol{i}}{\hbar} S^\hbar\bigr)$ of the Feynman path integral, the comultiples are copies of the (quantum BV-)\/action functional~$S^\hbar$, whence the nominal non\/-\/associativity of structure~$\times$ is all the more negligible.}
\end{rem}

\begin{rem}
We conclude that the proof of all these assertions about the 
Batalin\/--\/Vil\-ko\-vi\-sky Laplacian and variational Schouten
bracket remains literally valid in the gra\-ded\/-\/commutative set-up. Indeed, when the proof is over, it suffices to let $N\mathrel{{:}{=}}0$ (so that there are no generators~$\vx_i^{\pm1}$)
and proclaim that the letters~$a^i_{\sigma}$ and~$b_{j,\tau}$ are graded\/-\/permutable; the proof itself
does not require that assumption.

Likewise, the formalism developed in Ch.~\ref{SecKinematics} survives arbitrary changes of cell decomposition for manifolds~$\bigl(M^n$,\ $\dvol({\cdot})\bigr)$, even though the tilings of newly produced spaces, whenever irregular, would make the alphabets~$\bvx^{\pm1}$ point\/-\/dependent.

We also conclude that 
by shrinking the substrate manifold $M^n$ to a point, so that $n=0$ and $N=0$, we recover the standard properties
of the parity-odd differential $\Delta_0=\overrightarrow{\dd}^2/\dd a^i\dd b_i$ and parity-odd Poisson bracket in the 
(formal non)commutative geometry of symplectic su\-per\-ma\-ni\-folds of su\-per\-di\-men\-sion~$(m|m)$. The locality of 
couplings~\eqref{EqTwoCouplings} still in force, our reasoning explains why the differentials of \emph{two} Hamiltonians and the Poisson bi\/-\/vector
are referred to \emph{the same} point when the Poisson bracket is constructed.
\end{rem}


\newpage
\section{Noncommutative variational Poisson formalism}\label{SecDynamics}
\noindent%
The noncommutative variational cotangent superspace, which we built in Ch.~\ref{SecKinematics} 
for the bundle~$\piNC$ from Ch.~\ref{SecStatic},
and the calculus of local functionals on jet spaces 
$J^{\infty}(\bpiNCZO 
)$, see Ch.~\ref{SecKinematics}, refer to the canonical \emph{symplectic}
structure 
encoded by~\eqref{EqTwoCouplings}. Let us now introduce a more narrow (sic!) class of variational noncommutative
geometries in which the \emph{Poisson} structures are defined.

\subsection{Noncommutative variational multivectors}\label{SecMultVect}
Let us recall that the notion of space of integral functionals $\bar{H}^n(\bpiNCZO 
)$ was based in 
Ch.~\ref{SecKinematics} on an obvious analytic idea to integrate the sections $\bs\in\Gamma(\piNC)$ 
over $\dvol(\bx)$ on the substrate manifold~$M^n$; the integrals take every such evaluation mapping to the cyclic word(s) written
in the edge alphabet $\bvx^{\pm1}$ (see~\eqref{EqFs} on p.~\pageref{EqFs}). When the $\BBZ_2$-\/valued parity
function was introduced, the parity\/-\/odd symbols~$\bb$ and extension $\bs^{\dagger}$ of~$\bs$ to maps defined on~$\cA^{(0|1)}$ were felt as the objects that make everything go much better as soon as one gets rid of them; we refer to Remark~\ref{RemABAB} in particular (see p.~\pageref{RemABAB}).

Taking this into account, let us describe a very different geometric approach to the use of $\BBZ_2$-parity graded
noncommutative integral functionals. Namely, we shall 
view the parity\/-\/odd symbols~$\bb$ and their derivatives as
\emph{placeholders} for (non)commutative variational covectors; such placeholders appear in the fully skew\/-\/symmetric
poly\/-\/linear maps on 
the space $\bar{H}^n(\bpiNCZO 
)$ of purely even Hamiltonian functionals. 
By making this construction precise, which forces us to narrow the class of graded-homogeneous functionals under study, we resolve the difficulty which is known from Remark~\ref{RemABAB}.

The key idea is that --~unlike it is the case for cyclic-word integral functionals of generic nature~-- the (non)commutative
variational \emph{multivectors} are organised in precisely the same way with respect to each parity-odd entry $\bb$, as long
as the shifts $\gotht$ around the circle and integrations by parts are allowed.

Let $P\in\bar{H}^n(\bpiNCZO 
)$ be a homogeneous functional of grading $|P|\mathrel{{=}{:}}k\geqslant0$.
If $k=0$, none of the cyclic words in $P$ contains any parity-odd symbols $b_{i,{\tau}}$. If $k=1$, then there is the
noncommutative linear total differential operator $A$ (that is, an operator which is polynomial in the total derivatives and
the coefficients of which are operators of left and right multiplication by functions of $\bx$ or by parity-even symbols
$\bvx^{\pm1}$ or $a^i_{\sigma}$ from the alphabet on $J^{\infty}(\bpiNCZO 
)$) such that
\[
P=\bigl(A(\bb)\bigr).
\]
Clearly, there remains nothing more to do; for the above key idea is already realized.

Suppose now $k=2$; pick \emph{one} parity\/-\/odd letter in every cyclic word of $P$ and throw all the derivations off every
such letter by using a suitable number of integrations by parts; then, if necessary, transport the letters around the
circle so that those $b_{i,\varnothing}$ stand immediately after the 
lock $\pmb{\infty}$ in the positive,
counterclockwise direction. This brings $P$ to the normal shape
\begin{equation}\label{EqNormBi}
P\cong\tfrac{1}{2}\bigl(\boldb\circ A(\bb)\bigr)\;;
\end{equation}
by construction, $A$ is the arising $(m\times m)$-size matrix linear noncommutative total differential operator of one
argument.

Arguing as above and picking \emph{some} parity-odd letter in every word of a given integral functional $P$ of grading
$k$, we transform it to the sum of cyclic words, each starting with $b_{j,\varnothing}$ for $1\le j\le m$,
\begin{equation}\label{EqNormMultVect}
P\cong\frac1{k!}\bigl(\bb\circ
A(\underbrace{\bb,\,\dots\,,\bb}_{k-1\text{ slots}})\bigr),
\end{equation}
where the noncommutative total differential operator $A$ is poly-linear in its $k-1$ arguments.%
\footnote{Of course, the notation for $A$ acting on the $m$-tuples $\bb$ is symbolic; in reality, every cyclic word of $P$
carries $k$ parity-odd entries $b_{{i_1},\varnothing}$,
$b_{{i_2},\sigma_2^{i_2}}$,\ $\dots$,\ $b_{{i_k},\sigma_k^{i_k}}$, where
$1\leqslant i_{\alpha}\leqslant m$ and the multi-indexes are word-dependent.
It is often the case that
$|\sigma^i_{\alpha}|\ne|\sigma^j_{\alpha}|$ for $i\ne j$ at some $\alpha$; for instance, recall the differential order of
entries in the matrix operator for the second Poisson structure of the renowned Boussinesq hierarchy.}

To make the construction of operator $A$ independent of our initial choice of \emph{some} parity-odd entries, let us 
analyse the properties such an operator must have. We 
consider the case $k=2$ because it will be essential in what
follows. Through the chain of integrations by parts and by carrying the parity-odd letters around the circle,
\begin{equation}\label{EqCyclicAdjoint}
P=\tfrac12(\bb\circ A(\bb))\cong\tfrac12\left((\bb)\overleftarrow{A}^{{}\,\dagger}\circ\bb\right)\sim
-\tfrac12\left(\bb\circ(\bb)\overleftarrow{A}^{{}\,\dagger}\right)\stackrel{\text{def}}=
-\tfrac12\left(\bb\circ A^{\dagger}(\bb)\right),
\end{equation}
one defines the \emph{adjoint operator} $A^{\dagger}$ that acts on its argument in the left-to-right direction.%
\footnote{Note that the \emph{left} multiplications in $A$ become the right multiplications in
$\overleftarrow{A}^{{}\!\dagger}$, and \textit{vice versa}. At the same time, the total derivative operators are reshaped by
$(\overrightarrow{\Id}/\Id\bx)^{\sigma}\circ\mapsto\circ(-\overleftarrow{\Id}/\Id\bx)^{\sigma}\mapsto
(-\overrightarrow{\Id}/\Id\bx)^{\sigma}\circ$, e.g., the adjoint to
$(aa\circ)\overrightarrow{D}_x(\,\cdot\,)(\circ a)$ is $(-\overrightarrow{D}_x)\circ((a\circ)(\,\cdot\,)(\circ aa))$.
Thirdly, the operator's matrix is transposed: $(A^{\dagger})^{ij}=(A^{ji})^{\dagger}$,
where the rightmost symbol~$\dagger$ denotes 
the adjoint of a scalar differential operator.}
The starting objects $P$ and the resulting functional are identically the same if 
we require that
\begin{equation}\label{EqASkew}
A=-A^{\dagger}.
\end{equation}
For example, let $n=1$, $m=1$ and consider $P=\tfrac12(b\circ b_x)$ with $A=\vec{\Id}/\Id x$, see~\cite{Treves2007}.

The requirements which the poly\/-\/linear operator~$A$ of $k-1$~arguments must satisfy are imposed for all $k\ge3$ in 
the same way as in~\eqref{EqCyclicAdjoint}.

In what follows, we shall consider only the grading-homogeneous functionals on $J^{\infty}(\bpiNCZO 
)$
for which the poly-linear operators $A$ are well defined, so that normalisation~\eqref{EqNormMultVect} can be attained by
starting from any parity-odd entry in every cyclic word of the functional at hand.

\begin{define}\label{DefkVector}
Homogeneous integral functionals $P\in\bar{H}^n(\bpiNCZO 
)$ of grading $k\geqslant0$ and such that either
$k\leqslant1$ or normalisation~\eqref{EqNormMultVect} is well defined are called \emph{noncommutative variational $k$-vectors}.

Let us denote by $\bar{H}^n_k(\bpiNCZO 
)\varsubsetneq\bar{H}^n(\bpiNCZO 
)$
the vector space of noncommutative variational $k$-vectors on $J^{\infty}(\bpiNCZO 
)$.
\end{define}

Note that by Remark~\ref{RemABAB}, the subspaces $\bar{H}^n_k(\bpiNCZO 
)$ do \emph{not} exhaust the
homogeneous components of grading $k$ in $\bar{H}^n(\bpiNCZO 
)$ for $k\geqslant2$.

\begin{rem}
We claim that the vector space $\bigoplus_{k\geqslant0}\bar{H}^n_k(\bpiNCZO 
)$ of all noncommutative
variational multivectors is closed under $\lshad\,,\,\rshad$, which endows it with the structure of Gerstenhaber algebra
with respect to the noncommutative variational Schouten bracket.

Definition~\ref{DefkVector} is constructive but implicit. It is instructive to see why the Schouten bracket 
$\lshad F,G\rshad$ of a $k$-vector $F$ and $\ell$-vector $G$ is a $(k+\ell-1)$-vector: this fact relies on a very 
distinguished structure --~of the local variational differential operators $\lshad F,\,\cdot\,\rshad$ or
$\lshad\,\cdot\,,G\rshad$~-- which normalization~\eqref{EqNormMultVect} provides for the geometric model of
$\lshad\,,\,\rshad$ in Remark~\ref{RemGeometricSchouten}.
\end{rem}

\begin{rem}
The price that one pays for the (non)commutative variational multivectors' realisation --~uniform 
with respect to every parity-odd entry $\bb$ under integration by parts and cyclic shifts~-- is precisely having that legal
possibility to integrate by parts. 
Yet we remember from \S\ref{SecElementsGVBV} that all such integration is postponed until
the ultimate end of every object's construction 
in the frames of the geometry of iterated variations. Therefore, the variational calculus of (non)commutative variational multivectors is
\emph{step-by-step} indeed; every intermediate object is let to exist as a well-defined notion.

For instance, Poisson bi-vectors $\bcP$ first take the Hamiltonians $F$ to the respective one-vectors $X_F$, which are also
known to us under the name of Hamiltonian evolution equations (e.\,g., of (non)commutative Korteweg\/--\/de Vries type).
In turn, the well-defined one-vector $X_F$ acts by the Schouten bracket $\lshad X_F,\,\cdot\,\rshad$ on a given
$0$-vector $H$, which defines the Poisson bracket $\{F,G\}_{\bcP}$, 
see \S\ref{SecPBr} below.

Notice that no multiplication of copies of the substrate manifold $M^n$ can be seen from this way of reasoning; in fact,
the on-the-diagonal restriction in the last phase of construction of the Schouten bracket becomes the immediate next to
the first step. This is why the Poisson framework of (non)commutative variational multivectors was not capable
of providing the intrinsic self-regularisation of 
the Batalin\/--\/Vilkovisky formalism with generic local functionals.
\end{rem}

\subsection{Derived brackets}\label{SecDB}
Let $P\in\bar{H}^n_k(\bpiNCZO 
)$ be a noncommutative variational $k$-vector. Consider $k$ integral
functionals $H_1,\,\dots\,,H_k\in\bar{H}^n_0(\bpiNCZO 
)$ of grading zero (that is, a $k$-tuple of $0$-vectors).

\begin{define}
The $k$-linear bracket 
$\{\,\cdot\,,\,\dots\,,\,\cdot\,\}_P\colon(\bar{H}^n_0\times\ldots\times\bar{H}^n_0)
(\bpiNCZO 
)\to\bar{H}^n_0(\bpiNCZO 
)$ 
is defined by the noncommutative variational $k$-vector $P$ as the derived bracket,%
\footnote{We refer to~\cite{YKSDB} for a review of the concept of derived brackets in the geometry of usual manifolds.
An algebraic classification of $N$-ary brackets is obtained in~\cite{VinNary1999}; by analysing the jet\/-\/bundle geometry in this context, in the paper~\cite{Wronskians} we developed the notion of Wronskian determinants for functions in many variables. In particular, we proved that every such structure~$W$ 
encodes a differential~$\boldsymbol{\Id}^2_W=0$.
}
\begin{equation}\label{EqDefDerivedK}
\{H_1,\,\dots\,,H_k\}_P\stackrel{\text{def}}=(-)^k\;
\underline{\lshad\ldots
\underline{\lshad P,H_1\rshad}
,\,\ldots\,,H_k\rshad}.
\end{equation}
The nested Schouten brackets are underlined
in order to emphasize that each of them produces an \emph{object}, i.\,e.\ the
noncommutative variational multivector with one parity-odd entry less than the two arguments had together. In consequence,
the integrations by parts are legitimate at every such step. This makes the Poisson formalism on 
jet spaces a science of steps and stops. 
\end{define}

\begin{example}
If $k=1$ and the noncommutative variational one-vector is the cyclic word $P=(A(\bb))$ for some total differential operator
$A$ (i.\,e.\ for a linear operator that is polynomial in the total derivatives), then
$$\{H_1\}_P=-\lshad P,H_1\rshad=(A(\delta H_1/\delta\ba)).$$
Likewise, if $k=2$ and, after a suitable number of integrations by parts, the noncommutative variational bi-vector is 
represented by the cyclic word(s) $P=\frac12(\bb\circ A(\bb))$, then it is readily seen that%
\footnote{The first equality tells us that the bracket $\{\,\cdot\,,\,\cdot\,\}_P$ which the bi-vector $P$ determines
is a bracket \emph{between} its arguments indeed.}
\begin{equation}\label{EqBr2}
\{H_1,H_2\}_P=\underline{\lshad\underline{\lshad H_1,P\rshad},H_2\rshad}\cong
\left(A\left(\frac{\delta H_1}{\delta\ba}\right)\circ
\frac{\delta H_2}{\delta\ba}\right)\sim
\left(\frac{\delta H_2}{\delta a^i}\circ A^{ij}\left(\frac{\delta H_1}{\delta a^j}\right)\right).
\end{equation}
Let us comment on every step in this construction. First, the variational one-vector $X_{H_1}$ is produced from $P$
and $H_1$; consider
$$\lshad H_1,\tfrac12(\bb\circ A(\bb))\rshad=\left(\frac{\delta H_1}{\delta\ba}\circ
\frac12\sum_{|\tau|}\left(-\frac{\overrightarrow{\Id}}{\Id\bx}\right)^{\tau}\frac{\overrightarrow{\dd}}{\dd\bb_{\tau}}
(\bb\circ A(\bb))\right).
$$
When $P=\frac12\bigl(b\circ A(\bb)\bigr)$ is varied with respect to $\bb$, the partial derivatives 
$\overrightarrow{\dd}/\dd b_{j,\tau}$ reach the first occurence $\bb_{\varnothing}$ with $\tau=\varnothing$ at once;
before they reach the argument $\bb$ of skew-adjoint operator $A$, let us integrate by parts:
$\frac12(\boldb\circ A(\underline{\bb}))\cong\frac12(-A(\boldb)\circ\underline{\bb})\sim\frac12(\underline{\bb}\circ A(\boldb))$. This shows that due to the particular structure of
bi-vectors --~if compared with generic functionals of grading two,~-- the second term doubles and absorbs $\frac12$. We get
the one-vector $(\delta H_1/\delta\ba\circ A(\bb))$; integrating by parts once again and using~\eqref{EqASkew}, we
obtain the object
$$X_{H_1}=\left(-A\left(\frac{\delta H_1}{\delta\ba}\right)\circ\bb\right).$$
Now the construction of the outer Schouten bracket in~\eqref{EqBr2} is elementary.
\end{example}

\begin{lemma}\label{LDBSkew}
Derived bracket~\eqref{EqDefDerivedK} is totally antisymmetric under permutations of its arguments:
$$\{H_{\omega(1)},\,\dots\,,H_{\omega(k)}\}_P=(-)^{\omega}\;\{H_1,\,\dots\,,H_k\}_P$$
for any $\omega\in S_k$ and any $H_1,\,\dots\,,H_k\in\bar{H}^n_0(\bpiNCZO 
)$.
\end{lemma}

\begin{rem}
The total skew-symmetry of object~\eqref{EqDefDerivedK} produced in $k$ separate steps --~with integration by parts and full
stop after each step~-- does not follow from the Jacobi identity for $\lshad\,,\,\rshad$,
which was established in Ch.~\ref{SecKinematics}. Rather, this is a manifestation
of the noncommutative variational $k$-vectors' intrinsic property to be structurally identical with respect to every two
graded entries $\bb$.
\end{rem}

\begin{proof}[Sketch of the proof]
It suffices to show that the derived bracket $\{\,\cdot\,,\,\dots\,,\,\cdot\,\}_P$ changes its sign under a swap of two
consecutive arguments $H_i$ and $H_{i+1}$\;:
$$\ldots\underline{\lshad\underline{\lshad Q,H_i\rshad},H_{i+1}\rshad}\ldots\cong
-\ldots\underline{\lshad\underline{\lshad Q,H_{i+1}\rshad},H_i\rshad}\ldots\ .$$
Consider the noncommutative variational multivector's necklace $Q$ and mark, by using $\otimes$ and $\oplus$, two
parity\/-\/odd entries $\bb$ (e.\,g., the two \emph{consecutive} ones for the sake of clarity), see the figure on the facing 
page.

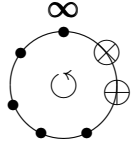
\begin{center}
\unitlength=1mm
\linethickness{0.4pt}
\begin{picture}(22.00,18.00)(5,0)
\put(15.00,10.00){\circle{14.00}}
\put(15.00,10.00){\makebox(0,0)[cc]{{\large$\circlearrowleft$}}}
\put(15.00,17.00){\circle*{1.33}}
\put(15.00,19.00){\makebox(0,0)[cb]{$\pmb{\infty}$}}
\put(9.33,14.33){\circle*{1.33}}
\put(8.33,7.33){\circle*{1.33}}
\put(12.00,3.67){\circle*{1.33}}
\put(18.00,3.67){\circle*{1.33}}
\put(20.67,14.33){\makebox(0,0)[cc]{{\large$\otimes$}}}
\put(22,9){\makebox(0,0)[cc]{{\large$\oplus$}}}
\end{picture}
\end{center}

\noindent%
This object's inner Schouten bracket with $H_i$ does basically the following: normalisation~\eqref{EqNormMultVect} throws all the
derivatives off the entry $\otimes$ and implants $\delta H_i/\delta\ba$ in its stead (the normalisation does exactly the same with every other entry $\bb$ by the definition of multivector, but let us focus on the term such that the
variation $\delta H_i/\delta\ba$ hits $\otimes$). Now reshape this output by making $\oplus$ free of derivatives falling
on it. Note that this session of integrations by parts again amounts to bringing  the multivector to normalized
shape~\eqref{EqNormMultVect}, --~only the neighbouring entry $\otimes$  is occupied now by $\delta H_i/\delta\ba$, not by
$\bb$. The outer Schouten bracket installs $\delta H_{i+1}/\delta\ba$ at $\oplus$ (or at any other parity-odd entry; we consider just one term, for definition).

\enlargethispage{\baselineskip}
On the other hand, consider the very same scenario of putting $\delta H_i/\delta\ba$ for $\otimes$ and
$\delta H_{i+1}/\delta\ba$ for $\oplus$, done in the reverse order. To reach $\oplus$ first in the construction of
(now, inner) Schouten bracket, the derivation $\overleftarrow{\dd}/\dd\bb$ has to overtake $\otimes$ currently occupied
by the parity-odd placeholder $\bb$; this overtake yields the sought\/-\/for minus sign. The variation
$\delta H_{i+1}/\delta\ba$ pasted for $\oplus$, we cast all the derivatives off the still-unused slot $\otimes$, 
leave $\delta H_i/\delta\ba$ there, and integrate by parts back, to 
isolate $\delta H_{i+1}/\delta\ba$ in the socket~$\oplus$. It is readily seen that the two algorithms produce the identical portraits of letters and derivatives, yet those two differ by the sign factor.
\end{proof}

\begin{rem}
Continuing this line of reasoning, we conclude that for a given noncommutative variational $k$-vector $P$, the value
$\{H_1,\,\dots\,,H_k\}_P$ of derived bracket~\eqref{EqDefDerivedK} at $k$ arguments
$H_1,\,\dots\,,H_k\in\bar{H}^n_0(\bpiNCZO 
)$ is equivalent, up to integration by parts, to the $0$-vector
\begin{equation}\label{EqPatExact}
(-)^{\frac{k(k-1)}2}\cdot\frac1{k!}\sum_{\omega\in S_k}(-)^{\omega}
\left(\frac{\delta H_{\omega(1)}}{\delta\ba}
\circ A\left(\frac{\delta H_{\omega(2)}}{\delta\ba},\,\dots\,,\frac{\delta H_{\omega(k)}}{\delta\ba}\right)
\right)
\cong\{H_1,\,\dots\,,H_k\}_P,
\end{equation}
where the alternating sum runs through the entire permutation group $S_k$; note that it is the parity-even arguments $H_i$
but not the slots for them which are shuffled.
\end{rem}

Observation~\eqref{EqPatExact} allows us to extend the \emph{mapping} $P$ from the geometry of exact (non)commutative
variational covectors $\delta H_i/\delta\ba$,
$$P\left(
{\delta H_1}\bigr/{\delta\ba},\,\dots\,,
{\delta H_k}\bigr/{\delta\ba}\right)\stackrel{\text{def}}=
\{H_1,\,\dots\,,H_k\}_P,$$
to $k$-tuples of arbitrary variational covectors 
$\bp_i=(p_{i,\alpha}\circ\delta a^{\alpha})$. 
(Let us think of the 
{variational covectors} $(\bp\circ\delta\ba)=\bigl(p_\alpha\bigl(\bx,\bvx^{\pm1},[\ba]\bigr)\circ\delta a^\alpha\bigr)$ on $J^\infty(\piNC
)$ as of (the formal sums of) necklaces equipped with the extra earrings~$\delta a^\alpha$, by which those cyclic words are handled.)
The case $k=1$ with $P(\bp_1)\mathrel{{:}{=}}(A(\bp_1))$ is elementary; 
for $k\geqslant2$, we put%
\footnote{\label{FootEarrings}%
The isomorphism $V^{\dagger}\simeq T_{\ba^{\dagger}}V^{\dagger}$ is used here to convert the placeholders $\bb$ for $\bp_i$
into the virtual offsets $\sum\limits_{\alpha=1}^m1\cdot\vec e^{{}\,\dagger,\alpha}$. The absorption of each argument
$\bp_i$ then goes closely to the lines of geometric construction of the Schouten bracket, see 
Remark~\ref{RemGeometricSchouten} on p.~\pageref{RemGeometricSchouten}.%
}
\begin{equation}\label{EvalkVector}
P(\bp_1,\,\dots\,,\bp_k)\mathrel{{:}{=}}(-)^{\frac{k(k-1)}2}\cdot\frac1{k!}\sum_{\omega\in S_k}(-)^{\omega}
\left(\bp_{\omega(1)}\circ A(\bp_{\omega(2)},\,\dots\,,\bp_{\omega(k)})\right).
\end{equation}
However, generic variational covectors, not necessarily exact, will not be studied in particular in what follows -- rather,
the converse can be assumed in view of the Substitution Principle. 


\label{SecSubst}
\begin{theor}[The Substitution Principle]\label{ThSubstPrinciple}
Suppose that a tuple of identities 
\[
\boldsymbol{I}\bigl((\bx,\bvx^{\pm1}),[\ba],
\bigl[\bp_1(\bx,\bvx^{\pm1})\bigr],\ldots,
\bigl[\bp_k(\bx,\bvx^{\pm1})\bigr]\bigr)\equiv0
\] 
holds on~$J^\infty(\piNC
)$ 
for every $k$-\/tuple of noncommutative variational \textup{(}co\textup{)}vectors the coefficients $p_{i,\alpha}(\bx,\bvx^{\,\pm1})$ 
of which 
depend only on points~$\bx\in M^n$
and letters from the edge alphabet~$\bvx^{\,\pm1}$.
Then the identities in total derivatives, 
\[
\boldsymbol{I}\bigl((\bx,\bvx^{\pm1}),[\ba],
\bigl[\bp_1\bigl(\bx,\bvx^{\pm1},[\ba]\bigr)\bigr],\ldots,
\bigl[\bp_k\bigl(\bx,\bvx^{\pm1},[\ba]\bigr)\bigr]\bigr)\equiv 0,
\]
viewed as identities with respect to~$\bp_i$, are valid 
on~$J^\infty(\piNC
)$ for all 
\textup{(}co\textup{)}\/vectors~$\bp_i$ depending not only 
on~$\bx$ and~$\bvx^{\pm1}$ but also admitting 
arbitrary, finite differential order dependence on the jet letters~$\ba_\sigma$, $|\sigma|<\infty$.
\end{theor}

\begin{rem}
At this moment it is legitimate to view the variational (co)\/vectors~$\bp_i=(p_{i,\alpha}\circ\delta a^\alpha)$ as bare collections of their indexed open\/-\/word components~$p_{i,\alpha}$ that are already built into the identities~$\boldsymbol{I}$. 
We emphasize that, unlike it is the case studied in~\S\ref{SecAlgebra} --\,the cyclic words in~$\cA$ do not carry any marked point,\,-- the earrings $\dd/\dd \ba_\sigma$ and~$\delta\ba$ are the only places where the (co)\/vectors can be unlocked. 
\end{rem}

\begin{cor}
If, under the assumptions of Theorem~\textup{\ref{ThSubstPrinciple},}
the identities in total derivatives~$\boldsymbol{I}\bigl(\bx,\bvx^{\pm1},[\ba],[\bp_i]\bigr)\equiv0$ 
with respect to~$\bp_1$,\ $\ldots$,\ $\bp_k$ hold on~$J^\infty(\piNC
)$ for every $k$-\/tuple of \emph{exact} variational covectors $\bp_i=(
{\delta\cH_i}/\delta\ba\circ\delta\ba)$ which are obtained by variation of arbitrary linear integral functionals $\cH\in\bar{H}^n(\piNC
)$, then these identities hold for \emph{all} covectors~%
$\bp_i$, 
i.e.\ not necessarily exact.
\end{cor}

\noindent%
Indeed, it is always possible to represent locally an $(\bx,\bvx^{\pm1})$-\/dependent cyclic word 
$\sum\nolimits_{\alpha=1}^m \bigl(p_{i,\alpha}(\bx,\bvx^{\pm1})\circ\delta a^\alpha\bigr)$ as the 
variation $
{\delta}\cH$ of the functional 
$\sum\nolimits_{j=1}^m \int \bigl(p_{i,\alpha}(\bx,\bvx^{\pm1})\circ a^\alpha\,\dvol(\bx)\bigr)$ and then apply Theorem~\ref{ThSubstPrinciple}.

\begin{proof}[Proof of Theorem~\textup{\protect\ref{ThSubstPrinciple}}]
For the sake of brevity, let each variational noncommutative \emph{co}vec\-tor~$\bp_i$ consist of just one word written in the alphabet of~$J^\infty(\piNC
)$.
The crucial idea is that the position of the locks~$\delta\ba$ is fixed on the circles which carry the words~$\bp_i$. This means that, whenever one declares an arbitrary differential dependence of~$\bp_i$ on~$\ba$, the words~$\boldsymbol{I}$ in principle lengthen but still, in the course of multiplications~$\times$ within the identities, each~$\bp_i$~is never torn in between any 
consecutive pair of letters~$\ba$. Namely, 
during the evaluation of~$\boldsymbol{I}$ at the words~$\bp_i$ those are unlocked, the letters 
and the words' overall coefficients depending on~$\bx$ are then stretched to  open strings (ordered counterclockwise). 
These strings of symbols are pasted into~$\boldsymbol{I}$ without splitting, i.e., the adjacent letters of~$\bp_i$ never become separated by any other symbols.%
\footnote{
This scenario is realised irrespectively of 
presence or absence of letters~$\ba$'s on the necklaces~$\bp_i$, 
which is in contrast with formula~\eqref{EqDefMult}.}

Total derivatives~\eqref{EqDefTD} 
then work according to their definition: under a restriction of~$\boldsymbol{I}$ (hence of all~$\bp_i$) to the jet 
of a 
mapping~$\ba=\bs(\bx,\bvx^{\pm})$, each symbol~$a^j$ is replaced with the respective sum of open strings~$s^j(\bx,\bvx^{\pm1})$ so that 
derivations~\eqref{EqDefTD} which act on~$
\ba_\sigma
$ occurring anywhere (either in~$\bp_i$ or in~$\boldsymbol{I}$ if
the identities explicitly depend on~$[\ba]$) 
then reduce to the 
derivations~$
{\dd}/\dd x^i$ 
of real\/-\/valued functions defined at~$\bx\in U\subseteq M^n$.
By the initial assumption of the theorem, its assertion is valid 
for all 
strings written in the basic alphabet~$(\bx,\bvx^{\pm1})$ that replace
the entries~$\bp_i$ in~$\boldsymbol{I}$. 
We conclude that the identities $\boldsymbol{I}\equiv0$ hold on~$J^\infty(\piNC
)$ for the full set of arguments of the (co)\/vec\-tors.%
\footnote{
One does not even have to postulate that the 
mappings~$\ba=\bs(\bx,\bvx^{\pm1})$ inserted in the explicit dependence of~$\boldsymbol{I}$ on~$[\ba]$ coincide with the mappings 
now standing for~$\ba$ in the implicit dependence $\bigl[\bp_i\bigl(\bx,\bvx^{\pm1},[\ba]\bigr)\bigr]$.
}
\end{proof}

\begin{rem}
The proof remains literally valid in the case of (evolutionary)
vector fields instead of variational \emph{co}vectors. This would be 
important for the description of variational noncommutative symplectic structures.
At the same time, the proof reveals \emph{why} this noncommutative phrasing 
of the Substitution Principle does \emph{not} hold for arbitrary 
cyclic words
$\bp_i\bigl(\bx,\bvx^{\pm1},[\ba]\bigr)$ of unspecified 
nature.
\end{rem}

\begin{rem}
Attempts to define the (non)commutative variational Schouten bracket of multivectors via a recursive procedure that involves
the use of the two arguments' values at test covectors are sometimes practiced in the literature (see discussion in~\cite{RingersProtaras} and references therein).
\end{rem}

\begin{OpenProblem}
Is there a way to detect that a given (non)\/commutative variational $0$-\/vector $H\in\bar{H}^n_0\bigl(\bpiNCZO 
\bigr)$ is the value of a (non)\/commutative variational $k$-\/vector at~$k$ 
zero\/-\/vectors\,?
\end{OpenProblem}


\subsection{Noncommutative variational Poisson structures}\label{SecPBr}
Now we analyse the construction of noncommutative variational Poisson brackets, recalling and re\/-\/proving several important facts
~--- here, under the coarse assumption of cyclic invariance (e.g., the Helmholtz lemma reveals  yet another mechanism for the differentials to anticommute).

\begin{rem}
Although the formalism is based on the noncommutative variational \emph{symplectic} geometry from Ch.~\ref{SecKinematics}, the presence of differential operators~$A$ in the definition of the Poisson bracket~$\{\,,\,\}_{\bcP}$ as derived with respect to a given Poisson bi\/-\/vector~$\bcP$, see~\eqref{EqDefDerivedK}, usually makes such brackets degenerate. Their Casimirs, forming the zeroth Poisson cohomology group with respect to~$\dd_{\bcP_1}=\lshad\bcP_1,\,\cdot\,\rshad$, start the Magri scheme for systems possessing the bi\/-\/Hamiltonian structures~$(\bcP_1$,\ $\bcP_2)$, 
see~\S\ref{SecMagri} 
and~\cite{DeSoleKacCMP2012,DeSoleKacJap2013}.  
\end{rem}

\subsubsection{The definition of Poisson bracket}
Consider a noncommutative variational bi\/-\/vector~$\bcP$ and let $H_1$,\ $H_2$,\ $H_3\in\bar{H}^n_0\bigl(\bpiNCZO 
\bigr)$ be any three noncommutative variational $0$-\/vectors.

\begin{define}
Bi\/-\/linear, skew\/-\/symmetric derived bracket~\eqref{EqBr2},
\[
\{H_i,H_j\}_{\bcP}=\underline{\lshad \underline{\lshad H_i,\bcP\rshad},
H_j\rshad},\qquad 1\leqslant i<j\leqslant 3,
\]
is called the noncommutative variational \emph{Poisson bracket} if it satisfies 
Jacobi identity,
\begin{equation}\label{JacPB}
\{\{H_1,H_2\}_{\bcP},H_3\}_{\bcP} +
\{\{H_2,H_3\}_{\bcP},H_1\}_{\bcP} +
\{\{H_3,H_1\}_{\bcP},H_2\}_{\bcP} \cong 0
\end{equation}
for all $H_1$,\ $H_2$,\ $H_3\in\bar{H}^n_0\bigl(\bpiNCZO 
\bigr)$, which are then called the \emph{Hamiltonians}.

If identity~\eqref{JacPB} holds, the noncommutative variational bi\/-\/vector $\bcP=\tfrac{1}{2}\bigl(\boldb\circ A(\boldb)\bigr)$ is called \emph{Poisson};
the skew\/-\/adjoint noncommutative linear operator~$A$ in total derivatives is then called a \emph{Hamiltonian operator}, and the noncommutative variational one\/-\/vectors $X_{H_i}\eqdef\lshad\bcP,H_i\rshad$ are the \emph{Hamiltonian one\/-\/vectors} (or \emph{one\/-\/vector fields}) specified by their Hamiltonians~$H_i$ and the Poisson bi\/-\/vector~$\bcP$.
\end{define}

\begin{criterion}\label{ThCriterionPoisson}
A noncommutative variational bi\/-\/vector~$\bcP$ is Poisson \textup{(}i.e.\ the derived bracket~$\{\,,\,\}_{\bcP}$ satisfies Jacobi identity~\eqref{JacPB}\textup{)} if the bi\/-\/vector~$\bcP$ satisfies the classical master\/-\/equation
\begin{equation}\label{ClassMaster}
\lshad\bcP,\bcP\rshad\cong 0\in\bar{H}^n_3\bigl(\bpiNCZO 
\bigr).
\end{equation}
The bi\/-\/vector~$\bcP$ is Poisson only if the value of $\lshad\bcP,\bcP\rshad$ at any triple $H_1$,\ $H_2$,\ $H_3$~of Hamiltonians is cohomologically trivial\textup{:}
\[
\lshad\bcP,\bcP\rshad(H_1,H_2,H_3)\cong 0\in\bar{H}^n_0\bigl(\bpiNCZO 
\bigr).
\]
\end{criterion}
\noindent%
The assertion is aimed to emphasize that the Poisson bi\/-\/vectors are the primary objects, whereas the Poisson brackets are the derived structures.

\begin{lemma}\label{LZeroQZeroValue}
If a noncommutative variational $k$-\/vector~$\bcQ$ represents the class of zero in~$\bar{H}^n_k\bigl(\bpiNCZO 
\bigr)$, then, $\bcQ$~viewed 
as the map $\bigl(\bar{H}^n_0\times\ldots\times\bar{H}^n_0\bigr)
\bigl(\bpiNCZO 
\bigr) \to \bar{H}^n_0\bigl(\bpiNCZO 
\bigr)$,
its value $\bcQ\bigl(\delta H_1/\delta\ba,\ldots,\delta H_k/\delta\ba\bigr)=
\{H_1,\ldots,H_k\}_{\bcQ}$ is cohomologically trivial for every $k$-\/tuple of the arguments~$H_1$,\ $\ldots$,\ $H_k\in\bar{H}^n_0\bigl(\bpiNCZO 
\bigr)$.
\end{lemma}

\begin{proof}[Sketch of the proof]
Indeed, whenever the cyclic word~$\bcQ=\Id_h\bcR(\boldb,\ldots,\boldb)$  carrying $k$~parity\/-\/odd entries~$\boldb$ is exact with respect to the lift~$\Id_h$ of the de Rham differential for~$M^n$ onto $J^\infty\bigl(\bpiNCZO 
\bigr)$, so is every term --\,in the sum over the $|S_k|=k!$~ways to permute the arguments~$H_1$,\ $\ldots$,\ $H_k$ by using~$\omega\in S_k$\,-- obtained by pasting whatever open string $\delta H_{\omega(i)}/\delta a^j$ of parity\/-\/even symbols instead of the $i$th copy of the symbol~$b_j$.
\end{proof}

\begin{rem}
\label{FootZeroQZeroValue}%
The gap between necessity, 
\begin{itemize}
\item a variational bi\/-\/vector~$\bcP$ is Poisson only if all the values of the variational  tri\/-\/vector $\lshad\bcP,\bcP\rshad$ are trivial in~$\bar{H}^n_0\bigl(\bpiNCZO\bigr)$,
\end{itemize}
and sufficience, 
\begin{itemize}
\item a variational bi\/-\/vector~$\bcP$ is Poisson if the variational tri\/-\/vector $\lshad\bcP,\bcP\rshad$ itself is trivial in the respective horizontal cohomology group~$\bar{H}^n_3\bigl(\bpiNCZO\bigr)\neq\bar{H}^n_0\bigl(\bpiNCZO\bigr)$,
\end{itemize}
is the statement 
that, whenever the value $\bcQ(\delta H_1/\delta\ba$,\ $\ldots$,\ $\delta H_k/\delta\ba)$ of a (non)\/com\-mu\-ta\-tive variational $k$-\/vector~$\bcQ$ at every $k$-\/tuple of exact variational covectors~$\delta H_i/\delta\ba$ is cohomologically trivial in~$\bar{H}^n_0\bigl(\bpiNCZO 
\bigr)$, the $k$-\/vector~$\bcQ$ itself is cohomologically trivial in~$\bar{H}^n_k\bigl(\bpiNCZO 
\bigr)$.
This claim proven, Criterion~\ref{ThCriterionPoisson} (and Lemma~\ref{LZeroQZeroValue}) would convert into an equivalence.
\end{rem}

\begin{lemma}\label{LCriterionStarShaped}
In fact, this is true,
\[
\bcP~\text{Poisson}\quad\Longleftrightarrow\quad\lshad\bcP,\bcP\rshad\cong0,
\]
over topologically trivial, star\/-\/shaped domains~$\subseteq M^n$.
\end{lemma}

Indeed, under the trivial topology assumption, the homotopy procedure (e.g., see~\cite{Olver1993} or~\cite{TwelveLectures})
in the constructive proof of the Poin\-ca\-r\'e\ lemma works both on the base, which we denote still by~$M^n$, and in the topologically trivial fibres of the Whitney sum of the (non)commutative bundle~$\piNC$ and $k$~copies of its dual~$\widehat{\pi}_{\nC}$.

\begin{proof}[Sketch of the proof]
Consider not the bundle~$\piNC$ such that $\bcQ\bigl(\delta H_1/\delta\ba$,\ $\ldots$,\ $\delta H_k/\delta\ba\bigr)\in\bar{H}^n(\piNC)$ but introduce the Whitney sum
$\piNC \mathbin{{\times}_{\MnnC}} \widehat{\pi}_{\nC} \mathbin{{\times}_{\MnnC}} \cdots \mathbin{{\times}_{\MnnC}} \widehat{\pi}_{\nC}$ with $k$~copies of the dual bundle (with the respective fibre variables $\bp_1$,\ $\ldots$,\ $\bp_k$ that imitate the variational covectors).
Now we have that the $n$th~degree horizontal cohomology classes $\bcQ(\bp_1$,\ $\ldots$,\ $\bp_k)$ are $k$-\/linear and totally skew\/-\/symmetric w.r.t.\ the new covector variables~$\bp_\alpha$. All these classes are known to be trivial by our initial assunmption.
The homotopy procedure then yields a $k$-\/linear w.r.t.\ $\bp_1$,\ $\ldots$,\ $\bp_k$, totally skew\/-\/symmetric horizontal $(n-1)$-\/form~$\bcR$ such that $\bcQ=\Id_h(\bcR)$ for all sections $\bp_\alpha(\bx,\bvx^{\,\pm1})$ of~$\widehat{\pi}_{\nC}$.
The Substitution Principle now works. Finally, replacing the $k$-\/linear skew terms over the Whitney sum by the variational $k$-\/vectors (with $k$~copies of the parity\/-\/odd~$\boldb$) over the superbundle~$\bpiNCZO$ is technical.
\end{proof}



\subsubsection{Noncommutative differential forms}
To approach the proof of Criterion~\ref{ThCriterionPoisson},
let us recall several classical 
structures that appear on the infinite jet spaces~$J^\infty\bigl(\piNC 
\bigr)$: in particular, in the context of the Vinogradov $\cC$-\/spectral sequence~\cite{VinogradovCspectral}.

By definition, put
\[
\vec{\dd}^{\,(\ba)}_{\vph(\bx,\bvx^{\pm1},[\ba])} = \sum_{i=1}^m \sum_{|\sigma|\geqslant0}
\Bigl((\vph^i)\Bigl(\frac{\overleftarrow{\Id}}{\Id\bx}\Bigr)^\sigma\Bigr)
\bigl(\bx,\bvx^{\pm1},[\ba]\bigr)\circ\frac{\overrightarrow{\dd}}{\dd a^i_\sigma}.
\]
It is readily seen that these \emph{evolutionary derivations} commute with the total derivatives on~$J^\infty\bigl(\piNC 
\bigr)$:
\[
\bigl[\vec{\dd}^{\,(\ba)}_{\vph},\vec{\Id}/\Id x^k\bigr] =0\qquad \text{for all~$k=1,\ldots,n$.}
\]
Consequently, for any operator~$A$ in total derivatives we have that
\[
\vec{\dd}^{\,(\ba)}_{\vph}\bigl(A(\bp)\bigr) =
\bigl(\vec{\dd}^{\,(\ba)}_{\vph}(A)\bigr)(\bp) + A\bigl(\vec{\dd}^{\,(\ba)}_{\vph}(\bp)).
\]
Next, define the \emph{linearization} $\ell^{(\ba)}_{\bp}$ of an object~$\bp$ over~$J^\infty\bigl(\piNC 
\bigr)$ by setting
\[
(\vph)\overleftarrow{\ell}^{(\ba)}_{\bp}=\overrightarrow{\dd}^{\!(\ba)}_{\!\vph}(\bp)
\]
whenever the right\/-\/hand side is well defined.

Thirdly, for each value of the index~$i$ running from~$1$ to~$m$ and for every multi\/-\/index~$\sigma$ let us introduce the symbol~$\Id a^i_\sigma$. Now define the \emph{Cartan differential} $\Id_\cC\colon a^i_\sigma\mapsto \Id a^i_\sigma$, $\Id a^i_\sigma\mapsto 0$, also setting its action equal to zero on~$\bx$ and~$\bvx^{\pm1}$ and postulating that~$\Id_\cC$ is a graded derivation.
By construction, let the differential~$\Id_\cC$ be correlated with other structures on~$J^\infty\bigl(\piNC 
\bigr)$ in the standard way: e.g., set $\vec{D}_{x^k}(\Id a^i_\sigma)=\Id a^i_{\sigma\cup\{k\}}$.

Let us explain what it means that the symbols~$\Id a^i_\sigma$ and~$\Id a^j_\tau$ ``anticommute.'' The key idea is that the precedence\/-\/succedence
relation of such symbols in a given cyclic word manifests that circle's orientation, which is provided by construction.

Consider a cyclic word that carries \emph{one} symbol~$\Id a^i_\sigma$; the word thus acquires a marked point. The derivation~$\Id_\cC$ acts on (the rest of) the word by starting at~$\Id a^i_\sigma$ and processing the letters~$a^j_\tau$ by going in the positive direction. We say that all the symbols~$\Id a^j_\tau$, newly produced by~$\Id_\cC$ from such~$a^j_\tau$ are \emph{succedent} with respect to the mark~$\Id a^i_\sigma$; in turn, the old symbol~$\Id a^i_\sigma$ is \emph{precedent} for each new object~$\Id a^j_\tau$. To change this precedence\/-\/succedence relation~$\Id a^i_\sigma\prec\Id a^j_\tau$ but still let the circle's orientation stay intact, the object~$\Id a^j_\tau$ is  proclaimed the new marked point~--- so that $\Id a^i_\sigma$~now 
\emph{succeeds} it with respect to the positive order of letters written along the oriented circle. By definition, such involution of the relative order~$\prec$ of the two symbols, $\Id a^i_\sigma$ and~$\Id a^j_\tau$, produces the factor $-1$ in front of the cyclic word that carries both of them. Clearly,~$\Id_\cC^2=0$.

\begin{lemma}[Helmholtz]\label{LHelmholtz}
The linearization $\vec{\ell}^{\,(\ba)}_{\delta H/\delta\ba}$ of an element in the image of variational derivative~$\delta/\delta\ba$ is self\/-\/adjoint\textup{:}
\begin{equation}\label{EqHelmholtz}
\vec{\ell}^{\,(\ba)}_{\delta H/\delta\ba} =
\vec{\ell}^{\,(\ba)\,\dagger}_{\delta H/\delta\ba}.
\end{equation}
\end{lemma}
\noindent%
Note that this half of 
Helmholtz' \emph{criterion} does not refer to the topology of the set\/-\/up.

\begin{proof}
Let $H$~be a noncommutative variational $0$-\/vector. Up to an integration by parts, we have that $\Id_\cC H\cong\bigl(\Id\ba\circ\delta H/\delta\ba\bigr)$.
By the above,
\[
0=\Id_\cC^2(H)\cong\bigl(\Id\ba\circ \overrightarrow{\ell}^{(\ba)}_{\delta H/\delta\ba}(\underline{\Id\ba})\bigr) \cong
\bigl((\Id\ba)\overleftarrow{\ell}^{(\ba)\,\dagger}_{\delta H/\delta\ba}\circ\underline{\Id\ba}\bigr) \sim
-\bigl(\underline{\Id\ba}\circ \overrightarrow{\ell}^{(\ba)\,\dagger}_{\delta H/\delta\ba}(\Id\ba)\bigr),
\]
whence~\eqref{EqHelmholtz}.
\end{proof}

\subsubsection{Proof of Criterion~\protect\ref{ThCriterionPoisson}}
First, let us recall the renowned cancellation mechanism in the left\/-\/hand side of Jacobi identity~\eqref{JacPB}. By definition, put $\bp_i=\delta H_i/\delta\ba$ for the three Hamiltonians.
Integrating by parts in the inner and outer Poisson brackets in~\eqref{JacPB} and using formula~\eqref{EqBr2}, we get
\begin{align}
{}&\vec{\dd}^{\,(\ba)}_{A(\bp_1)}\bigl(\bp_2\circ A(\bp_3)\bigr) +
\vec{\dd}^{\,(\ba)}_{A(\bp_2)}\bigl(\bp_3\circ A(\bp_1)\bigr) +
\vec{\dd}^{\,(\ba)}_{A(\bp_3)}\bigl(\bp_1\circ A(\bp_2)\bigr)\notag \\
{}&{}\quad{}=
\bigl(\vec{\dd}^{\,(\ba)}_{A(\bp_1)}(\bp_2)\circ A(\bp_3)\bigr) +
 \bigl(\bp_2\circ\vec{\dd}^{\,(\ba)}_{A(\bp_1)}(A)(\bp_3)\bigr) -
 \bigl(A(\bp_2)\circ\vec{\dd}^{\,(\ba)}_{A(\bp_1)}(\bp_3)\bigr)\notag \\
{}&{}\qquad{}
+\bigl(\vec{\dd}^{\,(\ba)}_{A(\bp_2)}(\bp_3)\circ A(\bp_1)\bigr) +
 \bigl(\bp_3\circ\vec{\dd}^{\,(\ba)}_{A(\bp_2)}(A)(\bp_1)\bigr) -
 \bigl(A(\bp_3)\circ\vec{\dd}^{\,(\ba)}_{A(\bp_2)}(\bp_1)\bigr)\notag \\
{}&{}\qquad{}
+\bigl(\vec{\dd}^{\,(\ba)}_{A(\bp_3)}(\bp_1)\circ A(\bp_2)\bigr) +
 \bigl(\bp_1\circ\vec{\dd}^{\,(\ba)}_{A(\bp_3)}(A)(\bp_2)\bigr) -
 \bigl(A(\bp_1)\circ\vec{\dd}^{\,(\ba)}_{A(\bp_3)}(\bp_2)\bigr).\label{EqExpandJac}
\end{align}
Applying Lemma~\ref{LHelmholtz} to the variational covectors~$\bp_i=\delta H_i/\delta\ba$ as follows,
\begin{multline*}
\bigl(\vec{\dd}^{\,(\ba)}_{A(\bp_1)}(\bp_2)\circ A(\bp_3)\bigr) \eqdef
\bigl(\vec{\ell}^{\,(\ba)}_{\bp_2}(A(\bp_1))\circ A(\bp_3)\bigr) =
\bigl(\vec{\ell}^{\,(\ba)\,\dagger}_{\bp_2}(A(\bp_1))\circ A(\bp_3)\bigr) \\
{}\cong
\bigl(A(\bp_1)\circ\vec{\ell}^{\,(\ba)}_{\bp_2}\bigl(A(\bp_3)\bigr)\bigr)\eqdef
\bigl(A(\bp_1)\circ\vec{\dd}^{\,(\ba)}_{A(\bp_3)}(\bp_2)\bigr),
\end{multline*}
we conclude that it is only the second column which survives the cancellation in~\eqref{EqExpandJac}. The left\/-\/hand side of Jacobi identity thus equals
\begin{equation}\label{EqLHSJac}
\Bigl(\tfrac{\delta H_1}{\delta\ba}\circ \vec{\dd}^{\,(\ba)}_{A(\delta H_3/\delta\ba)}(A)\bigl(\tfrac{\delta H_2}{\delta\ba}\bigr)\Bigr) +
\text{cyclic permutations}.
\end{equation}

On the other hand, consider the bi\/-\/vector~$\bcP=\tfrac{1}{2}\bigl(\boldb\circ A(\boldb)\bigr)$ and construct
\[
\lshad\bcP,\bcP\rshad \cong
\Bigl(\bigl(\boldb\circ A(\boldb)\bigr)\Bigl(\frac{\overleftarrow{\dd}}{\dd \ba_\sigma}\circ\Bigl(\frac{\overrightarrow{\Id}}{\Id\bx}\Bigr)^\sigma
\bigl(A(\boldb)\bigr)\Bigr)\Bigr);
\]
the right\/-\/hand side contains, for every multi\/-\/index~$\sigma$, the derivation that pastes its coefficient for each~$a^i_\sigma$ occurring in the coefficients of operator~$A$ within~$\bigl(\boldb\circ A(\boldb)\bigr)$.

The only thing which the evaluation of $\lshad\bcP,\bcP\rshad$ at~$H_1$,\ $H_2$,\ and~$H_3$ does,
\[
\lshad\bcP,\bcP\rshad\,\bigl(
{\delta H_1}/{\delta\ba},
{\delta H_2}/{\delta\ba},
{\delta H_3}/{\delta\ba}\bigr) =
(-)^3\,\underline{\lshad\underline{\lshad\underline{\lshad\underline{\lshad\bcP,\bcP\rshad},
H_1\rshad}, H_2\rshad}, H_3\rshad},
\]
is the speading of variational derivatives~$\delta H_i/\delta\ba$ over the three slots~$\boldb$ in the tri\/-\/vector~$\lshad\bcP,\bcP\rshad$.
In view of 
evaluation's total skew\/-\/symmetry (see Lemma~\ref{LDBSkew}),
it is enough to sum up over the cyclic (hence, even) permutations in the group~$S_3$, and then double. This yields the three terms
\begin{equation}\label{EqTriEval}
\Bigl(\tfrac{\delta H_1}{\delta\ba}\circ 
\bigl((A)\overleftarrow{\dd}^{\,(\ba)}_{A(\delta H_3/\delta\ba)}\bigr)
\bigl(\tfrac{\delta H_2}{\delta\ba}\bigl)\Bigr) +
\text{cyclic permutations}.
\end{equation}
Uniting the two parts of the reasoning, we conclude that the left\/-\/hand side~\eqref{EqLHSJac} of Jacobi identity~\eqref{JacPB} for the bracket~$\{\,,\,\}_{\bcP}$ and the value of tri\/-\/vector~$\lshad\bcP,\bcP\rshad$ at the same Hamiltonians~$H_1$,\ $H_2$,\ and~$H_3$ as in~\eqref{JacPB} are equal, hence simultaneously (non)\/trivial, as elements of the cohomology group~$\bar{H}^n_0\bigl(\bpiNCZO 
\bigr)$.%
\hfill$\square$


\smallskip
Referring to Remark
~\ref{FootZeroQZeroValue}
on p.~\pageref{FootZeroQZeroValue} 
and Lemma~\ref{LCriterionStarShaped}
and setting~$\bcQ=\lshad\bcP,\bcP\rshad$ there, 
we conclude that over star\/-\/shaped domains $\subseteq M^n$,
the bracket~$\{\,,\,\}_{\bcP}$ is Poisson if and only if the classical master\/-\/equation 
$\lshad\bcP,\bcP\rshad\cong0$ holds for~$\bcP$.


\subsubsection{Complete integrability}\label{SecMagri}
In the final section we address the 
cohomological structures of (non)\/commutative variational Poisson theory. We recall how the differential $\boldsymbol{\dd}_{\bcP}=\lshad\bcP,{\cdot}\rshad$ specified by a given Poisson bi\/-\/vector~$\bcP$ owes its property $\boldsymbol{\dd}_{\bcP}^2=0$ to a weak variant of the Jacobi identity for the variational Schouten bracket~$\lshad{\cdot},{\cdot}\rshad$. (We remember that the (non)\/commutative variational Poisson formalism is a science of steps and stops, so that calculations involving~$\lshad{\cdot},{\cdot}\rshad$ can be interrupted at every moment, to make legitimate the integrations by parts within every object. This makes the weak variant of Jacobi identity for~$\lshad{\cdot},{\cdot}\rshad$ different from~\eqref{JacobiCommutator} on p.~\pageref{JacobiCommutator}.)

\begin{state}\label{PropWeakJac}
Let $F,G,H\in\bar{H}^n_*\bigl(\bpiNCZO 
\bigr)$ be (non)\/commutative variational multivectors; suppose that $F$~and $G$~are homogeneous. Then the weak variant of Jacobi identity,
\begin{equation}\label{EqJacWeak}
\underline{\lshad F,\underline{\lshad G,H\rshad} }
-(-)^{(|F|-1)\cdot(|G|-1)}\,\underline{\lshad G,\underline{\lshad F,H\rshad} }
\cong \underline{\lshad \underline{\lshad F,G\rshad}, H\rshad },
\end{equation}
holds modulo integrations by parts in every Schouten bracket.

\noindent%
$\bullet$\quad Equivalently, for every homogeneous (non)commutative variational multivector~$Z$ define the shifted (by~$-1$) graded evolutionary vector field~$\bQ^Z$ on the jet space $J^\infty(\bpiNCZO 
)$: by definition, let $\lshad Z,\cH\rshad\cong \vec{\bQ}{}^Z(\cH)$ for all $\cH\in H^n_0(\bpiNCZO 
)$. In these terms, Jacobi identity~\eqref{EqJacWeak} is
\[
\bigl[\vec{\bQ}{}^F,\vec{\bQ}{}^G\bigr]\cong\vec{\bQ}{}^{\lshad F,G\rshad},
\]
that is, the graded commutator of adjoint actions $\lshad F,{\cdot}\rshad$ and~$\lshad G,{\cdot}\rshad$
is equivalent, modulo integrations by parts, to the adjoint action of the object~$\underline{\lshad F,G\rshad}$.
\end{state}

\begin{cor}
By satisfying the master\/-\/equation $\lshad\bcP,\bcP\rshad\cong0$, each (non)\/commutative variational Poisson bi\/-\/vector~$\bcP$ determines the Poisson differential~$\boldsymbol{\dd}_{\bcP}=\lshad\bcP,{\cdot}\rshad$.
\end{cor}

Indeed, Jacobi identity~\eqref{EqJacWeak} then reads $\boldsymbol{\dd}_{\bcP}^2(\cdot) = \underline{ \lshad\bcP,\underline{\lshad\bcP,{\cdot}\rshad}\rshad } \cong
\tfrac{1}{2}\lshad \underline{\lshad\bcP,\bcP\rshad},{\cdot}\rshad =0$.


\begin{proof}[Sketch of the proof \textup{(}of Proposition~\protect\ref{PropWeakJac}\textup{)}]
The graded derivation $\overleftarrow{\bQ}{}^{\!H}\cong\lshad{\cdot},H\rshad$
which acts clockwise (i.e.\ against the orientation) along the cyclic word~$\lshad F,G\rshad$ is permutable with the graded derivations~$\overrightarrow{\bQ}{}^{\!F}$ and~$\overrightarrow{\bQ}{}^{\!G}$ which act counterclockwise on~$G$ and, respectively, on~$-(-)^{(|F|-1)\cdot(|G|-1)} F$ in the object~$\lshad F,G\rshad$. Depending on the  origin --\,from either~$G$ or~$F$\,-- of an argument of~$\overleftarrow{\bQ}{}^{\!H}$ in the right\/-\/hand side of~\eqref{EqJacWeak}, the respective term in that Leibniz rule expansion is realised by using either
\begin{align*}
\bigl(\overrightarrow{\bQ}{}^{\!F}\,(G)\bigr)\,\overleftarrow{\bQ}^{\!H} &=
\overrightarrow{\bQ}{}^{\!F}\,\bigl((G)\,\overleftarrow{\bQ}^{\!H} \bigr)\\
\intertext{or}
-(-)^{(|F|-1)\cdot(|G|-1)} \bigl(\overrightarrow{\bQ}{}^{\!G}\,(F)\bigr)\,\overleftarrow{\bQ}^{\!H} &= -(-)^{(|F|-1)\cdot(|G|-1)} 
\overrightarrow{\bQ}{}^{\!G}\,\bigl((F)\,\overleftarrow{\bQ}^{\!H} \bigr),
\end{align*}
so that all 
terms (and only those terms) in the left\/-\/hand side of~\eqref{EqJacWeak} are recovered.
\end{proof}

For every~$\bcP$, the Poisson differential~$\boldsymbol{\dd}_{\bcP}$ gives rise to the Poisson\/(--\/Lichnerowicz) co\-ho\-mo\-lo\-gy groups~$\mathrm{H}^k_{\bcP}$, $k\ge0$.
\begin{itemize}
\item The group $\mathrm{H}^0_{\bcP}$ is composed by the Casimirs $\cH_0\in\bar H^n(\bpiNCZO 
)$
such that $\lshad\bcP,\cH_0\rshad\cong 0$.
\item The first Poisson cohomology group $\mathrm{H}^1_{\bcP}$ consists of 
cocycle variational one-\/vectors~$X$ without Hamiltonians:
$\lshad\bcP,X\rshad\cong0$ but $X\neq\lshad\bcP,\cH\rshad$ 
for any $\cH\in\overline H^n(\bpiNCZO 
)$.
\item The second group $\mathrm{H}^2_{\bcP}$ contains nontrivial deformations of the Poisson bi\/-\/vector~$\bcP$, \textrm{i.\,e.} those shifts $\bcP\mapsto\bcP+\veps\cdot\bcQ+\ov{o}(\veps)$
infinitesimally preserving the classical master\/-\/equation $\lshad\bcP,\bcP\rshad=0$ which are not generated by the bi\/-\/vector~$\bcP$ itself: $\bcQ\ne\lshad\bcP,{X}\rshad$ for any one\/-\/vector~${X}$.
\item The third group $\mathrm{H}^3_{\bcP}$ contains obstructions to the integrability of infinitesimal shifts 
$\bcP\mapsto\bcP+\veps\cdot\bcQ+\ov{o}(\veps)$ to genuine deformations $\bcP\mapsto\bcP(\veps)$ at~$\veps>0$.
\end{itemize}
These interpretations are standard~\cite{Gerstenhaber}; we also refer to~\cite{f16} for an illustration of classical Poisson deformation theory in the commutative set\/-\/up (in this context, see Problem~\ref{PrbTetraNC} at the end of this chapter).

Likewise, the vanishing of some extra 
cohomological obstructions implies the existence of infinitely many Hamiltonians in involution and the presence of hierarchies of commuting
flows. This is the renowned (Lenard\/--)\/Magri scheme~\cite{Magri1979}.

\begin{theor}[The Magri scheme]\label{ThMagri}
Let $\bcP_1$ and $\bcP_2$ be two \textup{(}non\textup{)}\/commutative variational Poisson bi\/-\/vectors on the jet space $J^{\infty}(\bpiNCZO 
)$.
Suppose they are compatible\textup{:} $\lshad\bcP_1,\bcP_2\rshad\cong0$\textup{,} and assume that the first Poisson\/--\/Lichnerowicz
cohomology group $\mathrm{H}^1_{\bcP_1}$ with respect to the differential $\boldsymbol{\dd}_{\bcP_1}=\lshad\bcP_1,\cdot\,\rshad$ vanishes.
Let $\cH_0\in\mathrm{H}^0_{\bcP_1}\subseteq\bar{H}^n(\bpiNCZO 
)$ be a Casimir of~$\bcP_1$.

Then for any integer $k>0$ there is a Hamiltonian functional $\cH_k\in\bar H^n(\bpiNCZO 
)$ such that
\begin{equation}\label{MagriResolve}
\lshad\bcP_2,\cH_{k-1}\rshad=\lshad\bcP_1,\cH_k\rshad.
\end{equation}
Moreover\textup{,} let $\cH_0^{(\alpha)}$ and $\cH_0^{(\beta)}$ be any two distinct Casimirs 
for the bi\/-\/vector~$\bcP_1$
and construct the two infinite sequences of the functionals $\cH_i^{(\alpha)}$ and $\cH_j^{(\beta)}$
by using~\eqref{MagriResolve}\textup{,} here $i,j\ge0$.
Let $
{\vph_i^{(\alpha)}}\mathrel{{:}{=}}\lshad\bcP_1,\cH_i^{(\alpha)}\rshad$ and similarly\textup{,}
$
{\vph_j^{(\beta)}}\mathrel{{:}{=}}\lshad\bcP_1,\cH_j^{(\beta)}\rshad$. Then for all~$i,j$ and~$\alpha,\beta$\textup{,}
\begin{itemize}
\item
the Hamiltonians~$\cH_i^{(\alpha)}$ and~$\cH^{(\beta)}_j$ Poisson\/-\/commute with respect to either of the Poisson brackets\textup{,}
$\{\,,\,\}_{\bcP_1}$ and $\{\,,\,\}_{\bcP_2}$\textup{;}
\item
the one\/-\/vectors 
$
{\vph_i^{(\alpha)}}$ and~$
{\vph_j^{(\beta)}}$ commute\textup{;}
\item
the density of $\cH_i^{(\alpha)}$ is conserved\textup{,}
$\lshad\cH_i^{(\alpha)},{\vph_j^{(\beta)}}
\rshad\cong 0$\textup{,}
by virtue of each one\/-\/vector
~$
{\vph_j^{(\beta)}}
$.
\end{itemize}
\end{theor}


\begin{proof}[Existence proof.]
Main homological equality~\eqref{MagriResolve} is established by
induction on~$k$. 
Consider the bi\/-\/vectors~$\bcP_1$ and~$\bcP_2$ and a Hamiltonian~$\cH_0$. The steps\/-\/and\/-\/stops variant of Jacobi identity, see~\eqref{EqJacWeak} above, acquires the form
\begin{equation}\label{Jacobi4Bivectors}
\underline{\lshad \bcP_1,\underline{\lshad\bcP_2,\cH_0\rshad} } 
+ \underline{\lshad \bcP_2,\underline{\lshad\bcP_1,\cH_0\rshad} } \cong
\underline{ \lshad \underline{ \lshad\bcP_1,\bcP_2\rshad},\cH_0 \rshad }.
\end{equation}
Hence by starting with a Casimir for a given Poisson bi\/-\/vector~$\bcP_1$, we obtain that
\[
0=[\![\boldsymbol{\cP}_2,0]\!]\cong[\![\boldsymbol{\cP}_2,[\![\boldsymbol{\cP}_1,\mathcal{H}_0]\!]]\!] \cong
 -[\![\boldsymbol{\cP}_1,[\![\boldsymbol{\cP}_2,\mathcal{H}_0]\!]]\!]\mod[\![\boldsymbol{\cP}_1,\boldsymbol{\cP}_2]\!]\cong0,
\]
using Jacobi identity~\eqref{Jacobi4Bivectors}.
The first Poisson cohomology $\mathrm{H}^1_{\bcP_1}
=0$ is trivial by an assumption of the theorem, hence the closed element
$[\![\boldsymbol{\cP}_2,\mathcal{H}_0]\!]$ in the kernel of $[\![\boldsymbol{\cP}_1,\cdot]\!]$
is exact: $[\![\boldsymbol{\cP}_2,\mathcal{H}_0]\!]\cong[\![\boldsymbol{\cP}_1,\mathcal{H}_1]\!]$
for some~$\mathcal{H}_1$.
For~$k\ge1$ we have that
\[
[\![\boldsymbol{\cP}_1,[\![\boldsymbol{\cP}_2,\mathcal{H}_k]\!]]\!] \cong
-[\![\boldsymbol{\cP}_2,[\![\boldsymbol{\cP}_1,\mathcal{H}_k]\!]]\!] \cong
-[\![\boldsymbol{\cP}_2,[\![\boldsymbol{\cP}_2,\mathcal{H}_{k-1}]\!]]\!] \cong0
\]
using~\eqref{Jacobi4Bivectors}
and by $[\![\boldsymbol{\cP}_2,\boldsymbol{\cP}_2]\!]\cong 0$.
Consequently, by~$\mathrm{H}^1_{\cP_1}
=0$ we have that $[\![\boldsymbol{\cP}_2,\mathcal{H}_k]\!]\cong[\![\boldsymbol{\cP}_1,\mathcal{H}_{k+1}]\!]$, 
and we thus proceed indefinitely.
\end{proof}

\begin{define}
Bi\/-\/Hamiltonian evolutionary differential equations which satisfy the hypotheses of Theorem~\ref{ThMagri} and possess as many non\/-\/extendable sequences of local Hamiltonians in involution as the number of the unknowns are called the (infinite\/-\/di\-men\-sio\-nal) \emph{completely integrable systems}.
\end{define}

The (non)\/commutative Korteweg\/--\/de Vries equation~\cite{Magri1979,OlverSokolovCMP1998} is the best\/-\/known example 
of an infinite\/-\/dimensional completely integrable system. 

\begin{rem}\label{RemRestrictHamOp}
The inductive step, that is, the existence of the next, $(k+1)$\/th Hamiltonian functional in involution
with all the preceding ones, 
is possible if and only if the seed~$\cH_0$ is a Casimir,%
\footnote{The Magri scheme starts from any two Hamiltonians
$\mathcal{H}_{k-1},\mathcal{H}_k\in\bar{H}^n(\piNC 
)$ that satisfy~\eqref{MagriResolve},
but we 
operate with the maximal subspaces of the space of functionals
such that the sequence~$\{\mathcal{H}_k\}$ 
cannot be extended with any local quantities at~$k<0$.}
and therefore the Hamiltonian operators~$A_i$ in the bi\/-\/vectors $\bcP_i=\frac12\langle\bb,A_i(\bb)\rangle$
are restricted onto the linear subspace which is spanned in the space of 
variational covectors
by the Euler derivatives of the descendants of~$\mathcal{H}_0$, 
i.\,e.\ of the Hamiltonians of the hierarchy. 
We note that 
the image under~$A_2$ of a generic element 
from the domain of operators~$A_1$ and~$A_2$
cannot be resolved w.r.t.\ $A_1$ by~\eqref{MagriResolve}.

For example, the image $\img A_2^{\KdV}$ of the second Hamiltonian operator for the purely commutative Kor\-te\-weg\/--\/de~Vries equation is not entirely contained in the image of the first structure for the generic values of the arguments. But on the linear subspace of descendants~$\cH_k$ of the Casimir $\int a\,\mathrm{d}x$ for $A_1^{\KdV}$, the inclusion $\img A_2^{\KdV}\subseteq\img A_1^{\KdV}$ is attained.
\end{rem}

\begin{OpenProblem}[The Kontsevich tetrahedral flows]\label{PrbTetraNC}
Does the construction from~\cite{Ascona96,KontsevichFormality} and~\cite{f16} of the quartic\/-\/nonlinear flow $\dot{\cP}=\cQ_{1:\frac{6}{2}}([\cP])$ on spaces of Poisson bi\/-\/vectors~$\cP$ over affine $m$-\/dimensional manifolds~$N^m$ extend --\,in the frames of cyclic word calculus\,-- to the finite\/-\/dimensional%
\footnote{The construction of tetrahedral flow is known~\cite{f16} to have no universal extension to the purely commutative variational set\/-\/up: the flows $\dot{\bcP}=\bcQ_{1:\frac{6}{2}}([\bcP])$ do not always preserve --\,even infinitesimally\,-- the property of Cauchy data~$\bcP$ to be variational Poisson bi\/-\/vectors. Consequently, a search for non\-com\-mu\-ta\-ti\-ve \emph{and} variational generalisation for the existing flow $\dot{\cP}=\cQ_{1:\frac{6}{2}}([\cP])$ is not in order.}
formal noncommutative Poisson geometry\,?

Is such cyclic\/-\/word generalisation also possible for the flow of nonlinearity degree six 
which is built in~\cite{sqs17} from the pentagon\/-\/wheel cocycle in the graph complex\,?
\end{OpenProblem}

\newpage
\subsubsection*{Acknowledgements}
The author is grateful to M.\,Kontsevich for very helpful discussions
and to the 
referee for constructive criticism.
A~part of this research was done while the author was visiting at the $\smash{\text{IH\'ES}}$ (Bures\/-\/sur\/-\/Yvette, France) and MPIM Bonn, Germany; 
financial support and warm hospitality of both institutions are gratefully acknowledged.
This research was also supported in part 
by 
JBI~RUG project~106552 (Groningen, The Netherlands). 




\begin{thebibliography}{77}\normalsize

\bibitem{AKZS}
\by{Alexandrov M., Schwarz A., Zaboronsky O., Kontsevich M.} (1997)
The geometry of the master equation and topological quantum field
theory, \jour{Int.~J.\ Modern Phys.} \vol{A12}:7, 1405--1429.

\bibitem{BV1981}
\by{Batalin I., Vilkovisky G.} (1981)
Gauge algebra and quantization, \jour{Phys.\ Lett.} \vol{B102}:1, 27--31.

\bibitem{BV1983}
\by{Batalin I. A., Vilkovisky G. A.} (1983)
Quantization of gauge theories with linearly dependent generators, 
\jour{Phys.\ Rev.} \vol{D29}:10, 2567--2582. 

\bibitem{Bar2007}
\by{Barannikov S.} (2007)
Modular operads and Batalin\/--\/Vilkovisky geometry,
\jour{Int.\ Math.\ Res.\ Not.} IMRN~\vol{19}, Article~
rnm075, 31~p. 

\bibitem{BarCRM2010}
\by{Barannikov S.} (2010)
Noncommutative Batalin\/--\/Vilkovisky geometry and matrix integrals,
\jour{C.\,R.\,Math.\ Acad.\ Sci.\ Paris} \vol{348}:7--8, 359--362. 

\bibitem{f16}
\by{Bouisaghouane A., Buring R., Kiselev A.} (2017)
The Kontsevich tetrahedral flow revisited,
\jour{J.~Geom.\ Phys.}~\vol{119}, 272--285.\ %
\jour{Preprint} \texttt{arXiv:1608.01710} 
[q-alg]

\bibitem{cpp}
\by{Buring R., Kiselev A.\,V.} (2017) 
The expansion $\star$ mod~$\bar{o}(\hbar^4)$
and computer\/-\/assisted proof schemes in the Kon\-tse\-vich deformation quantization,
50${}+{}$xvi~p.,
\jour{Preprint} $\smash{\text{IH\'ES}}$/M/17/05, 
\texttt{arXiv:1702.00681}~
[math.CO]

\bibitem{JNMP2017}
\by{Buring R., Kiselev A.\,V., Rutten N.\,J.} (2017)
The heptagon\/-\/wheel cocycle in the Kontsevich graph complex,
\jour{J.~Nonlin.\ Math.\ Phys.} \vol{24} 
Suppl.~1, 
157--173.\ 
\texttt{arXiv:1710.00658}~[math.CO]

\bibitem{sqs17}
\by{Buring R., Kiselev A.\,V., Rutten N.\,J.} (2017)
Poisson brackets symmetry from the pentagon\/-\/wheel cocycle in the graph complex,
4$+$v~p.\ \jour{Preprint} \texttt{arXiv:1712.05259} \mbox{[math-ph]}

\bibitem{CattaneoFelderCMP2000}
\by{Cattaneo A. S., Felder G.} (2000) 
A~path integral approach to the Kontsevich quantization formula, \jour{Comm.\ Math.\ Phys.}
\vol{212}:3, 591--611.

\bibitem{ConnesBook}
\by{Connes A.} (1994)
\book{Noncommutative geometry.}
Academic Press Inc., San Diego,~CA. 

\bibitem{ConwaySloane}
\by{Conway J.\,H., Sloane N.\,J.\,A.} (1999)
\book{Sphere packings, lattices and groups.} 3rd ed., 
Grundlehren der Mathematischen Wissenschaften 
\vol{290}, Springer--Verlag, NY.

\bibitem{LoopReview}
\by{Creutz M.} (1983)
\book{Quarks, gluons and lattices.} (Cambridge Univ.\ Press, Cambridge).

\bibitem{DeSoleKacCMP2012}
\by{De Sole A., Kac V. G.} (2012)
Essential variational Poisson cohomology,
\jour{Comm.\ Math.\ Phys.} \vol{313}:3, 837--864. 

\bibitem{DeSoleKacJap2013}
\by{De~Sole A., Kac V.~G.} (2013)
The variational {P}oisson cohomology,
\jour{Jpn.\ J.~Math.} {\bf 8}, 1--145.

\bibitem{DVVWDVV}
\by{Dijkgraaf R., Verlinde H., Verlinde E.} (1991)
Topological strings in $d<1$, \jour{Nuclear Phys.} \vol{B352}:1, 59--86.

\bibitem{DubrZhang2001}
\by{Dubrovin B., Zhang Y.} (2001) 
Normal forms of hierarchies of integrable {P}{D}{E}s, Fro\-be\-ni\-us manifolds and {G}romov\/--\/{W}itten invariants,
\textit{Preprint} \texttt{math.DG/0108160}, 187~pp.

\bibitem{LoopQG}
\by{Gambini R., Pullin J.} (1996)
\book{Loops, knots, gauge theories and quantum gravity.}
(Cambridge Univ.\ Press, Cambridge). 

\bibitem{GelfandShilov}
\by{Gel'fand I. M., Shilov G. E.} (1964) 
\book{Generalized functions.}~\vol{1}.
Properties and operations. 
Academic Press, 
NY\/--\/London.\\ 
\quad
\by{Gel'fand I. M., Shilov G. E.} (1968) 
\book{Generalized functions.}~\vol{2}.
Spaces of fundamental and generalized functions. 
Academic Press, 
NY\/--\/London. 

\bibitem{Gerstenhaber}
\by{Gerstenhaber M., Schack S. D.} (1988) Algebraic cohomology and
deformation theory. \book{Deformation theory of algebras and structures
and applications} (Il~Ciocco, 1986), 
NATO Adv.\ Sci.\ Inst.\ Ser.~C Math.\ Phys.\ Sci.~\vol{247} 
(M.~Gerstenhaber and M.~Hazelwinkel, eds.) Kluwer,
Dordrecht, 11--264.

\bibitem{GitmanTyutin}
\by{Gitman D. M., Tyutin I. V.} (1990)
\book{Quantization of fields with constraints}.
Springer Ser.\ Nucl.\ Part.\ Phys., 
Springer\/-\/Verlag, Berlin. 

\bibitem{GomisParisSamuel}
\by{Gomis J., Par\'\i s J., Samuel S.} (1995)
Antibracket, antifields and gauge\/-\/theory quantization,
\jour{Phys.\ Rep.} \vol{259}:1-2, 1--145.

\bibitem{Looijenga}
\by{Hain R., 
Looijenga E.} (1997)
Mapping class groups and moduli of curves, 
\book{Algebraic Geometry~-- Santa Cruz~1995.} 
\jour{Proc.\ Symp.\ Pure Math.} \vol{62}.2, 97--142.

\bibitem{HenneauxTeitelboim}
\by{Henneaux M., Teitelboim C.} (1992)
\book{Quantization of gauge systems}. 
Princeton University Press, Princeton, NJ.

\bibitem{Wronskians}
\by{Kiselev A. V.} (2005) 
Associative homotopy Lie algebras and Wronskians,
\jour{Fundam.\ Appl.\ Math.} \vol{11}:1, 159--180
(
English transl.:
\jour{J.~Math.\ Sci.} (2007) \vol{141}:1, 1016--1030).\ %
\texttt{arXiv:math.RA/0410185}


\bibitem{TwelveLectures}
\by{Kiselev A. V.} (2012) The twelve lectures in the (non)\/commutative
geometry of 
dif\-fe\-ren\-ti\-al equations, 
\jour{Preprint} $\smash{\text{IH\'ES}}$/M/12/13 (Bures\/-\/sur\/-\/Yvette, 
France), 140~p.

\bibitem{Prague2011}
\by{Kiselev A. V.} (2012)
Homological evolutionary vector fields in Korteweg\/--\/de Vries,
Liouville, Maxwell, and several other models,
\jour{J.~Phys.\textup{:}\ Conf.\ Ser.} \vol{343}, Proc.\ 7th Int.\ workshop QTS-7
`Quantum Theory and Symmetries' (7--13 August 2011, CVUT Prague,
Czech Republic), Paper 012058, 1--20.\ %
\texttt{arXiv:1111.3272} [math-ph]

\bibitem{gvbv}
\by{Kiselev A. V.} (2013) 
The geometry of variations in Batalin\/--\/Vilkovisky formalism, 
\jour{J.~Phys.\textup{:}\ Conf.\ Ser.} \vol{474}, Proc.\ XXI Int.\ conf.\
`Integrable Systems and Quantum Symmetries' (11--16 June 2013, CVUT Prague,
Czech Republic), Paper 012024, 1--51.\ %
\texttt{arXiv:1312.1262} [math-ph]

\bibitem{sqs13}
\by{Kiselev A.\,V.} (2014) The Jacobi identity for graded\/-\/commutative variational Schouten bracket revisited,
\jour{Physics of 
Particles and 
Nuclei Letters} 
\vol{11}:7,
Proc.\ Int.\ workshop SQS'13 `Supersymmetry and Quantum Symmetries'
(29~July -- 3~August 2013, JINR Dubna, Russia), 
950--953.
\ \texttt{arXiv:1312.4140} [math-ph]


\bibitem{dq15}
\by{Kiselev A.\,V.} (2017) The deformation quantization mapping of Poisson\/-\/{} to associative structures in field theory, 
\jour{Banach Center\ Publ.} \vol{113} 
50th Seminar `Sophus Lie' 
(in press), 24~p.\ %
\jour{Preprint} \texttt{arXiv:1705.01777} [q-alg]
%

\bibitem{prg15}
\by{Kiselev A.\,V.} (2016)
The right\/-\/hand side of the Jacobi identity: to be naught or not to be\,?
\jour{J.~Phys.\textup{:}\ Conf.\ Ser.} \vol{670}, Proc.\ XXIII Int.\ conf.\
`Integrable Systems and Quantum Symmetries' (23--27 June 2015, CVUT Prague,
Czech Republic), Paper~012030, 1--17.\ %
\texttt{arXiv:1410.0173} 
[math-ph]



\bibitem{Norway} 
\by{Kiselev A. V., Krutov A. O.} (2014)
Non\/-\/Abelian Lie algebroids over jet spaces,
\jour{J.~Nonlin.\ Math.\ Phys.} \vol{21}:2, 188--213.\ \texttt{arXiv:1305.4598} [math.DG]


\bibitem{Galli10}
\by{Kiselev A. V., van de Leur J. W.} (2011)
Variational Lie algebroids and homological evolutionary vector fields,
\jour{Theor.\ Math.\ Phys.} \vol{167}:3, 772--784.\ 
\texttt{arXiv:1006.4227} [math.DG]

\bibitem{RingersProtaras}
\by{Kiselev A. V., Ringers S.} (2013)
A comparison of definitions for the Schouten bracket on jet spaces,
Proc.\ 6th Int.\ workshop `Group analysis of differential equations and
integrable systems' (18--20 June 2012, Protaras, Cyprus),
127--141.\ \texttt{arXiv:1208.6196} [math.DG]

\bibitem{KontsevichCyclic}
\by{Kontsevich M.} (1993)
Formal (non)\/commutative symplectic geometry, 
\book{The Gel'fand Mathematical Seminars, 1990-1992}
(L.~Corwin, I.~Gelfand, and J.~Le\-pow\-sky, eds),
Birk\-h\"au\-ser, Boston MA, 173--187.

\bibitem{Ascona96}
\by{Kontsevich M.} (1997)
Formality conjecture. 
\book{Deformation theory and symplectic geometry} (Ascona 
1996, D.\,
Sternheimer, J.\,
Rawnsley and S.\,
Gutt, eds), 
Math.\ Phys.\ Stud.~\vol{20}, Kluwer Acad.\ Publ., Dordrecht, 139--156.

\bibitem{KontsevichFormality}
\by{Kontsevich M.} (2003)
Deformation quantization of {P}oisson manifolds,
\jour{Lett.\ Math.\ Phys.}~\vol{66}:3, 157--216.
(\textit{Preprint} \texttt{q-alg/9709040})

\bibitem{YKSDB}
\by{Kosmann\/-\/Schwarzbach Y.} (2004)
Derived brackets,
\jour{Lett.\ Math.\ Phys.} \vol{69}, 61--87. 


\bibitem{LiuZhang2009AdvM}
\by{Liu Si\/-\/Qi, Zhang Y.} (2011)
Jacobi structures of evolutionary partial differential equations,
\jour{Adv.\ Math.} \vol{227}:1, 73--130. 

\bibitem{LiuZhangCMP2013}
\by{Liu Si\/-\/Qi, Zhang Y.} (2013)
Bihamiltonian cohomologies and integrable hierarchies~I: A special case, 
\jour{Comm.\ Math.\ Phys.} \vol{324}:3, 897--935. 

\bibitem{Magri1979}
\by{Magri F.} (1978) A simple model of the integrable equation,
\jour{J.~Math.\ Phys.} \vol{19}:5, 1156--1162.

\bibitem{ManinFmanifolds}
\by{Manin Yu.\ I.} (2005)
$F$-\/manifolds with flat structure and Dubrovin's duality,
\jour{Adv.\ Math.} \vol{198}:1, 5--26.

\bibitem{ManinSchemesQG}
\by{Manin Yu.\ I.} (2012)
\book{Vvedenie v teoriyu shem i kvantovye gruppy} (in Russian),
MCCME, Moscow. 

\bibitem{MumfordRedBook}
\by{Mumford D.} (1999)
\book{The red book of varieties and schemes.} 2nd ed.,
Lect.\ Notes in Math. \vol{1358}, Springer\/--\/Verlag, Berlin.

\bibitem{Olver1993}
\by{Olver P. J.} (1993) \book{Applications of Lie groups to differential
equations}, Grad.\ Texts in Math.\ \vol{107} (2nd ed.), 
Springer\/--\/Verlag, NY. 

\bibitem{OlverSokolovCMP1998}
\by{Olver P. J., Sokolov V. V.} (1998) Integrable evolution
equations on associative algebras, \jour{Comm.\ Math.\ Phys.}
\vol{193}:2, 245--268.\\
\by{Olver P. J., Wang J. P.} (2000)
Classification of integrable one\/-\/component systems on associative algebras,
\jour{Proc.\ London Math.\ Soc.~(3)} \vol{81}:3, 566--586. 

\bibitem{SchwarzCMP1993}
\by{Schwarz A.} (1993)
Geometry of Batalin\/--\/Vilkovisky quantization,
\jour{Comm.\ Math.\ Phys.} \vol{155}:2, 249--260.

\bibitem{Treves2007}
\by{Treves F.} (2007) Noncommutative KdV hierarchy,
\jour{Rev.\ Math.\ Phys.} \vol{19}:7, 677--724;
Errata (2008) \jour{Rev.\ Math.\ Phys.} \vol{20}:1, 117--118. 

\bibitem{Vaintrob}
\by{Vaintrob A. Yu.} (1997)
Lie algebroids and homological vector fields,
\jour{Russ.\ Math.\ Surv.} \vol{52}:2, 428--429. 

\bibitem{VinogradovCspectral}
\by{Vinogradov A. M.} (1984)
The {${\mathcal{C}}$}-\/spectral sequence, {L}agrangian formalism,
              and conservation laws.~{I, II},
\jour{J.~Math.\ Anal.\ Appl.} \vol{100}:1, 1--
129.

\bibitem{VinNary1999}
\by{Vinogradov A., Vinogradov M.} (1998)
On multiple generalizations of Lie algebras and Poisson manifolds.
\book{Secondary calculus and cohomological physics (Moscow, 1997)},
Contemp.\ Math. \vol{219} (AMS, Providence, RI), 273--287. 

\bibitem{Wilson1975}
\by{Wilson K. G.} (1977) Quantum chromodynamics on a lattice. 
\book{New developments in quantum field theory and statistical mechanics.}
Proc.\ Carg\`ese Summer Inst., Carg\`ese, 1976,
edited by M.~L\'evy et al. 
(NATO Adv.\ Study Inst.\ Ser.~B: Physics \textbf{26}, Plenum, New York\/--\/London), 143--172.

\bibitem{WittenWDVV}
\by{Witten E.} (1990)
On the structure of the topological phase of two\/-\/dimensional gravity,
\jour{Nuclear Phys.} \vol{B340}:2-3, 281--332.

\bibitem{WittenAntibracket}
\by{Witten E.} (1990)
A note on the antibracket formalism,
\jour{Modern Phys.\ Lett.} \vol{A5}:7, 487--494.

\bibitem{ZinnJustin1975}
\by{Zinn\/-\/Justin J.} (1975)
Renormalization of gauge theories.
\book{Trends in elementary particle theory}
(Lect.\ Notes in Phys.~\vol{37},
H.~Rollnick and K.~Dietz eds), Springer, Berlin, 2--39.\\
\by{Zinn\/-\/Justin J.} (1976)
\book{M\'ethodes en th\'eorie des champs / 
Methods in field theory.} 
(\'Ecole d'\'Et\'e de Physique Th\'eorique, Session XXVIII, tenue \`a Les Houches, 28 Juillet\,--\,6 Septembre 1975;
R.~Balian and J.~Zinn\/-\/Justin, eds),
North\/-\/Holland Publ.\ Co., 
Amsterdam~etc. 

\end{thebibliography}
\end{document}